\documentclass[12pt, leqno, a4paper]{article}
\usepackage{latexsym}
\usepackage{amscd, amsmath}
\usepackage{amssymb}
\usepackage[all]{xy}
\def\noproof{{\unskip\nobreak\hfill\penalty50\hskip2em\hbox{}%
     \nobreak\hfill$\Box$\parfillskip=0pt%
     \finalhyphendemerits=0\par}}
\def\enddemo{\ifmmode\eqno\Box\else\noproof\vskip0.8truecm\fi}

\newtheorem{theorem}{Theorem}[section]
\newtheorem{definition}[theorem]{Definition}

\newtheorem{corollary}[theorem]{Corollary}

\newtheorem{assumption}[theorem]{Assumption}
\newtheorem{remark}[theorem]{Remark}
\newtheorem{abschnitt}[theorem]{}
\newtheorem{remarks}[theorem]{Remarks}
\newtheorem{lemma}[theorem]{Lemma}
\newtheorem{proposition}[theorem]{Proposition}

\newcommand{\lra}{\longrightarrow}
\newcommand{\hra}{\hookrightarrow}

\DeclareMathOperator{\Spec}{Spec}
\DeclareMathOperator{\Spf}{Spf}
\DeclareMathOperator{\inv}{inv}
\DeclareMathOperator{\id}{id}
\DeclareMathOperator{\Hom}{Hom}
\DeclareMathOperator{\Aut}{Aut}

\DeclareMathOperator{\uEnd}{{\underline{\End}}}
\DeclareMathOperator{\Ext}{Ext}
\DeclareMathOperator{\uExt}{{\underline{\Ext}}}
\DeclareMathOperator{\End}{End}
\DeclareMathOperator{\Ker}{Ker}

\DeclareMathOperator{\Pic}{Pic}
\DeclareMathOperator{\Gal}{G}
\DeclareMathOperator{\Nrd}{Nrd}
\DeclareMathOperator{\Mod}{Mod}
\DeclareMathOperator{\PMod}{PMod}
\DeclareMathOperator{\Vect}{Vect}
\DeclareMathOperator{\PVect}{PVect}
\DeclareMathOperator{\PCoh}{PCoh}

\DeclareMathOperator{\cAVect}{_{\cA}Vect}
\DeclareMathOperator{\VectcA}{Vect_{\cA}}
\DeclareMathOperator{\cAVectr}{_{\cA}Vect^r}
\DeclareMathOperator{\cAVects}{_{\cA}Vect^s}
\DeclareMathOperator{\VectcAr}{Vect_{\cA}^r}

\DeclareMathOperator{\Sh}{Sh}
\DeclareMathOperator{\wSh}{\widehat{\Sh}}
\DeclareMathOperator{\ChDr}{Sh}
\DeclareMathOperator{\wCD}{\widehat{\ChDr}}
\DeclareMathOperator{\Coh}{Coh}
\DeclareMathOperator{\Mor}{Mor}
\DeclareMathOperator{\Cht}{Sht}
\DeclareMathOperator{\Sch}{Sch}

\DeclareMathOperator{\Br}{Br}
\DeclareMathOperator{\Gl}{GL}
\DeclareMathOperator{\PGL}{PGL}
\DeclareMathOperator{\Bad}{\underline{Disc}}
\DeclareMathOperator{\Image}{Im}
\DeclareMathOperator{\rank}{rank}
\DeclareMathOperator{\rankcA}{rank_{\cA}}
\DeclareMathOperator{\Rad}{Rad}
\DeclareMathOperator{\Coker}{Coker}

\DeclareMathOperator{\Frob}{Frob}
\DeclareMathOperator{\Div}{Div}
\DeclareMathOperator{\DivcA}{Div(\cA)}
\DeclareMathOperator{\detcA}{det_{\cA}}
\DeclareMathOperator{\degcA}{deg_{\cA}}
\DeclareMathOperator{\degcB}{deg_{\cB}}
\DeclareMathOperator{\Jac}{Jac}
\DeclareMathOperator{\dvsor}{div}
\DeclareMathOperator{\lcm}{lcm}
\DeclareMathOperator{\Supp}{Supp}

\DeclareMathOperator{\El}{Ell}
\DeclareMathOperator{\SE}{SE}
\DeclareMathOperator{\Inje}{Inj}

\DeclareMathOperator{\tor}{tor}

\DeclareMathOperator{\stab}{stab}

\DeclareMathOperator{\zar}{zar}
\DeclareMathOperator{\opp}{opp}

\DeclareMathOperator{\an}{an}
\DeclareMathOperator{\nr}{nr}
\DeclareMathOperator{\spe}{sp}
\DeclareMathOperator{\zero}{{\it char}}
\DeclareMathOperator{\pole}{{\it pole}}
\DeclareMathOperator{\ord}{ord}

\newcommand{\fP}{{\mathfrak{P}}}

\newcommand{\fM}{{\mathfrak{M}}}
\newcommand{\fA}{{\mathfrak{A}}}
\newcommand{\fp}{{\mathfrak{p}}}

\newcommand{\barX}{{\overline{X}}}
\newcommand{\barD}{{\overline{D}}}

\newcommand{\Fqbar}{{\overline{\mathbb F}_q}}
\newcommand{\kbar}{{\overline{k}}}

\newcommand{\cA}{{\cal A}}
\newcommand{\cB}{{\cal B}}
\newcommand{\cC}{{\cal C}}
\newcommand{\cD}{{\cal D}}
\newcommand{\cE}{{\cal E}}
\newcommand{\cF}{{\cal F}}
\newcommand{\cG}{{\cal G}}
\newcommand{\cH}{{\cal H}}
\newcommand{\cI}{{\cal I}}
\newcommand{\cJ}{{\cal J}}
\newcommand{\cK}{{\cal K}}
\newcommand{\cL}{{\cal L}}
\newcommand{\cM}{{\cal M}}
\newcommand{\cN}{{\cal N}}
\newcommand{\cO}{{\cal O}}
\newcommand{\cP}{{\cal P}}

\newcommand{\cS}{{\cal S}}
\newcommand{\cT}{{\cal T}}
\newcommand{\cV}{{\cal V}}
\newcommand{\cW}{{\cal W}}
\newcommand{\cX}{{\cal X}}
\newcommand{\cY}{{\cal Y}}

\newcommand{\cLd}{{\cal L_{\star}}}
\newcommand{\cMd}{{\cal M_{\star}}}

\newcommand{\cEd}{{\cal E_{\star}}}
\newcommand{\cFd}{{\cal F_{\star}}}

\newcommand{\cId}{{\cal I_{\star}}}
\newcommand{\cJd}{{\cal J_{\star}}}
\newcommand{\cKd}{{\cal K_{\star}}}

\newcommand{\wcO}{{\widehat{\cal O}^{\nr}}}
\newcommand{\wcF}{{\widehat{F}^{\nr}}}

\newcommand{\Cd}{{C_{\star}}}
\newcommand{\Dd}{{D_{\star}}}

\makeatletter
\renewcommand{\subsubsection}{\@startsection{subsubsection}{2}%
        {\z@}{-3.25ex plus -1ex minus-.2ex}{-0em}{\bf}}
\makeatother

\newcommand{\barcA}{{\overline{\cal A}}}
\newcommand{\barcB}{{\overline{\cal B}}}

\newcommand{\barcE}{{\overline{\cal E}}}

\newcommand{\barcL}{{\overline{\cal L}}}

\newcommand{\barcT}{{\overline{\cal T}}}

\newcommand{\tax}{{^{\tau}x}}
\newcommand{\tcE}{{^{\tau}\cal E}}

\newcommand{\tcL}{{^{\tau}\cal L}}
\newcommand{\ta}{{^{\tau}}}
\newcommand{\sM}{{^{\sigma}\cM}}
\newcommand{\sig}{{^{\sigma}}}
\newcommand{\scM}{{^{\sigma}{\cal M}}}
\newcommand{\sncM}{{^{{\sigma^n}}{\cal M}}}
\newcommand{\snmecM}{{^{{\sigma^{n-1}}}{\cal M}}}

\newcommand{\srcM}{{^{{\sigma^r}}{\cal M}}}
\newcommand{\srmecM}{{^{{\sigma^{r-1}}}{\cal M}}}
\newcommand{\srfP}{{^{{\sigma^r}}{\mathfrak P'}}}

\newcommand{\SEl}{{\cS\cE}}
\newcommand{\Ell}{{\cE\ell\ell}}
\newcommand{\PEll}{{\cP\cE\ell\ell}}

\newcommand{\bFq}{{{\mathbb F}_q}}

\newcommand{\bN}{{\mathbb N}}
\newcommand{\bQ}{{\mathbb Q}}

\newcommand{\bR}{{\mathbb R}}
\newcommand{\bZ}{{\mathbb Z}}
\newcommand{\bA}{{\mathbb A}}
\newcommand{\bF}{{\mathbb F}}

\newcommand{\bG}{{\mathbb G}}

\newcommand{\PP}{{\mathbb P}}
\newcommand{\io}{{\iota}}

\newcommand{\ioinf}{{\infty_S}}

\newcommand{\la}{{\lambda}}
\newcommand{\al}{{\alpha}}
\newcommand{\wxi}{{\widetilde{\xi}}}

\newcommand{\seq}{{\subseteq}}
\newcommand{\speq}{{\supseteq}}

\newcommand{\noi}{\noindent}

\begin{document}

\title{Twists of Drinfeld--Stuhler modular varieties}
\author{By Michael Spie{\ss}}
\date{}
\maketitle

\section{Introduction}

In \cite{drinfeld1}, V.\ Drinfeld has introduced the analogues of Shimura varieties for $\Gl_d$ over a global field $F$ of positive characteristic. Following a suggestion of U.\ Stuhler the corresponding varieties for an inner form of $\Gl_d$, i.e.\ the group of invertible elements $A^*$ of a central simple algebra $A$ of dimension $d^2$ over $F$, have been introduced by Laumon, Rapoport and Stuhler in \cite{lrs}. For $d=2$ these are the analogues of Shimura curves. In this paper we show that some of these varieties (for different $A$) are twists of each other. This result can be viewed as a global variant of the Cherednik-Drinfeld Theorem for Shimura curves. 

Let us recall the latter in the simplest case (i.e.\ over $\bQ$ and by neglecting level structure). Let $D$ be an indefinite quaternion algebra over $\bQ$ and $\cD$ a maximal order in $D$. The Shimura curve $S_D$ is the (coarse) moduli space corresponding to the moduli problem
\[
(S \to \Spec \bZ)\mapsto \mbox{abelian surfaces over $S$ with $\cD$-action}.
\]
By fixing an isomorphism $D\otimes \bR \cong M_2(\bR)$ the group of units $\cD^*$ acts on the symmetric space $\cH_{\infty}\colon =\PP^1_{/\bR}-\PP^1(\bR)$ (the upper and lower half plane) through linear transformations. The curve  $S_D\otimes_{\bQ}\bR$ admits the following concrete description 
\begin{equation}
\label{eqn:unic}
S_D\otimes_{\bQ}\bR = \cD^*\backslash \cH_{\infty}.
\end{equation}
If $p$ is a prime number which is ramified in $D$ then there is a similar explicit description over $\bQ_p$. For that let $\overline{D}$ be the definite quaternion algebra over $\bQ$ given by the local data $\overline{D}\otimes \bQ_{\ell} \cong D\otimes \bQ_{\ell}$ for all prime numbers $\ell$ different from $p$ and $\overline{D}\otimes \bQ_p \cong M_2(\bQ_p)$. Let $\overline{\cD}$ denote a maximal $\bZ[\frac{1}{p}]$-order in $\overline{D}$ and denote by $\bQ_p^{\nr}$ the quotient field of the ring of Witt vectors of $W(\overline{\bF}_p)$. The Theorem of Cherednik-Drinfeld asserts that
\begin{equation}
\label{equation:cherednik}
S_D\otimes_{\bQ}\bQ_p = \overline{\cD}^*\backslash (\cH_p\otimes_{\bQ_p} \bQ_p^{\nr}).
\end{equation}
(see \cite{cherednik}, \cite{drinfeld2} or \cite{bc}). Here $\overline{\cD}^*$ acts on $\bQ_p^{\nr}$ via $\gamma\mapsto \Frob_p^{-\ord_p(\Nrd(\gamma))}$ and on the $p$-adic upper half plane $\cH_p \colon =\PP^1_{/\bQ_p}-\PP^1(\bQ_p)$ via linear transformations.

Now let $F$ be a global field of positive characteristic, i.e.\ $F$ is the function field of a smooth proper curve $X$ over a finite field $\bFq$. The analogues of Shimura curves over $F$ are the moduli spaces of $\cA$-elliptic sheaves as introduced in \cite{lrs}. In this paper we generalise this notion slightly by making systematically use of hereditary orders. Let $\infty\in X$ be a fixed closed point. For simplicity we assume in the introduction that $\deg(\infty) =1$. Let $A$ be a central simple $F$-algebra of dimension $d^2$ and let $\cA$ be a locally principal hereditary $\cO_X$-order in $A$. The condition {\it locally principal} means that the radical $\Rad(\cA_x)$ of $\cA_x \colon = \cA\otimes_{\cO_X} \widehat{\cO}_{X,x}$ is a principal ideal for every closed point $x\in X$. There exists a positive integer $e= e_x(\cA)$ such that $\Rad(\cA_x)^e$ is the ideal $\cA_x\varpi_x$ generated by a uniformizer $\varpi_{x}$ of $X$ at $x$. The number $e_x(\cA)$ divides $d$ for all $x$ and is equal to $1$ for almost all $x$. We assume in the following that $e_{\infty}(\cA) = d$. If $A$ is unramified at $\infty$ then this amounts to require that $\cA_{\infty}$ is isomorphic to the subring of matrices in $M_d(\widehat{\cO}_{X,\infty})$ which are upper triangular modulo $\varpi_{\infty}$. 

Roughly, an $\cA$-elliptic sheaf with pole $\infty$ is a locally free $\cA$-module of rank 1 together with a meromorphic $\cA$-linear Frobenius having a simple pole at $\infty$ and a simple zero. The precise definition is as follows. 

\noi {\it An $\cA$-elliptic sheaf over an $\bFq$-scheme $S$ is a pair $E = (\cE, t)$ consisting of a locally free right $\cA\boxtimes \cO_S$-module of rank 1 and an injective homomorphism of $\cA\boxtimes \cO_S$-modules
\[
t: (\id_X\times \Frob_S)^*(\cE\otimes_{\cA} \cA(-\frac{1}{d}\infty)) \lra \cE
\]
such that the cokernel of $t$ is supported on the graph $\Gamma_z\subseteq X \times_{\Spec \bFq} S$ of a morphism $z: S \to X$ (called the zero) and is -- when considered as a sheaf on $S\cong \Gamma_z$ -- a locally free $\cO_S$-module of rank $d$.}

\noi Here $\cA(-\frac{1}{d}\infty)$ denotes the two-sided ideal in $\cA$ given by $\cA(-\frac{1}{d}\infty)_x = \cA_x$ for all $x\ne \infty$ and $\cA(-\frac{1}{d}\infty)_{\infty} = \Rad(\cA_{\infty})$. This definition differs, but, as will be proved in the appendix, is equivalent to the one given in \cite{lrs}\footnote{In \cite{lrs}, the authors work with hereditary orders $\cA$ with $\cA_{\infty}\cong M_d(\widehat{\cO}_{X,\infty})$ and parabolic structures at $\infty$ on $\cE$ instead.}. Unlike in loc.\ cit.\ we do not require the zero $z$ to be disjoint from the pole $\infty$ nor from the closed points which are ramified in $A$. Also we allow $\infty$ to be ramified in $A$. For an arbitrary effective divisor $I$ on $X$ there is the notion of a level-$I$-structure on $E$. The moduli stack of $\cA$-elliptic sheaves with level-$I$-structure and fixed degree $\deg(\cE) = \deg(\cA)$ admits a coarse moduli scheme $\El^{\infty}_{\cA, I}$. It is a fine moduli scheme if $I\ne 0$. We 
 will show (Theorem \ref{theorem:gutesmodell}) that $\El^{\infty}_{\cA,I}$ is a semistable quasiprojective scheme of relative dimension $d-1$ over $X-I$. Previously it had been known (see \cite{lrs}, sects.\ 4, 5 and 6) that $\El^{\infty}_{\cA, I}$ is smooth and quasiprojective over $X'-I$ where $X'$ denotes the complement of set of closed points $x\in X$ with $e_x(\cA)>1$. 

Let $B$ be another central simple $F$-algebra of dimension $d^2$ and assume that there exists a closed point $\fp\in X-\{\infty\}$ such that the local invariants of $B$ are given by $\inv_{\infty}(B) = \inv_{\infty}(A)+\frac{1}{d}$, $\inv_{\fp}(B) = \inv_{\fp}(A)-\frac{1}{d}$ and $\inv_x(B)= \inv_x(A)$ for all $x\ne \infty, \fp$. Let $\cB$ be a locally principal hereditary $\cO_X$-order in $B$ with $e_x(\cB) = e_x(\cA)$ for all $x$. Our main result is that the moduli space $\El^{\fp}_{\cB, I}$ is a twist of $\El^{\infty}_{\cA, I}$. To state this more precisely we assume for simplicity that $\deg(\fp) =1$ and $I=0$ (see \ref{theorem:main}  and \ref{theorem:main2} for the general statement). We have 
\begin{equation}
\label{equation:twist}
\El^{\fp}_{\cB}\quad \cong \quad (\El^{\infty}_{\cA}\otimes_{\bFq} \Fqbar)/<w_{\fp}\otimes \Frob_q>.
\end{equation}
Here $w_{\fp}$ is a certain modular automorphism of $\El^{\infty}_{\cA}$ (in the case $d=2$ it is the analogue of the Atkin-Lehner involution at $p$ for a modular or a Shimura curve). 

We explain briefly our strategy for proving (\ref{equation:twist}). We consider invertible $\cA$-$\cB$-bimodules $\cL$ together with a meromorphic Frobenius $t$ having a simple zero at $\infty$ and simple pole at $\fp$. More precisely, for an $\bFq$-scheme $S$, we consider pairs $L= (\cL, t)$ where $\cL$ is an invertible $\cA\boxtimes \cO_S$-$\cB\boxtimes \cO_S$-bimodule and $t$ is an isomorphism of bimodules
\[
t: (\id_X\times \Frob_S)^*(\cL\otimes_{\cA} \cA(-\frac{1}{d}\fp)) \lra \cL \otimes_{\cA}\cA(-\frac{1}{d}\infty).
\]
These will be called {\it invertible Frobenius bimodules of slope $D = \frac{1}{d}\infty-\frac{1}{d}\fp$} and their moduli space will be denoted by $\SE^D_{\cA, \cB}$. We will show in section \ref{subsection:torsor} that it is a torsor over $\Spec \bFq$ of the finite group of modular automorphisms of $\El^{\infty}_{\cA}$ and compute it explicitly (it is instructive to view $\SE^D_{\cA, \cB}$ as an analogue of the moduli space of supersingular elliptic curves with a fixed ring of endomorphisms). In section \ref{subsection:main} we construct a canonical tensor product $\El^{\infty}_{\cA}\times \SE^D_{\cA, \cB}\to \El^{\fp}_{\cB}, (E,L) \mapsto E\otimes_{\cA} L$. The existence of (\ref{equation:twist}) is then a simple consequence.

As already mentioned, we regard (\ref{equation:twist}) as a global form (in the function field case) of the Cherednik-Drinfeld theorem. In fact an analogue of the uniformization result (\ref{eqn:unic}) for the moduli spaces $\El^{\infty}_{\cA,I}$ has been proved by Blum and Stuhler in \cite{stuhler} (in case where the level $I$ is prime to $\infty$). On the other hand Hausberger has shown in \cite{hausberger} (again under the assumption that $\infty$ does not divide the level $I$) that there is also an analogue of the Cherednik-Drinfeld theorem. We explain in section \ref{subsection:uniformization} that, by using (\ref{equation:twist}), it is possible to deduce one type of uniformization from the other one.

We describe briefly the contents of each section. In part \ref{section:lochereditary} and \ref{section:globhereditary} we discuss hereditary orders in central simple algebras over local fields and global function fields. We show in particular that any hereditary order is Morita equivalent to a (locally) principal hereditary order. This is the reason why it suffices to consider $\cA$-elliptic sheaves for locally principal $\cA$. In section \ref{subsection:special} we introduce the notion of a {\it special $\cA$-module}. If $E=(\cE,t)$ is an $\cA$-elliptic sheaf then $\Coker(t)$ is special. The stack $\Coh_{\cA,\spe}$ of special $\cA$-modules plays a key role in the study of the bad fibers of the characteristic morphism $\zero: \El^{\infty}_{\cA,I} \to X$ in section \ref{subsection:coarsemoduli}. In fact $\Coh_{\cA,\spe}$ is an Artin stack and $\zero$ admits a canonical factorization $\El^{\infty}_{\cA,I} \to \Coh_{\cA,\spe}\to X$. In fact $\Coh_{\cA,\spe}$ is an Artin stack and $\zero$ admits a canonical factorization $\El^{\infty}_{\cA,I} \to \Coh_{\cA,\spe}\to X$. We shall show that the first map is smooth and the second semistable. In sections \ref{subsection:aell} -- \ref{subsection:coarsemoduli} we introduce $\cA$-elliptic sheaves and study their moduli spaces and sections \ref{subsection:fqalg} and \ref{subsection:torsor} are devoted to invertible Frobenius bimodules. In section \ref{subsection:main} we construct the tensor product of an $\cA$-elliptic sheaf (with level-$I$-structure) and an invertible Frobenius bimodules (with level-$I$-structure) and prove our main result (Theorems \ref{theorem:main}  and \ref{theorem:main2}). Finally, in section \ref{subsection:uniformization} we discuss the application to uniformization of $\El^{\infty}_{\cA,I}$ by Drinfeld's symmetric spaces and its coverings. 

{\it Acknowledgements.} I thank E.\ Lau, V.\ Paskunas and T.\ Zink for helpful conversations. Parts of this work were done while the author stayed at the Max-Planck-Institut f{\"u}r Mathematik in Bonn in Winter 2004/05. So I am very grateful to this institution for the generous hospitality.

\paragraph{Notation}

As an orientation for the reader we collect here a few basic notations which are used in the entire work. However most notations listed below will be introduced again somewhere in this work.

For a scheme $S$ we let $|S|$ be the set of closed points of $S$. The category of $S$-schemes is denoted by $\Sch/S$. If $S = \Spec k$ for a field $k$ then we also write $\Sch/k$.

The algebraic closure of a field $k$ is denote by $\kbar$. If $k$ is finite then $k_n\subset \kbar$ denotes the extension of degree $n$ of $k$. 

In chapters \ref{section:globhereditary}, \ref{section:aell} and in \ref{subsection:aelllrs}, $X$ denotes a smooth proper curve over some base field $k$. In chapter \ref{section:globhereditary}, $k$ is an arbitrary perfect field of cohomological dimension 1, whereas in chapter \ref{section:globhereditary} $k$ is the finite field $\bFq$. The function field of $X$ is denoted by $F$. For $Y,Z\in \Sch/k$ we write $X\times Y$ for their product over $k$.

For a closed point $x\in X$ we denote by $k(x)$ its residue field and by $\deg(x)$ the degree $[k(x):k]$. If $S$ is a $k(x)$-scheme, then $x_S$ will denote the morphism $S \to \Spec k(x) \hookrightarrow X$. If $S= \Spec k'$ is a field then we also write $x_{k'}$ instead of $x_{\Spec k'}$.

For a non-zero effective divisor $I$ on $X$, we denote the corresponding closed subscheme of $X$ by $I$ as well. If $\cM$ is a sheaf of $\cO_X$-modules then we use $\cM_I$ for $\cM\otimes_{\cO_X} \cO_I$.  

In chapter \ref{section:aell}, for $S\in \Sch/\bFq$ we denote by $\Frob_S$ its Frobenius endomorphism (over $\bFq$). In the case where $S= \Spec k'$ for some algebraic extension field $k'$ of $\bFq$ we also sometimes write $\Frob_q$ for $\Frob_{\Spec k'}$ and Frobenius in the Galois group $\Gal(k'/\bFq)$. If $S\in \Sch/\bFq$ and $\cE$ is a sheaf of $\cO_{X\times S}$-modules then $\tcE$ denotes the sheaf $(\id_X\times \Frob_S)^*(\cE)$. 

We denote by $\bA$ the Adele ring of $F$ and for a finite set of closed points $T$ of $X$ we let $\bA^T$ denote the Adele ring outside of $T$.

\tableofcontents

\section{Local theory of hereditary orders}
\label{section:lochereditary} 

\subsection{Basic definitions}

Let $X$ be a scheme and $\cA$ a quasi-coherent sheaf of $\cO_X$-algebras. We denote by $_{\cA}\Mod$ (resp.\ $\Mod_{\cA}$) the category of sheaves of left (resp.\ right) $\cA$-modules. Let $\cB$ be another quasi-coherent $\cO_X$-algebra. An $\cA$-$\cB$-bimodule $\cI$ is an $\cO_X$-module with a left $\cA$- and right $\cB$-action which are compatible with the $\cO_X$-action. 

$\cA$ and $\cB$ are said to be {\it (Morita) equivalent} (notation: $\cA\simeq \cB$) if there exists a quasi-coherent $\cA$-$\cB$-bimodule $\cI$ and a quasi-coherent $\cB$-$\cA$-bimodule $\cJ$ such that the following equivalent conditions hold: 

\noi (i) There exists bimodule isomorphisms 
\[
\cI\otimes_{\cB} \cJ \lra \cA, \qquad \cJ\otimes_{\cA} \cI \lra \cB.
\]

\noi (ii) The functors
\[
\cdot\otimes_{\cA}\cI: \Mod_{\cA} \lra \Mod_{\cA},\qquad
\cdot\otimes_{\cB}\cJ: \Mod_{\cB} \lra \Mod_{\cA}
\]
are equivalences of categories and mutually quasi-invers.

In this case $\cI$ and $\cJ$ are called {\it invertible} bimodules and $\cJ$ is called the {\it inverse} of $\cI$. The group of isomorphism classes of invertible $\cA$-$\cA$-bimodules will be denoted by $\Pic(\cA)$.

Now assume that $X$ is a Dedekind scheme that is a one-dimensional connected regular noetherian scheme with function field $K$, i.e.\ $\Spec K \to X$ is the generic point. Let $A$ be a central simple algebra over $K$. An {\it $\cO_X$-order} in $A$ is a sheaf of $\cO_X$-algebras $\cA$ with generic fiber $A$ which is coherent and locally free as an $\cO_X$-module. If $\cB$ is an $\cO_X$-order in another central simple $K$-algebra then it is easy to see that an invertible $\cA$-$\cB$-bimodule is a coherent and locally free $\cO_X$-module.

The $\cO_X$-order $\cA$ in $A$ is called {\it maximal} if for any open affine $U = \Spec R \subseteq X$ the set of sections $\Gamma(U, \cA)$ is a maximal $R$-order in $A$. $\cA$ is called {\it hereditary} if its sections $\Gamma(U, \cA)$ over any open affine $U = \Spec R \subseteq X$ is a hereditary $R$-order in $A$ that is any left ideal in $\Gamma(U, \cA)$ is projective (equivalently any right ideal is projective; compare (\cite{reiner}, (10.7)). Let $\cE$ be a locally free $\cO_X$-module of finite rank which has a left or right $\cA$-action compatible with the $\cO_X$-action. Then the set of sections of $\cE$ over any affine open $U = \Spec R \subseteq X$ are a projective $\Gamma(U, \cA)$-module. 

If $X$ is affine, i.e.\ the spectrum of a Dedekind ring $\cO$ we usually identify $\cA$ with its sections $\Gamma(X, \cA)$. An {\it $\cO$-lattice} is a finitely generated torsionfree (hence projective) $\cO$-module.  A (left or right) {\it $\cA$-lattice} is a (left or right) $\cA$-module which is an $\cO$-lattice. By (\cite{reiner}, (10.7)) $\cA$ is hereditary if and only if every (left or right) $\cA$-lattice is projective.

\subsection{Structure theory}
\label{subsection:localstructure}

Let $\cO$ be a henselian discrete valuation ring with maximal ideal $\fp$ and residue field $k= \cO/\fp$. Let $\varpi\in \fp$ be a fixed prime element. We will recall the structure theory of hereditary $\cO$-orders in central simple $K$-algebras (a reference for what follows is \cite{reiner}, section 39).  Since we are only interested in applications to the case where $\cO$ is the henselisation or completion of a local ring in a global field we will assume for simplicity that $k$ is perfect and of cohomological dimension $\le 1$. 


Let $\cA$ be a hereditary $\cO$-order in a central simple $K$-algebra $A$ of dimension $d^2$. Its Jacobson radical will be denoted by $\fP = \fP_{\cA}$. By (\cite{reiner}, 39.1 and exercise 6 on p.\ 365) $\fP$ is an invertible two-sided ideal and any other two-sided invertible fractional ideal is an integral power of $\fP$. Let $B$ be a central simple $K$-algebra equivalent to $A$ and $\cB$ be a maximal order in $B$. We denote its radical by $\fM = \fP_{\cB}$. Let $I$ be an invertible $A$-$B$-bimodule. Its inverse is $J \colon = \Hom_K(I, K)$. Let $\cI$ be a $\cA$-$\cB$-stable lattice in $I$, i.e.\ $a \cI, \cI b \seq \cI$ for all $a\in \cA, b\in \cB$. Such a lattice exists. In fact if $\cL\seq I$ is any $\cO$-lattice then the $\cO$-module generated by the set $\{ a x b\mid a\in \cA, x\in \cL, b\in \cB\}$ is a $\cA$-$\cB$-stable lattice. There exists a positive integer $t$ -- called the {\it type} of $\cA$ -- such that $\fP^t \cI = \cI\fM$ (see \cite{reiner}, 39.18 (i)). It is also equal to the number of isomorphism classes of indecomposable left (or right) $\cA$-lattices. If $\cM$ is an indecomposable left $\cA$-lattice then $\{\fP^i\cM|\, i=0,1,\ldots, t-1\}$ is a full set of representatives of the set of indecomposable left $\cA$-lattices. For $i\in \bZ$ we set $\cI_i \colon = \fP^{-i}\cI$ and $\cJ_i \colon = \Hom_{\cO}(\cI_{-i}, \cO)$. The sequences $\{\cI_i\mid\,i\in\bZ\}$ and $\{\cJ_i\mid\,i\in\bZ\}$ satisfy the following conditions:
\begin{itemize}
\label{itemize:morher}
\item[(i)]
$\fP \cI_i = \cI_{i-1}, \,\,\, \cI_{i}\fM = \cI_{i-t}, \,\,\, \cJ_i\fP = \cJ_{i-1}, \,\,\, \fM\cJ_{i} = \cJ_{i-t}\quad \mbox{ for all $i\in \bZ$}.$
\item[(ii)]
Let $\cA_i \colon = \{ x\in A\mid \, x \cI_i \seq \cI_i\} = \{ x\in A\mid \, \cJ_{-i}x \seq \cJ_{-i}\}$. Then $\cA_1, \ldots, \cA_t$ are the different maximal orders containing $\cA$ and we have $\cA = \cA_1 \cap \ldots \cap \cA_t$ (note that $\cA_i = \cA_j$ if $i\equiv j \mod t$). The lattice $\cI_i$ is an invertible $\cA_i$-$\cB$-bimodule with inverse $\cJ_{-i}$. Note that $\cA_i = \cA_j$ if $i\equiv j \mod t$.
\item[(iii)]
Let $\barcA \colon = \cA/\fP, \barcB \colon = \cB/\fM$ and let
\[
\barcA^{(i)} \colon = \Image(\barcA\to\End_{\barcB}(\cI_i/\cI_{i-1})) \cong \Image(\barcA\to\End_{\barcB}(\cJ_{-i}/\cJ_{-i-1}))
\]
for $i=1, \ldots, t$. Then, considered as a $\barcA^{(i)}$-$\barcB$-bimodule, $\cI_i/\cI_{i-1}$ is invertible with inverse $\cJ_{-i}/\cJ_{-i-1}$. We have
\[
\barcA \cong \barcA^{(1)} \times\ldots \times \barcA^{(t)}
\]
and $\barcA^{(i)} \cong M_{n_i}(k')$ for $i=1, \ldots, t$. Here $k'$ is the center of $\barcB$ and $n_i= \rank_{\barcB}(\cI_i/\cI_{i-1})$. The numbers $(n_1, \ldots, n_t)$ are called the {\it invariants} of $\cA$. They are well-defined up to cyclic permutation.
\end{itemize}

\begin{definition}
\label{definition:index}
The positive integer $e = e(\cA)$ with $\fP^e = \varpi \cA$ will be called the index of $\cA$. 
\end{definition}

We will see below (Lemma \ref{lemma:hered2}) that $e(\cA)$ does not change under finite {\'e}tale base change. If $h^2$ is the dimension of a central division algebra equivalent to $A$ (thus $h$ is the order of $[A]$ in $\Br(F)$ if $k$ is finite) and $t$ is the type of $\cA$ then $e = d t$. 

Recall (\cite{bushnell}, p.\ 216) that $\cA$ is said to be {\it principal} if every two-sided invertible ideal of $\cA$ is a principal ideal or equivalently if there exists $\Pi\in \fP$ with $\cA \Pi = \Pi \cA = \fP$. For example $\cA$ is principal if it is a maximal $\cO$-order in $A$ or if $e(\cA)=d$. This is a consequence of the following characterization of principal orders.

\begin{lemma} 
\label{lemma:principal}
Let $\cA$ be a hereditary order in a central simple $K$-algebra $A$ of dimension $d^2$. The following conditions are equivalent.

\noi (i) $\cA$ is a principal order.

\noi (ii) If $(n_1, \ldots, n_t)$ are the invariants of $\cA$ then $n_1= \ldots =n_t$. 

\noi (iii) Let $\cM_1, \ldots, \cM_t$ be a full set of representatives of the isomorphism classes of indecomposable right $\cA$-lattices. Then there exists an integer $f\in \bN$ such that
\[
\cA \quad \cong \quad (\cM_1\oplus \ldots\oplus \cM_t)^f
\]
as right $\cA$-modules. In this case we have $f = n_1= \ldots =n_t$ and $d = e f$.
\end{lemma}

{\em Proof.} (i) $\Leftrightarrow$ (ii) see (\cite{bushnell}, Theorem 1.3.2, p.\ 217).

\noi (i) $\Leftrightarrow$ (iii) Since $\cA$ is principal if and only if $\cA \cong \fP$ as right $\cA$-modules this follows from the fact that the map $[\cM]\mapsto [\cM\fP]$ is a cyclic permutation of the set isomorphism classes of indecomposable right $\cA$-lattices (\cite{reiner}, 39.23).

For the last assertion note that if $\cA$ is principal then on the one hand
\[
\dim_k(\cA/\fP) = \sum_{j=1}^t \dim_k(\barcA^{(j)}) = \sum_{j=1}^t \, d n_j^2 = t d n_i^2 = e n_i^2
\]
for $i \in \{1,\ldots, t\}$. On the other hand since $\cM_j/\cM_j\fP$ is an irreducible $\barcA^{(j)}$-module we have
\[
\dim_k(\cA/\fP) = f \sum_{j=1}^t \, \dim_k(\cM_j/\cM_j\fP) = f t d n_i = f e n_i\]
Therefore we get $f = n_i$. Finally because of 
\[
d^2 = \dim_k(\cA/\varpi\cA) = \sum_{i=0}^{e-1} \, \dim_k(\fP^i/\fP^{i+1}) = e\dim_k(\cA/\fP)
\]
we obtain $e f =d$.\enddemo

Suppose that $\cA$ is principal. We denote the subgroup of $A^*$ of elements $x\in A^*$ with $x\cA = \cA x$ by $N(\cA)$. For $x\in N(\cA)$ there exists a unique $m\in \bZ$ with $x\cA = \fP^m$ and we set $v_{\cA}(x) = \frac{m}{e}$. We have a commutative diagram with exact rows
\[
\begin{CD}
1 @>>> \cO^* @>>> K^* @> v_K >> \bZ @>>> 0\\
@.@VVV@VVV@VVV\\
1 @>>> \cA^* @>>> N(\cA) @> v_{\cA} >> \frac{1}{e}\bZ @>>> 0
\end{CD}
\]
where $v_K$ denoted the normalized valuation of $K$ and the vertical maps are the natural inclusions.

Next we consider the special case where $A = \End_K(V)$ for a finite-dimensional $K$-vector space $V$ (i.e.\ $A$ is split). A {\it lattice chain} in $V$ is a sequence of $\cO$-lattices $\cLd=\{\cL_i\mid i \in \bZ\}$ such that 
\begin{itemize}
\item[--] $\cL_i \subsetneq \cL_{i+1}$ for all $i\in \bZ$.

\item[--] There exists a positive integer $e$, the {\it period} of $\cLd$, such that $\cL_{i-e} = \varpi \cL_i$ for all $i\in \bZ$.
\end{itemize}
The ring 
\begin{equation}
\label{eqn:herlat}
\cA = \End(\cLd) \colon = \{f\in A\mid\, f(\cL_i) \seq \cL_i \,\,\forall i\in \bZ\}
\end{equation}
is a hereditary $\cO$-order in $A$ of index (= type) $e$ with invariants $n_i = \dim_k(L_i/L_{i-1})$. We have:
\begin{equation}
\label{eqn:radlat}
\fP_{\cA}^{-m} = \End^m(\cLd) \colon =\{f\in A\mid\, f(\cL_i) \seq \cL_{i+m} \,\,\forall i\in \bZ\}.
\end{equation}
Any hereditary $\cO$-order in $A$ is of the form (\ref{eqn:herlat}) for some lattice chain. 

\subsection{\'Etale base change}

We keep the notation and assumption of the last section. Let $A$ be a central simple algebra and $\cA$ an $\cO$-order in $A$ with radical $\fP$.

\begin{lemma} 
\label{lemma:hered1}
The following conditions are equivalent:

\noi (i) $\cA$ is hereditary.

\noi (ii) There exists a two-sided invertible ideal $\fM$ in $\cA$ such that $\cA/\fM$ is semisimple and $\fM^e = \varpi \cA$ for some $e\ge 1$. 

\noi Moreover if $\fM$ is as in (ii) then $\fM = \fP$.
\end{lemma}

{\em Proof.} (i) $\Rightarrow$ (ii) follows from (\cite{reiner}, (39.18) (iii) ) (for $\fM = \fP$).

\noi (ii) $\Rightarrow$ (i) In view of (\cite{reiner}, (39.1)) it suffices to show that $\fM = \fP$. The inclusion $\fM \speq \fP$ is a consequence of the assumption that $\cA/\fM$ is semisimple. The converse inclusion follows from (\cite{reiner}, exercise 1).\enddemo

\begin{lemma} 
\label{lemma:hered2}
Let $K'/K$ be a finite unramified extension and $\cO'$ the integral closure of $\cO$ in $K'$. Then $\cA$ is hereditary (resp.\ principal) if only if $\cA\otimes_{\cO} \cO'$ is hereditary (resp.\ principal). In this case $\fP\otimes_{\cO} \cO'$ is the radical of the latter.
\end{lemma}

{\em Proof.} (compare also \cite{janusz}) We will prove only the statement for hereditary orders and leave the case of principal orders to the reader. If $\cA$ is hereditary then $\fM \colon = \fP\otimes _{\cO} \cO'$ satisfies the condition (ii) of Lemma \ref{lemma:hered1}. Hence $\cA\otimes_{\cO} \cO'$ is hereditary. 

To prove the converse let $\cP$ be a left $\cA$-lattice. We have to show that 
\[
\Hom_{\cA}(\cP, \cdot): \Mod_{\cA}\to \Mod_{\cO}
\]
is an exact functor or -- since $\cO'$ is a faithfully flat $\cO$-algebra -- that
\[
\Hom_{\cA}(\cP, \cdot)\otimes_{\cO} \cO' \cong \Hom_{\cA\otimes_{\cO} \cO'}(\cP\otimes_{\cO} \cO', \cdot\otimes_{\cO} \cO')
\]
is exact. However the assumption implies that $\cP\otimes_{\cO} \cO'$ is a projective $\cA\otimes_{\cO} \cO'$-module.\enddemo

\subsection{Morita equivalence}

Let $A$ be a central simple algebra and $\cA$ a hereditary $\cO$-order in $A$ with radical $\fP$. If $\cA'$ is another $\cO$-order in $A$ containing $\cA$ then $\cA'$ is hereditary as well and $\fP_{\cA'}\seq \fP$. 

\begin{lemma} 
\label{lemma:hererad}
Let $\cA_1,\ldots, \cA_s$ be a collection of $\cO$-orders in $A$ containing $\cA$ with radicals $\fP_1, \ldots, \fP_s$. If $\cA_1\cap\ldots\cap \cA_s = \cA$ then $\fP_1+ \ldots+ \fP_s= \fP$.
\end{lemma} 

{\em Proof.} Clearly $\fP_1+ \ldots+ \fP_s\seq \fP$. By Lemma \ref{lemma:hered2} to prove equality we may pass to a finite unramified extension $K'/K$. Hence we can assume $A = \End_K(V)$ for some finite-dimensional $K$ vector space $V$ and that there exists a lattice chain $\cLd=\{\cL_i\mid i \in \bZ\}$ in $V$ with period $e = e(\cA)$ such that  
\[
\cA = \{f\in A\mid\, f(\cL_i) \seq \cL_i \,\,\forall i\in \bZ\},\quad \fP = \{f\in A\mid\, f(\cL_i) \seq \cL_{i-1} \,\,\forall i\in \bZ\}.
\]
Clearly it is enough to consider the case where $s=e$ and $\cA_i = \{f\in A\mid\, f(\cL_i) \seq \cL_i\}$, $i=1, \ldots, e$ are the different maximal orders containing $\cA$. We proceed by induction on $e$ so we can assume that $e>1$ and that the radical 
\[
\fP' = \{f\in A\mid f(\cL_i) \seq \cL_{i-1} \,\,\forall\, i\not\equiv 0, 1 \,\mbox{mod}\, e\,\mbox{ and } f(\cL_i) \seq \cL_{i-2} \,\,\forall\, i\equiv 1 \,\mbox{mod}\, e\},
\]
of $\cB \colon = \cA_1\cap \ldots\cap\cA_{e-1}$ is $= \fP_1 + \ldots + \fP_{e-1}$. 

Let $f\in \fP$. Consider the diagram of $k$-vector spaces
\[
  \xymatrix{
    \cL_1/\cL_0 \ar[r]^{\overline{f}}\ar@{^{(}->}[d] &  \cL_0/\cL_{-1}\\
    \cL_e/\cL_0 \ar@{-->}[r]^{\overline{g}} & \cL_0/\cL_{-e}\ar@{->>}[u]
}
\]
where the vertical maps are induced by $\cL_1 \hookrightarrow \cL_e$ and $id: \cL_0 \to \cL_0$ respectively and the upper horizontal map by $f$. There exists a dotted arrow $\overline{g}$ making the diagram commutative. Let $g\in \Hom_{\cO}(\cL_e, \cL_0) = \fP_e$ be a ``lift'' of $\overline{g}$. Then $g(x) \equiv f(x) \mod \cL_{-1}$ for every $x\in \cL_1$. Therefore $(f-g)(\cL_1) \seq \cL_{-1}$ and $(f-g)(\cL_{i})\seq f(\cL_i) + g(\cL_e) \seq \cL_{i-1} +\cL_0 = \cL_{i-1}$ for $i = 2, \ldots, e-1$ and consequently $f-g\in \fP'$. This proves $\fP \seq \fP' + \fP_e = \fP_1 + \ldots + \fP_{e-1} + \fP_e$.\enddemo

\begin{corollary} 
\label{corollary:hererad}
Let $\cA$ be a hereditary $\cO$-order of type $t$ with radical $\fP$ in the central simple $K$-algebra $A$ and let $\cA_1, \ldots, \cA_t$ denote the different maximal orders containing $\cA$. Then
\[
\cA_1 + \ldots + \cA_t = \fP^{-t+1}
\]
is a two-sided invertible ideal.
\end{corollary} 

{\em Proof.} By (\cite{reiner}, section 39, exercise 10) and Lemma \ref{lemma:hererad} above we have 
\[
\fP = \fP_1 + \ldots + \fP_t = \fP^t\cA_1 + \ldots + \fP^t\cA_t = \fP^t(\cA_1 + \ldots + \cA_t)
\]
hence $\cA_1 + \ldots + \cA_t = \fP^{-t+1}$.\enddemo

Let $B$ be another central simple $K$-algebra which is equivalent to $A$ and let $\cB$ be a maximal order in $B$ with radical $\fM$. Let $I$ be an invertible $A$-$B$-bimodule, $J \colon = \Hom_K(I, K)$ and let $\{\cI_i\mid\,i\in\bZ\}$ and $\{\cJ_i\mid\,i\in\bZ\}$ be as in section \ref{subsection:localstructure}.

\begin{lemma} 
\label{lemma:heremax}
Consider $\cI_i\otimes_{\cB} \cJ_j$ (resp.\ $\cJ_i\otimes_{\cA} \cI_j$) as a submodule of $(\cI_i\otimes_{\cB} \cJ_j)\otimes_{\cO} K = I \otimes_B J$ (resp.\ $J \otimes_A I$).

\noi (a) $\sum_{i+j = -t +1} \, \cI_i\otimes_{\cB} \cJ_j \cong \cA$ as an $\cA$-$\cA$-bimodule.

\noi (b) $\cJ_i\otimes_{\cA} \cI_j$ is an invertible bimodule. If $i +j=0$ then $\cJ_i\otimes_{\cA} \cI_j \cong \cB$ (as $\cB$-$\cB$-bimodule). We have
\[
\cJ_{i+1}\otimes_{\cA} \cI_j = \cJ_i\otimes_{\cA} \cI_{j+1} = \left\{ \begin{array}{lll}
                  \fM^{-1}(\cJ_i\otimes_{\cA} \cI_j) & \mbox{if $i+j \equiv 0 \mod t$;}\\
                  \cJ_i\otimes_{\cA} \cI_j & \mbox{if $i+j \not\equiv 0 \mod t$}
                \end{array}
\right.
\]
\end{lemma}

{\em Proof.} (a) Under the identification $I \otimes_B J= \Hom_K(J, K)\otimes_B J = \Hom_B(J, J)$ the submodule $\cI_i\otimes_{\cB} \cJ_j$ corresponds to $\Hom_{\cB}(\cJ_{-i}, \cJ_j)$. Hence if we fix an $A$-$A$-bimodule isomorphism $I \otimes_B J \cong A$ so that $\Hom_{\cB}(\cJ_{0}, \cJ_0)$ is mapped to $\cA_0$ then for arbitrary $i,j\in \bZ$ with $i +j = 0$ the module $\Hom_{\cB}(\cJ_{-i}, \cJ_j)$ is mapped to $\cA_i$. It follows $\sum_{i+j = 0} \, \cI_i\otimes_{\cB} \cJ_j \cong \cA_1 + \ldots + \cA_r$ hence together with Lemma \ref{lemma:hererad} the assertion.

\noi (b) The proof of the first two statements is similar and will be left to the reader. For the last statement note that
\[
\Coker(\cJ_i\otimes_{\cA} \cI_j\to \cJ_{i+1}\otimes_{\cA} \cI_j) \cong \cJ_{i+1}/\cJ_i\otimes_{\cA} \cI_j\cong \cJ_{i+1}/\cJ_i\otimes_{\barcA} \cI_j/\cI_{j-1}
\]
By (iii) above we have 
\[
\cJ_{i+1}/\cJ_i\otimes_{\barcA} \cI_j/\cI_{j-1} \cong \cJ_{i+1}/\cJ_i\otimes_{\barcA^{(j)}} \cI_j/\cI_{j-1} \cong \barcB
\]
if $i+j \equiv 0 \mod t$ and $\cJ_{i+1}/\cJ_i\otimes_{\barcA} \cI_j/\cI_{j-1} = 0$ if $i+j \not\equiv 0 \mod t$.\enddemo

\begin{corollary} 
\label{corollary:hermax}
The assignment
\[
\cM \mapsto \{\cM\otimes_{\cA}\cI_i\mid\,i\in\bZ\}
\]
defines an equivalence between the category of right $\cA$-lattices and the category of increasing chains $\{\cM_i\mid\,i\in\bZ\}$ of right $\cB$-lattices such that $\cM_i \fM = \cM_{i-t}$ for all $i\in \bZ$. A quasi-inverse is given by
\[
\{\cM_i\mid\,i\in\bZ\} \mapsto\sum_{i+j = -t +1} \, \cM_i\otimes_{\cB} \cJ_j.
\]
Here the sum is taken inside of $(\bigcup_{i\in \bZ} \cM_i)\otimes_B J$.
\end{corollary}
 
\begin{proposition}
\label{proposition:hermor}
Let $\cA_1$ and $\cA_2$ be hereditary $\cO$-orders in central simple $K$-algebras $A_1$ and $A_2$. The following conditions are equivalent:

\noi (i) $\cA_1$ and $\cA_2$ are Morita equivalent.

\noi (ii) $A_1$ and $A_2$ are equivalent and $\cA_1$ and $\cA_2$ have the same index.
\end{proposition}

{\em Proof.} We will show only that (ii) implies (i). The proof of the converse is easier and will be left to the reader. Suppose that (ii) holds. Let $D$ be a central division algebra over $K$ equivalent to $A_1$ and $A_2$ and $\cD$ be the maximal $\cO$-order in $D$. For $\nu = 1,2$ we fix increasing sequences of $\cA_{\nu}$-$\cD$- and $\cD$-$\cA_{\nu}$-bimodules $\{\cI^{(\nu)}_i\mid i\in \bZ\}$ and $\{\cI^{(\nu)}_i\mid i\in \bZ\}$ as in \ref{subsection:localstructure}. Put $I^{(\nu)} = \bigcup_{i\in \bZ}\cI^{(\nu)}_i$ and $J^{(\nu)} = \bigcup_{i\in \bZ}\cJ^{(\nu)}_i$. The assumption implies that $\cX\colon = \sum_{i+j = -t +1} \, \cI^{(1)}_i\otimes_{\cD} \cJ^{(2)}_j$ is an $\cA_1$-$\cA_2$-lattice (the summation takes place in $I^{(1)}\otimes_D J^{(2)}$) and $\cY\colon = \sum_{i+j = -t +1} \, \cI^{(2)}_i\otimes_{\cD} \cJ^{(1)}_j$ a $\cA_2$-$\cA_1$-lattice. By Corollary \ref{corollary:hermax} above the assignment $\cM \mapsto \cM\otimes_{\cA_1} \cX$ defines an equivalence between the category of right $\cA_1$-lattices and the category of right $\cA_2$-lattices. A quasi-inverse is given by $\cN \mapsto \cN\otimes_{\cA_2} \cY$. This implies that $\cA_1$ and $\cA_2$ are Morita equivalent. In fact using Lemma \ref{lemma:heremax} it is easy to see that the $\cX\otimes_{\cA_2} \cY \cong \cA_1$ and  $\cY\otimes_{\cA_1} \cX \cong \cA_2$.\enddemo

Recall that a right $\cA$-lattice $\cM$ is called {\it stably free} if there exists integers $r\ge 1, s\ge 0$ such that $\cM^r \cong \cA^s$ as right $\cA$-modules. 

\begin{lemma} 
\label{lemma:stablyfree1}
Let $\cA$ is a principal $\cO$-order of index $e$ in a central simple $K$-algebra of dimension $d^2$. Let $\cM_1, \ldots, \cM_t$ be representatives of isomorphism classes of indecomposable right $\cA$-lattices. For a right $\cA$-lattice $\cM\ne 0$ the following conditions are equivalent.

\noi (i) $\cM$ is stably free.

\noi (ii) $\cM\cong (\cM_1\oplus \ldots\oplus \cM_t)^r$ for some positive integer $r$.

\noi (iii) $\cD\colon = \End_{\cA}(\cM)$ is a principal $\cO$-order of index $e$ in a central simple $K$-algebra $D$. 

\noi Moreover in this case $\cA$ and $\cD$ are Morita equivalent and $\cM$ is an invertible $\cD$-$\cA$-bimodule. If $\rank_{\cO} \cM = r d e$ then $\dim_K(D) = (er)^2$.
\end{lemma}

{\em Proof.} The equivalence of (i) and (ii) follows immediately from Lemma \ref{lemma:principal}.

(ii) $\Leftrightarrow$ (iii) By Lemma \ref{lemma:hered2} we may pass to a finite unramified extension $K'/K$. Therefore we can assume that $A = \End_K(V)$ for a $d$-dimensional $K$-vector space $V$ and $\cA = \End(\cLd)$ for a lattice chain $\cLd$ with period $e$ in $V$. There exists $r_1, \ldots, r_e\ge 0$ with 
\[
\cM\cong \cL_1^{r_1}\oplus \ldots\oplus \cL_e^{r_1}
\]
Since $\Hom_{\cA}(\cL_i, \cL_j) \cong \fp^{\mu}$ with $i-j \le \mu e < i-j + e$ we have
\[
\End_{\cA}(\cM)\cong 
\left( \begin{matrix} M_{r_1,r_1}(\cO) & M_{r_1,r_2}(\cO) & \ldots & M_{r_1,r_e}(\cO)\\
M_{r_2,r_1}(\fp) & M_{r_2,r_2}(\cO) & \ldots & M_{r_2,r_e}(\cO)\\
\vdots &\vdots &\ddots &\vdots &\\
M_{r_e,r_1}(\fp) & M_{r_e,r_2}(\fp) & \ldots & M_{r_e,r_e}(\cO)\\
\end{matrix}\right)
\]
By (\cite{reiner}, 39.14) the order on the right is a hereditary order in $M_m(K)$ where $m = \sum_{i=1}^e r_i$. Its index is $=e$ if and only if $r_i\ge 1$ for all $i\in \{1, \ldots, e\}$ and in this case the invariants are $(r_1, \ldots, r_e)$. The equivalence of (ii) and (iii) follows. The proof of the last assertion will be left to the reader.\enddemo

\begin{corollary} 
\label{corollary:stablyfree2}
Let $\cA$ be as in \ref{lemma:stablyfree1} and let $\cM$ be a stably free $\cA$-module. We have:

\noi (a) $\rank_{\cO} \cM$ is a multiple of $e d$. 

\noi (b) $\cM$ is free if and only if $\rank_{\cO} \cM$ is a multiple of $d^2$. In particular if $e = d$ then $\cM$ is free.
\end{corollary}

{\em Proof.} If $A \cong M_m(D)$ where $D$ is the central division algebra equivalent to $A$ then $\rank_{\cO} \cM_i = m h^2$ with $h^2 = \dim_K(D)$. Hence if $\cM\cong (\cM_1\oplus \ldots\oplus \cM_t)^r$ for $r\in \bN$ then $\rank_{\cO} \cM =  r t m h^2 = r e d$. The second assertion is obvious.\enddemo

\begin{corollary}
\label{corollary:invlocfree} Let $\cA$ and $\cB$ be principal orders of the same index $e$ in central simple $K$-algebras $A$ and $B$ respectively both of dimension $d^2$. Let $\cI$ be an $\cA$-$\cB$-bimodule. The following conditions are equivalent:

\noi (i) $\cI$ is an invertible $\cA$-$\cB$-bimodule.

\noi (ii) $\cI$ is a free left $\cA$-module of rank 1.

\noi (iii) $\cI$ is a free right $\cB$-module of rank 1.
\end{corollary}

{\em Proof.} (i) $\Rightarrow$ (iii) We show first that $\cI$ is a lattice. Let $\cJ$ be an inverse of $\cI$ and $_{\cB-\tor}\cJ$ its $\cB$-torsion ($\cA$-)submodule. Since $_{\cB-\tor}\cJ\otimes_{\cA}\cI\hookrightarrow \cJ\otimes_{\cA}\cI \cong \cB$ we have $_{\cB-\tor}\cJ\otimes_{\cA}\cI = 0$ and therefore $_{\cB-\tor}\cJ= {_{\cB-\tor}\cJ} \otimes_{\cA}\cI\otimes_{\cB}\cJ =0$. For $m\in \cJ, m\ne0$ we get $\cB m \cong \cB$ as left $\cB$-module and therefore
\[
\cI \cong \cI\otimes_{\cB}\cB m\hookrightarrow \cI\otimes_{\cB}\cJ \cong \cA.
\]
Hence $\cI$ is a lattice. Let $\cD \colon = \End_{\cB}(\cI)\supseteq \cA$. Thus $\cI$ is a $\cD$-$\cB$-bimodule and so $\cI\otimes_{\cB}\cJ\cong \cA$ is a $\cD$-$\cA$-bimodule. But $\End_{\cA}(\cA_{\cA}) = \cA$\footnote{Here $\cA_\cA = \cA$ considered as a right $\cA$-module} and therefore $\cD = \cA$. By \ref{lemma:stablyfree1} and \ref{corollary:stablyfree2}, $\cI$ is a free $\cB$-module of rank 1.

(iii) $\Rightarrow$ (i) By \ref{lemma:stablyfree1}, $\cD$ is a principal order of index $e$ in a central simple $K$-algebra $D$ of dimension $d^2$ and $\cI$ is an invertible $\cD$-$\cB$-bimodule. Since $\cD \supseteq \cA$ this implies $\cD = \cA$.\enddemo

\begin{corollary}
\label{corollary:principal}
Let $\cA$ be a hereditary $\cO$-order of index $e$ in central simple $K$-algebra $A$. Then there exists a principal $\cO$-order $\cD$ in a central simple $K$-algebra which is Morita equivalent to $\cA$. In fact that $\cD$ can be chosen such that $\rank_{\cO}(\cD) = e^2$.
\end{corollary}

{\em Proof.} Let $\cA'$ be a principal $\cO$-order of index $e$ in $A' \colon = M_e(A)$ (since $e^2$ divides $\dim_K(A')$ and $e$ is a multiple of the order of $[A'] = [A]$ in $\Br(F)$ such an order clearly exists). By Proposition \ref{proposition:hermor}, $\cA'$ is Morita equivalent to $\cA$. The second assertion follows immediately from \ref{lemma:stablyfree1}.\enddemo

\subsection{Maximal tori}

Let $A$ be a central simple $K$-algebra of dimension $d^2$ and $\cA$ a hereditary $\cO$-order in $A$ with radical $\fP$. In this section we consider commutative {\'e}tale $\cO$-subalgebras of $\cA$. Note that a finite flat $\cO$-algebra $\cT$ is {\'e}tale if and only if $\Rad(\cT) = \varpi \cT$.

\begin{lemma} 
\label{lemma:etalerad}
Let $\cT$ be a commutative {\'e}tale $\cO$-subalgebra of $\cA$. Then we have $\Rad(\cT) = \cT \cap \fP$. 
\end{lemma}

{\em Proof.} Since $\cT$ is a direct product of local $\cO$-algebras $\cT = \prod \cT_i$ and $\Rad(\cT) =  \prod \Rad(\cT_i)$ it suffices to prove the assertion for each factor. Thus we may assume that $\cT$ is a local ring. Hence $\Rad(\cT)$ is the maximal ideal of $\cT$ which implies $\cT \cap \fP\subseteq \Rad(\cT)$. On the other hand, by assumption, we have $\Rad(\cT)= \varpi \cT$ hence $\Rad(\cT)\subseteq\cT \cap \fP$.\enddemo

A commutative {\'e}tale $\cO$-subalgebra $\cT$ of $\cA$ is called {\it maximal torus} if $\rank_{\cO} \cT =d$. It follows immediately from the structure theory for hereditary $\cO$-orders in central simple $K$-algebras (\cite{reiner}, 39.14) that there exists a maximal torus in $\cA$. We have the following characterization of maximal tori:

\begin{lemma} 
\label{lemma:maxtoroidal1}
Let $\cT$ be a commutative {\'e}tale $\cO$-subalgebra of $\cA$. The following conditions are equivalent. 

\noi (i) $\cT$ is a maximal torus.

\noi (ii) $\cT$ is a maximal element in the set of commutative {\'e}tale $\cO$-subalgebras of $\cA$.

\noi (iii)  $\cT = Z_{\cA}(\cT) = \{x\in \cA\mid\, x t = tx \quad\forall t\in\cT\}$.

\noi (iv) $\cT/\Rad(\cT)$ is a maximal commutative separable $k$-subalgebra of $\cA/\fP$.
\end{lemma}

{\em Proof.} The simple proof of the equivalence of the first three conditions will be left to the reader. 

(i) $\Leftrightarrow$ (iv) By \ref{lemma:etalerad} above we have $\Rad(\cT) = \cT \cap \fP = \varpi \cT$. Thus it follows from Lemma \ref{lemma:maxcommsepsubalgebra1} of the appendix that (iv) holds if and only if $\rank_{\cO} \cT= \dim_k( \cT/\varpi \cT) =n$.\enddemo

\begin{lemma} 
\label{lemma:maxtoroidal2}
(a) If $k = \bFq$ and $\cA$ is a maximal order in $A$ then $\cA$ admits a maximal torus isomorphic to $\cO_d$, the ring of integers of the unramified extension of degree $d$ of $K$.

\noi (b) Let $\cO'$ be a finite {\'e}tale local $\cO$-algebra and $\cT$ be a maximal torus in $\cA$. Then $\cT\otimes_{\cO} \cO'$ is a maximal torus in $\cA\otimes_{\cO} \cO'$.
 
\noi (c) For any two maximal tori $\cT, \cT'$ of $\cA$ there exists a finite {\'e}tale local $\cO$-algebra $\cO'$ such that $\cT\otimes_{\cO} \cO'$ and $\cT'\otimes_{\cO} \cO'$ are conjugated (by some $a\in (\cA\otimes_{\cO} \cO')^*$).
\end{lemma}

{\em Proof.} (a) and (b) are obvious.

To prove (c) we may pass to a finite unramified extension of $K$ if necessary so that $A = \End_K(V)$ and $\cA = \End(\cLd)$ where $V$ is a finite-dimensional $K$-vector space and $\cLd$ is a lattice chain in $V$. We may also assume that $\cT \cong \cO^n\cong \cT'$ where $d = \dim(V)$. Let $e$ be the period of $\cLd$ and let $\barcL_i \colon = \cL_i/\cL_{i-1}$. Consider the $\barcA \colon = \cA/\fP$-module $\barcL \colon =  \bigoplus_{i=1}^e \, \barcL_i$. As a $\barcT \colon = \cT/\Rad(\cT)$- and $\barcT' \colon = \cT'/\Rad(\cT')$-module it is free of rank 1 (by Lemma \ref{lemma:maxcommsepsubalgebra3} of the appendix). Hence there exists an isomorphism $\overline{\Theta}: \barcT\to \barcT'$ such that $\overline{\Theta}(\bar{t}) x = \bar{t} x$ for all $\bar{t}\in \barcT, x\in \barcL$. We choose a lifting $\Theta$ of $\overline{\Theta}$, i.e.\ an isomorphism of $\cO$-algebras $\Theta: \cT\to \cT'$ which reduces to $\overline{\Theta}$ modulo $\varpi$. Then for any $i\in \bZ$ we have
\begin{equation}
\label{eqn:maxtor1}
\Theta(t) x = t x\qquad \mbox{for all $t\in \cT, x\in \barcL_i$}
\end{equation}
Since $\cL_0$ is a free $\cT$- and $\cT'$-module of rank 1 there exists $f\in \Aut_{\cO}(\cL_0)\subseteq A^*$ such that $f(tx) = \Theta(t) f(x)$ for all $t\in \cT, x\in \cL_0$. Hence $\Theta(t) = f t f^{-1}$ for all $t\in \cT$ and therefore $\cT' = f \cT f^{-1}$. We claim that $f\in\cA^*$, i.e.\ $f(\cL_i)= \cL_i$ for all $i\in \bZ$. For that it is enough to see that $f(\cL_i)\subseteq \cL_i$ for all $i=1,2, \ldots e$ and in fact for $i=1$ (by induction). Note that $f(\cL_1) \subseteq f(\cL_e) =\varpi^{-1} f(\cL_0) = \cL_e$. Choose $i\in \{1,2,\ldots,e\}$ minimal with $f(\cL_1) \subseteq \cL_i$ and assume that $i \ge 2$. Then $f$ induces a nontrivial $\cT$-linear homomorphism $\bar{f}: \barcL_1 \lra \barcL_i$ such that
\[
\bar{f}(tx) = \Theta(t) \bar{f}(x) = t \bar{f}(x)\qquad \mbox{for all $t\in \cT, x\in \barcL_i$}.
\]
On the other hand since $\barcL$ is a free $\barcT$-module of rank 1 we have $\Hom_{\cT}(\barcL_1, \barcL_i)\break =0$, a contradiction. This proves $f\in\cA^*$.\enddemo

We need the following two simple Lemmas in section \ref{subsection:special}. 

\begin{lemma} 
\label{lemma:stablyfree4}
Suppose that $\cA$ is principal and let $\cT$ be a maximal torus in $\cA$. Let $\cM$ be an $\cA$-lattice and put $\barcT \colon = \cT/\Rad(\cT)$. The following conditions are equivalent.

\noi (i) $\cM$ is stably free.

\noi (ii) $\cM/\fP\cM$ is a free $\barcT$-module.
\end{lemma}

The proof will be left to the reader.

\begin{lemma} 
\label{lemma:stablyfree5}
Assume that $\cA$ is principal and let $\cT$ be a maximal torus in $\cA$. Let $0\to \cM'\to \cM \to \cN\to 0$ be a short exact sequence of $\cA$-modules and assume that $\cM$ is a stably free $\cA$-lattice and $\varpi \cN = 0$. The following conditions are equivalent.

\noi (i) $\cM'$ is stably free.

\noi (ii) $\cN$ is a free $\barcT$-module.
\end{lemma}

{\em Proof.} By using the exact sequence
\[
0  \lra \Ker(\cN \otimes_{\cA} \fP \to \cN) \to \cM'/\fP\cM'\lra \cM/\fP\cM \lra \cN/\fP\cN \lra 0 
\]
we see that
\begin{align*}
[\cM'/\fP\cM'] & =  [\cM/\fP\cM] + [\Ker(\cN \otimes_{\cA} \fP \to \cN)] - [\cN/\fP\cN]\\ 
& = [\cM/\fP\cM] + [\cN \otimes_{\cA} \fP] - [\cN]
\end{align*}
in the Grothendieck group $K_0(\barcT)$. Note that $[\cN]= [\cN \otimes_{\cA} \fP]$ if and only if $\cN$ is a free $\barcT$-module. Hence (ii) is equivalent to the equality $[\cM'/\fP\cM'] = [\cM/\fP\cM]$ in $K_0(\barcT)$. The assertion follows from \ref{lemma:stablyfree4}.\enddemo

\subsection{Local theory of invertible Frobenius bimodules}
\label{subsection:fqalg}

Let $\cO$ be a henselian discrete valuation ring with quotient field $K$, maximal ideal $(\varpi)= \fp$ and residue field $k= \cO/\fp$. We assume that $k$ is finite of characteristic $p$. Let $v_K$ be the normalized valuation of $K$. We denote by  $\inv$ the canonical isomorphism $\Br(K) \to \bQ/\bZ$ of class field theory. Let $\cO'$ be a finite {\'e}tale local $\cO$-algebra with quotient field $K'$. By $\sigma\in \Gal(K'/K)$ we denote Frobenius isomorphism (i.e.\ $\sigma(x) \equiv x^{\sharp(k)} \mod \fp$). For an $\cO'$-module $\cM$ we write $\cM_{K'}$ for $\cM\otimes_{\cO'} K'$ and $\sM$ for $\cM\otimes_{\cO',\sigma} \cO'$ (or equivalently $\sM = \cM$ with the new $\cO'$-action $x\cdot m = \sigma(x) m$). 

Let $A$ be a central simple $K$-algebra of dimension $d^2$ and $\cA$ a principal $\cO$-order in $A$ with radical $\fP$ and index $e = e(\cA)$ (note that we have $e \inv(A) = 0$). Let $\cM$ be a free right $\cA_{\cO'}$-module of rank 1 together with an isomorphism of $\cA_{\cO'}$-modules
\[
\phi: \scM\fP^m \lra \cM
\]
for some $m\in \bZ$. We set
\[
\cB \colon = \End_{\cA_{\cO'}}(\cM, \phi) = \{ f\in \End_{\cA_{\cO'}}(\cM)\mid\, \phi \circ \sig{f} = f \circ \phi\}.
\]
\begin{lemma} 
\label{lemma:frobinv1} 
The $\cO$-algebra $\cB$ is a principal order of index $e$ in the central simple $K$-algebra $B\colon = \cB_K$ of dimension $d^2$. We have 
\begin{equation}
\label{eqn:frobinvtwist}
\inv(B) = \inv(A) + \frac{m}{e} \mod\bZ
\end{equation}
\end{lemma}

{\em Proof.} Let $\phi_{K'} \colon = \phi\otimes_{\cO} \id_{K}:\sig(M_{K'})\to M_{K'}$. By Lemma \ref{lemma:stablyfree1} the $\cO'$-algebra $\cB' \colon = \End_{\cA_{\cO'}}(\cM)$ is a principal $\cO'$-order of index $e$ in $B' \colon = \End_{A_{K'}}(\cM_{K'})$. Define a $\sigma$-linear isomorphism $\psi: B' \to B'$ by $\psi(f) \colon = \phi_{K'}\circ \sig{f} \circ \phi_{K'}^{-1}$. We have 
\[
\cB = \{b\in \cB'\mid\, \psi(b) = b\} \quad\mbox{and}\quad \cB' \cong \cB_{\cO'}.
\] 
Together with Lemma \ref{lemma:hered2} this implies the first statement. The proof of the second assertion will be left to the reader.\enddemo

Conversely suppose that we have given a second central simple $K$-algebra $B$ of dimension $d^2$ and a principal $\cO$-order $\cB$ in $B$ of index $e$. We also assume that $[K':K]$ is a multiple of the order of $[B\otimes A^{\opp}]$ in $\Br(K)$. Let $m$ be any integer such that (\ref{eqn:frobinvtwist}) holds. 

\begin{lemma} 
\label{lemma:frobinv2} 
There exists an invertible $\cB_{\cO'}$-$\cA_{\cO'}$-bimodule $\cM$ and an isomorphism of bimodules
\[
\phi: \scM\fP^m \lra \cM.
\]
\end{lemma}

{\em Proof.} By Proposition \ref{proposition:hermor} the principal orders $\cB_{\cO'}$ and $\cA_{\cO'}$ are Morita equivalent. Let $\cM$ be an invertible $\cB_{\cO'}$-$\cA_{\cO'}$-bimodule. Then $\scM$ is invertible as well. Hence there exists an isomorphism $\phi': \scM\fP^{m'} \to \cM$ for some $m'\in \bZ$. By \ref{lemma:frobinv1} we have
\[
\inv(B) = \inv(A) + \frac{m'}{e} \mod\bZ
\]
and therefore $m\equiv m' \mod e$. Put $\phi \colon = \varpi^{\frac{m'-m}{e}} \phi'$.\enddemo

For the rest of this section we assume that $\cO$ is an $\bFq$-algebra ($q = p^r$ for some $r\in \bZ$) and let $k'$ be an (possibly infinite) algebraic extension of $k$ whose degree (over $k$) is a multiple of $e$. Let $\cO'\colon = \cO\otimes_{\bFq} k'$ and $\sigma\colon = \id_{\cO}\otimes \Frob_q\in \Gal(\cO'/\cO)$. For $\rho\in\Hom_{\bFq}(k,k')\cong \Hom_{k'}(k\otimes_{\bFq} k',k')$ we denote the kernel of $\cO' \to k\otimes_{\bFq} k'\to k'$ by $\fp'_{\rho}$ and we set $\cO'_{\rho} \colon = \cO'_{\fp'_{\rho}}$. Then $\cO'_{\rho}$ is a (pro-)finite (pro-){\'e}tale local $\cO$-algebra whose degree is a multiple of $e$ and $\cO' \cong \bigoplus_{\rho} \, \cO'_{\rho}$. Similarly 
\[
\cA_{\cO'} = \bigoplus_{\rho} \, \cA'_{\rho} \qquad\mbox{with $\cA'_{\rho} = \cA_{\cO_{\rho}'}$}
\]
and $\fP_{\cO'} = \Rad(\cA_{\cO'})$ is equal to the product $\prod_{\rho} \, \fP'_{\rho}$ where $\fP'_{\rho}$ denotes the maximal invertible two-sided ideal $\Ker(\cA_{\cO'}\to \cA'_{\rho}/\Rad(\cA'_{\rho}))$ of $\cA_{\cO'}$. For the distinguished element $\iota \colon = \mbox{incl}\,:k\hookrightarrow k'$ in $\Hom_{\bFq}(k,k')$ we put $\fp' = \fp'_{\rho}$ and $\fP' \colon = \fP'_{\iota}$. Let $\cM$ be a free right $\cA_{\cO'}$-module of rank $1$. For $m\in \bZ$ the $\cA_{\cO'}$-module $\sig(\cM(\fP')^m)$ is also free of rank $1$. Hence there exists an isomorphism 
\[
\phi:\sig(\cM(\fP')^m)\lra \cM.
\]
If we set 
\[
\cB\colon = \End_{\cA_{\cO'}}(\cM, \phi) = \{ f\in \End_{\cA_{\cO'}}(\cM)\mid\, \phi \circ \sig{f} = f \circ \phi\}.
\]
then one can deduce easily from Lemma \ref{lemma:frobinv1} that $\cB$ is a principal $\cO$-order of index $e$ in the central simple $K$-algebra $B = \cB_K$ and that equation (\ref{eqn:frobinvtwist}) holds. 

Conversely given such a principal $\cO$-order $\cB$ of index $e$ and $m\in \bZ$ such that (\ref{eqn:frobinvtwist}) holds there exists a pair $(\cM, \phi)$ as above with $\cB = \End_{\cA_{\cO'}}(\cM, \phi)$. To see this let $\cM$ be any invertible $\cB_{\cO'}$-$\cA_{\cO'}$-bimodule. Since $\scM$ is invertible as well we have
\begin{equation}
\label{eqn:frobinvtwist2}
\sig(\cM\prod_{\rho}(\fP'_{\rho})^{m_{\rho}}) \cong \cM
\end{equation}
for certain $m_{\rho}\in \bZ$. Since $\sig(\fP'_{\rho}) \cong \fP'_{\Frob_q\circ \rho}$ we may assume -- after replacing $\cM$ by $\cM\fA$ for a suitable invertible two-sided $\cA_{\cO'}$-ideal $\fA$ -- that $m_{\rho}=0$ for all $\rho\in\Hom_{\bFq}(k,k')$ except $\rho = \iota$. As in the proof of Lemma \ref{lemma:frobinv2} we deduce 
\[
\inv(B) = \inv(A) + \frac{m_{\iota}}{e}\mod\bZ
\]
hence $m_{\iota}\cong m \mod e$ and therefore $(\fP')^{m_{\iota}} \cong (\fP')^m$. Hence there also exists an isomorphism $\sig(\cM(\fP')^m)\cong\cM$.

\begin{definition}
\label{definition:invfrobloc}
A pair $(\cM, \phi)$ consisting of an invertible $\cB_{\cO'}$-$\cA_{\cO'}$-bimodule $\cM$ and an isomorphism $\phi:\sig(\cM(\fP')^m)\to \cM$ is called an invertible $\phi$-$\cA$-$\cB$-bimodule of slope $-\frac{m}{e}$ over $\cO'$.
\end{definition}

We have seen that an invertible $\phi$-$\cB_{\cO'}$-$\cA_{\cO'}$-bimodule of a given slope $r\in \bQ$ exists if and only if $r = \inv(A)-\inv(B) \mod\bZ$. It is also easy to see that any two invertible $\phi$-$\cB_{\cO'}$-$\cA_{\cO'}$-bimodules of the same slope differ (up to isomorphism) by a fractional $\cA$-ideal. This implies that if $k''$ is an algebraic extension of $k'$ and $\cO''= \cO\otimes_{\bFq} k''$ then any $\phi$-$\cA$-$\cB$-bimodule over $\cO''$ is obtained by base change from an $\phi$-$\cA$-$\cB$-bimodule over $\cO'$.

\begin{remark}
\label{remark:phief}
\rm Assume that $[k':k]= e$ and let $n = [k':\bFq]$. 
 Let $(\cM, \phi)$ be an invertible $\phi$-$\cB_{\cO'}$-$\cA_{\cO'}$-bimodule of slope $-\frac{m}{e}$. For $r\in \bZ/n\bZ$ we put $\fP'_r=  \srfP$. We have $\prod_{r\in \bZ/n\bZ}\, \fP'_r = \fP_{\cO'}^e = \fp \cA_{\cO'}$. For each two-sided invertible ideal $\fA'$ of $\cA_{\cO'}$ and $r\in \bZ/n\bZ$, the map $\phi$ induces isomorphisms
$(\srcM) \fA' \to (\srmecM) \fA'{\fP'_r}^m$ which will be also denoted by $\phi$. Consider the map
\begin{equation}
\label{eqn:phief}
\phi^n: \cM = (\sncM)\stackrel{\phi}{\lra} (\snmecM) {\fP'_n}^m \stackrel{\phi}{\lra}\ldots \lra \cM \prod_{r\in \bZ/n\bZ}\, {\fP'}_r^m = \cM\fp^m
\end{equation} 
Since (\ref{eqn:phief}) is $\cB_{\cO'}$-$\cA_{\cO'}$-bilinear and commutes with $\phi$ there exists an element $x\in K$ with $v_K(x) =m$ such that (\ref{eqn:phief}) is given by multiplication with $x$. This fact will be used later when we discuss level structure at the pole of $\cA$-elliptic sheaves.
\end{remark}

\section{Global theory of hereditary orders}
\label{section:globhereditary}

In this section we study hereditary orders in a central simple algebras over a function field of one variable (though most results hold also for number fields). We shall show that two hereditary orders are Morita equivalent if their generic fibers are equivalent and all their local indices are the same. Furthermore any such hereditary order is Morita equivalent to a locally principal one. We will then study the Picard group of a locally principal order $\cA$ and introduce the notion of $\cA$-degree of a locally free $\cA$-module of finite rank. In the final part we will introduce the notion of a special $\cA$-module. 

In this chapter $k$ denotes a fixed perfect field of cohomological dimension $\le 1$ and $X$ a smooth projective geometrically connected curve over $k$ with function field $F$. For $x\in |X|$ we denote by $\cO_x$ the completion of $\cO_{X,x}$ and by $F_x$ the quotient field of $\cO_x$. The maximal ideal of $\cO_x$ will be denoted by $\fp_x$. If $\cV$ is a coherent $\cO_X$-module then we set $\cV_x=\cV\otimes_{\cO_X} \cO_x$ and if $V$ is a finite-dimensional $F$-vector space we put $V_x = V\otimes_F F_x$. 

\subsection{Morita equivalence.} 

Let $V$ be a finite-dimensional $F$-vector space. The set of locally free coherent $\cO_X$-modules $\cV$ with generic fiber $\cV_{\eta} = V$ is in one-to-one correspondence with the set of $\cO_x$-lattices $\cV_x$ in $V_x$ for all $x\in |X|$ such that there exists an $F$-basis $B$ of $V$ with $\cV_x = \sum_{b\in B} \cO_x b$ for almost all $x$. Consequently if $U\seq X$ is an open subscheme then there is a one-to-one correspondence between coherent and locally free $\cO_X$-modules $\cV$ and coherent and locally free $\cO_U$-module $\cV_U$ and together with an $\cO_x$-lattice $\cV_x$ in $\cV_U\otimes F_x$ for all $x\in X-U$. 

Let $A$ be a central simple $F$-algebra and $\cA$ a hereditary $\cO_X$-order in $A$. We put $e_x(\cA) \colon = e(\cA_x)$. There are only finitely many points $x\in |X|$ with $e_x(\cA)>1$. Define the divisor $\Bad(\cA)$ as $\Bad(\cA)\colon = \sum_{x\in |X|}\, (e_x(\cA)-1) x$. If $k$ is finite and $x\in |X|$ then $\inv_x(A)$ denotes the image of the class of $A_x$ under the canonical isomorphism of class field theory $\Br(F_x) \to \bQ/\bZ$. 

\begin{proposition}
\label{proposition:moritaglobal}
Let $A_1,A_2$ be central simple algebras over $F$ and let $\cA_1$ and $\cA_2$ be hereditary $\cO_X$-orders in $A_1$ and $A_2$ respectively. The following conditions are equivalent.

\noi (i) $\cA_1$ and $\cA_2$ are equivalent.

\noi (ii) $A_1$ and $A_2$ are equivalent and $(\cA_1)_x$ and $(\cA_2)_x$ are equivalent for all $x\in |X|$.

\noi (iii) $A_1$ and $A_2$ are equivalent and $\Bad(\cA_1) = \Bad(\cA_2)$.

Moreover if $k$ is a finite field then the above conditions are also equivalent to:

\noi(iv) $\inv_x(A_1) = \inv_x(A_2)$ for all $x\in |X|$ and $\Bad(\cA_1) = \Bad(\cA_2)$.
\end{proposition}

{\em Proof.} (i) $\Rightarrow$ (ii) is clear.
\noi (ii) $\Leftrightarrow$ (iii) follows from Proposition \ref{proposition:hermor} and (iii) $\Leftrightarrow$ (iv) from the Theorem of Brauer--Hasse--Noether. It remains to show that (ii) implies (i). Let $U$ be an affine open subscheme of $X$ contained in the complement of $\Bad(\cA_1) = \Bad(\cA_2)$ in $X$. By (\cite{reiner}, 21.7) ${\cA_1}|_U$ and ${\cA_2}|_U$ are Morita equivalent. Let $\cI_U$ be an invertible $\cA_1|_U$-$\cA_2|_U$-bimodule and let $\cI_x$ be an invertible $(\cA_1)_x$-$(\cA_2)_x$-bimodule for each $x\in X-U$. Since there is only one invertible $(A_1)_x$-$(A_2)_x$-bimodule up to isomorphism we may assume that $\cI_x\otimes F_x = \cI_U\otimes F_x$ i.e.\ that $\cI_x$ is a lattice in $\cI_U\otimes F_x$. It is easy to see that the locally free $\cO_X$-module $\cI$ corresponding to $\cI_U$ and the $\cI_x$, $x\in X-U$ is then an invertible $\cA_1$-$\cA_2$-bimodule.\enddemo

 A {\it locally principal $\cO_X$-order $\cA$} is a hereditary $\cO_X$-order in a central simple $F$-algebra $A$ such that $\cA_x$ is principal for all $x\in |X|$. The rank of $\cA$ is its rank as an $\cO_X$-module, hence $=\dim_F(A)$. If $\cA$ is a hereditary $\cO_X$-order in $A$ then it is locally principal if for example $\cA_x$ is either maximal or $e_x(\cA) =d$ for all $x\in|\Bad(\cA)|$.

Suppose that $\cA$ is a locally principal $\cO_X$-order of rank $d^2$.
We define two positive integers $e(\cA), \delta(\cA)$ by
\begin{equation}
\label{eqn:e}
e(\cA) \colon = \lcm\{e_x(\cA)\mid\, x\in |X|\}
\end{equation}
\[
\delta(\cA) \colon = \lcm\{\,\mbox{numerator of $\frac{e_x(\cA)}{\deg(x)}$}\,\mid\, x\in |X|\}
\]
According to Lemma \ref{lemma:principal} we have $\delta(\cA)\mid e(\cA)\mid d$. If $\cA$ is locally principal then one can easily see that 
\[
\deg(\cA) = - \frac{d^2}{2} \sum_{x\in |X|}\, (1-\frac{1}{e_x(\cA)})\deg(x).
\]
In particular if $\cB$ is a second locally principal $\cO_X$-order of rank $d^2$ with $\Bad(\cA) = \Bad(\cB)$ then 
\begin{equation}
\label{eqn:degprinc}
\deg(\cA) = \deg(\cB).
\end{equation}

\begin{corollary}
\label{corollary:gloprincipal}
Let $\cA$ be a hereditary $\cO_X$-order in a central simple $F$-algebra $A$. Then there exists a locally principal $\cO_X$-order $\cD$ which is Morita equivalent to $\cA$. In fact $\cD$ can be chosen such that $\rank_{\cO_X}(\cD) = e(\cA)^2$.
\end{corollary}

{\em Proof.} That $\cA$ is equivalent to a locally principal $\cO_X$-order follows easily from the corresponding local statement \ref{corollary:principal}. In fact if $B\colon = M_e(A)$ then for all $x\in |\Bad(\cA)|$ we can pick a principal $\cO_x$-order $\cB_x$ in $B_x$ equivalent to $\cA_x$. If $U \colon = X-|\Bad(\cA)|$ and $\cB_U$ is a maximal $\cO_U$-order in $B$ then there exists a uniquely determined hereditary $\cO_X$-order $\cB$ in $B$ with $\cB\otimes_{\cO_X} \cO_x = \cB_x$ for all $x\in |\Bad(\cA)|$ and $\cB|_U = \cB_U$. The order $\cB$ is locally principal and equivalent to $\cA$ by \ref{proposition:moritaglobal}.

Thus to prove the second statement we may assume that $\cA$ is locally principal. Let $\cI$ be a locally stably free $\cA$-module which is of rank $de$ as an $\cO_X$-module. By Lemma \ref{lemma:stablyfree1} and \ref{proposition:moritaglobal} above it follows that $\cD\colon = \uEnd_{\cA}(\cI)$ is a locally principal $\cO_X$-order in $\End_A(\cI_{\eta})$. Moreover $\cD$ is equivalent to $\cA$ and $\rank_{\cO_X}(\cD) = e(\cA)^2$.\enddemo
 
\subsection{Locally free $\cA$-modules}

\paragraph{The Picard group of a locally principal order.} 

In this section $\cA$ denotes a locally principal $\cO_X$-order of rank $d^2$. We are going to compute the Picard group of $\cA$. Define
\[
\DivcA \colon = \{\sum_{x\in |X|}\, n_x x\in \Div(X)\otimes \bQ\mid \, e_x(\cA) n_x\in \bZ\, \,\,\forall\, x\in |X|\}.
\]
Note that $\deg(\DivcA) = \frac{1}{\delta(\cA)}\bZ$. For a divisor $D = \sum_{x\in |X|}\, n_x x \in \DivcA)$ we denote by $\cA(D)$ the invertible $\cA$-$\cA$-bimodule given by $\cA(D)|_{X-|D|} = \cA|_{X-|D|}$ and $\cA(D)_{x} = \fP_{\cA_{x}}^{-n_x e_x(\cA)}$ for all $x\in |X|$. If $D\in \Div(X)$ then $\cA(D) = \cA\otimes_{\cO_X} \cO_X(D)$.

\begin{proposition}
\label{proposition:picprinc}
The sequence 
\[
0 \lra F^*/k^* \stackrel{\dvsor}{\lra} \DivcA\stackrel{D\mapsto \cA({\mathbf D})}{\lra} \Pic(\cA)\lra 0
\]
is exact. 
\end{proposition}

{\em Proof.} This follows from (\cite{reiner}, 40.9).\enddemo

We also need to consider the group of isomorphism classes of invertible $\cA$-$\cA$-bimodules with level structure and give a description of it as an idele class group. Let $I = \sum_x \, n_x x$ be an effective divisor on $X$. The corresponding finite closed subscheme of $X$ will be also denoted by $I$. A {\it level-$I$-structure} on an invertible $\cA$-$\cA$-bimodule $\cL$ is an isomorphism $\beta: \cA_I \to \cL_I$ of right $\cA_I$-modules. We denote by $\Pic_I(\cA)$ the set of isomorphism classes of invertible $\cA$-$\cA$-bimodules with level-$I$-structure. If $(\cL_1, \beta_1),(\cL_2, \beta_2)$ are invertible $\cA$-$\cA$-bimodules with level-$I$-structures we define the level-$I$-structure $\beta_1\beta_2$ on $\cL_1\otimes_{\cA}\cL_2$ as the composite
\begin{equation}
\label{eqn:tensorlevel}
\begin{CD}
\beta_1\beta_2: \cA_I @> \beta_2 >> (\cL_2)_I = \cA_I\otimes_{\cA_I} (\cL_2)_I @> \beta_1\otimes \id >> (\cL_1\otimes_{\cA}\cL_2)_I
\end{CD}
\end{equation}
thus defining a group structure on $\Pic_I(\cA)$. Note that unlike $\Pic(\cA)$, $\Pic_I(\cA)$ is in general not abelian. In fact we have a short exact sequence 
\begin{equation}
\label{eqn:pici}
0\lra\cA_I^*/k^*\lra \Pic_I(\cA)\lra\Pic(\cA)\lra 0
\end{equation}
where the first map is given by $a\in \cA_I^*\mapsto (\cA, l_a: \cA_I \stackrel{a\cdot}{\lra} \cA_I)$ . 

Let $U_I(\cA) \colon = \Ker(\prod_{x\in |X|} \, \cA_x^*\to \prod_{x\in |X|}\,(\cA_x/\fp_x^{n_x}\cA_x)^* = \cA_I^*)$ and let
\[
\cC_I(\cA)\colon = \mbox{$(\prod_{x\in |X|}'\, N(\cA_x))/U_I(\cA)F^*$}
\]
where $\prod'_{x\in |X|}\, N(\cA_x))$ denotes the restricted direct product of the groups $\{N(\cA_x))\}_{x\in |X|}$ with respect to $\{\cA_x^*\}_{x\in |X|}$. Given $a=\{a_x\}_x\in \prod'_{x\in |X|}\, N(\cA_x)$ we put $\dvsor(a) = \sum_{x\in |X|} v_{\cA_x}(a_x) x$. Left multiplication by $a$ induces a level-$I$-structure $\beta_a: \cA_I \to \cA(\dvsor(a))_I$.

\begin{corollary}
\label{corollary:piclevel}
The assignment $a\mapsto (\cA(\dvsor(a)), \beta_a)$ induces an isomorphism $\cC_I(\cA)\cong \Pic_I(\cA)$.
\end{corollary}

\paragraph{Relative divisors and invertible bimodules.} 

Let $S$ be a $k$-scheme and let $\pi: X\times S \to S$ be the projection. We need to define the bimodule $\cA(D)$ also for elements of a certain group  of relative divisors $\Div(\cA\boxtimes \cO_S)$. For the latter we use the following ad hoc definition. Assume first that $S$ is of finite type over $k$. Let $\cal S$ be the collection of all connected components of $x\times S$ where $x$ runs through all closed points of $X$. Thus if $S'\in \cal S$ there exists a unique closed point $x \colon = \pi(S')$ with $S'\subseteq x\times S$. We set
\[
\Div(\cA\boxtimes \cO_S)\colon = \bigoplus_{S'\in \cS} \, \frac{1}{e_{\pi(S')}(\cA)} \bZ.
\]
Let $R$ be the integral closure of $k$ in $\Gamma(S, \cO_S)$. Note that for $x\in |X|$ the set of open and closed subschemes of $x\times S$ corresponds to the set of idempotents in $k(x)\otimes_k R$. If $f: S_1\to S_2$ is a morphism of $k$-schemes there is an obvious notion of a pull-back $f^*:\Div(\cA\boxtimes \cO_{S_2})\to\Div(\cA\boxtimes \cO_{S_1})$. For an arbitrary $k$-scheme we define $\Div(\cA\boxtimes \cO_S)$ as the direct limit of $\Div(\cA\boxtimes \cO_{S'})$ over the category of pairs $(S', g)$ consisting of a $k$-scheme $S'$ of finite type and a morphism $g: S \to S'$ in $\Sch/k$.

Let $S\in \Sch/k$. A $k$-morphism $x_S: S\to X$ which factors as $S\to \Spec k(x) \to X$ for some $x\in |X|$ yields an element -- denoted by $x_S$ as well -- of the group $\Div(\cA\boxtimes \cO_S)$. For that we can assume that $S$ is of finite type. Since the graph $\Gamma_{x_S} = (x_S, \id_S): S \lra X\times S$ is an open and closed subscheme of $x\times S$ it is a disjoint union of connected components and we define $x_S\in \Div(\cA\boxtimes \cO_S)$ to be the sum of these components.

There exists a unique homomorphism 
\begin{equation}
\label{eqn:divpic}
\Div(\cA\boxtimes \cO_S)\to \Pic(\cA\boxtimes \cO_S),\quad D \mapsto (\cA\boxtimes \cO_S)(D)
\end{equation}
compatible with pull-backs which agrees with the previously defined map in case $S = \Spec k'$ for a finite extension $k'/k$. It suffices to define (\ref{eqn:divpic}) for $\frac{1}{e_x(\cA)}D$, where $D$ is a connected component of $x\times S$ for some $x\in |X|$. It is also enough to consider the case where $S$ is connected and of finite type over $k$. Let $R$ be the integral closure of $k$ in $\Gamma(S, \cO_S)$. Then $\Spec R$ is connected and finite over $\Spec k$, i.e.\ $R$ is an artinian finite local $k$-algebra. Let $k'$ denotes the residue field of $R$. Since $k$ is perfect the canonical projection $R \to k'$ has a unique section. Therefore the structural morphism $S\to \Spec k$ factors as $S\to \Spec k' \to \Spec k$. Thus by replacing $k$, $X$ and $\cA$ by $k'$ and $X_{k'}$ and $\cA\boxtimes k'$ respectively we can assume that the residue field of $R$ is $k$. However, in this case, $x\times S$ is connected for all $x\in |X|$, hence $D = x\times S$ with $x = \pi(D)$. So we are forced to define $(\cA\boxtimes \cO_S)(\frac{1}{e_x(\cA)}D)\colon = \pi^*(\cA(\frac{1}{e_x(\cA)}x))$. 

\paragraph{$\cA$-rank and $\cA$-degree.}

\noi Let $f:S \to X$ be a morphism. For $\cE$ in $_{f^*(\cA)}\Mod$ and $\cF$ in $\Mod_{f^*(\cA)}$ we put $\cE\otimes_{\cA} \cF \colon = \cE\otimes_{f^*(\cA)} \cF$. If $D= \sum_{x\in |X|}\, n_x x\in \DivcA$ we set $\cE(D) \colon = \cE \otimes_{\cA} f^*(\cA(D))$ and $\cF(D) \colon = f^*(\cA(D)\boxtimes\cO_S) \otimes_{\cA} \cF$.

Let $S$ be a $k$-scheme. We denote by $\cAVect(S)$ (resp.\ $\VectcA(S)$) the category coherent and locally free left (resp.\ right) $\cA\boxtimes \cO_S$-modules. For $\cF$ in $\cAVect(S)$ or $\VectcA(S)$ let $\rankcA \cF$ be the locally constant function $s\mapsto \rank_{\cA\boxtimes k(s)}(\cF|_{X\times s})$ on $S$ (hence $\rankcA \cF$ can be viewed as an element of $\bZ^{\pi_0(S)}$). For a positive integer $r$ we denote by $\cAVectr(S)$ (resp.\ $\VectcAr(S)$) the subcategory of $\cF\in \cAVect(S)$ (resp.\ $\cF\in \VectcA(S)$) with $\rankcA \cF = r$.

Let $\cF$ be a locally free $\cA\boxtimes \cO_S$-module of rank  $r$. Define $\detcA \cF$ as the image of the isomorphism class of $\cF$ (viewed as an element of $H^1(X\times S, \Gl_r(\cA\boxtimes \cO_S)$) under the map
\[
H^1_{\zar}(X\times S, \Gl_r(\cA\boxtimes \cO_S)) \lra  H^1_{\zar}(X\times S,\cO^*)= \Pic(X\times S)
\]
induced by the reduced norm $\Nrd:M_r(A) \to F$. We obtain a locally constant function 
\[
\degcA(\cF): S \to \frac{1}{d}\bZ, s\mapsto \deg((\detcA \cF)|_{X\times s})
\]
It is easy to see that 
\[
\degcA(\cF) = \frac{1}{d^2}(\deg(\cF) - \rankcA(\cF) \deg(\cA)).
\]
In particular since $\deg(\cA(D)) = \deg(\cA) + d^2\deg(D)$ we have 
\[
\degcA(\cA(D)) = \deg(D)
\]
for all $D\in \DivcA$.

\begin{lemma}
\label{lemma:principaldegree}
(a) Let $0\to \cF_1\to\cF_2\to\cF_3\to 0$ be a short exact sequence of coherent and locally free $\cA\boxtimes \cO_S$-modules. Then 
\[
\degcA(\cF_2)= \degcA(\cF_1)+ \degcA(\cF_3).
\]
(b) Let $\cE$ be an object of $\VectcAr(S)$ and $\cF$ be an object of $\cAVects(S)$. Then 
\[
 \frac{1}{d^2}(\deg(\cE\otimes_{\cA} \cF) - r s \deg(\cA)) = r \degcA(\cF) + s\degcA(\cE).
\]
(c) Let $\cB$ be a second locally principal $\cO_X$-order of rank $d^2$ equivalent to $\cA$. Let $\cE$ be an object of $\VectcAr(S)$ and let $\cI$ be an invertible $\cA$-$\cB$-bimodule. Then 
\[
\degcB(\cE\otimes_{\cA} \cI) = \degcA(\cE) + r \degcA(\cI).
\]
(d) Let $\cE$ be an object of $\cAVectr(S)$ and $D\in \DivcA$. Then
\[
\degcA(\cE(D)) = \degcA(\cE) + r\deg(D)
\]
\end{lemma}

{\em Proof.} (a) is obvious, (c) follows from (b) and (\ref{eqn:degprinc}) and (d) is a special case of (c). Note that by \ref{corollary:invlocfree} the bimodule $\cI$ in (c) is a locally-free left $\cA$- and right $\cB$-module of rank 1.

For (b) it is enough to consider the case when $S$ is a connected $k$-scheme of finite type and therefore -- by choosing a fixed closed point $s\in S$ and taking the base change $\Spec k(s) \to \Spec k$ -- to consider the case $S = \Spec k$. If 
\begin{equation}
\label{eqn:modification}
\xymatrix@-1.5pc{
&& \cE \\
\cE''\ar[urr]  \ar[drr] \\
&& \cE'
}
\end{equation}
is a diagram of locally free $\cA\boxtimes \cO_S$-modules of the same rank $r$ and injective $\cA\boxtimes \cO_S$-linear homomorphisms then it is easy to see that (b) holds for $\cE$ if and only if it holds for $\cE'$. Since $\cE|_U \cong \cA^r|_U$ for some non-empty open subscheme $U\seq X$ there exists a diagram (\ref{eqn:modification}) with $\cE' = \cA^r$. The assertion follows.\enddemo

It follows from \ref{proposition:picprinc} or \ref{lemma:principaldegree} (b) that $\degcA: \Pic(\cA) \to \bQ$ is a homomorphism. We denote its kernel by $\Pic_0(\cA)$. Also if $I\in \Div(X)$, $I\ge 0$ we let $\Pic_{I,0}(\cA)$ be the subgroup of $(\cL, \beta)\in\Pic_I(\cA)$ with $\degcA(\cL)=0$. The image $\degcA(\Pic(\cA))$ is equal to $\frac{1}{\delta(\cA)}\bZ$. 

\begin{remark}
\label{remark:stronglyequivalent} 
{\rm Let $\cA, \cB$ be locally principal $\cO_X$-order of rank $d^2$ and suppose that $\cA$ and $\cB$ are equivalent. The set of isomorphism classes of invertible $\cA$-$\cB$-bimodule has a simple transitive left $\Pic(\cA)$-action. Hence for any two invertible $\cA$-$\cB$-bimodule $\cI, \cJ$ the degrees $\degcA(\cJ)$ and $\degcA(\cI)$ differ by a multiple of $\frac{1}{\delta(\cA)}$. Call $\cA$ and $\cB$ strongly Morita equivalent if there exists an invertible $\cA$-$\cB$-bimodule $\cI$ with $\degcA(\cI) =0$. It is easy to see that a given Morita equivalence class of locally principal $\cO_X$-orders of rank $d^2$ decomposes into $\frac{e}{\delta}$ strong equivalence classes (where $e$ and $\delta$ are defined in (\ref{eqn:e})).}
\end{remark}

\subsection{Special $\cA$-modules}
\label{subsection:special}

Let $\cA$ be a locally principal $\cO_X$-order of rank $d^2$. If $g:U \to X$ is an {\'e}tale morphism then a maximal torus in $\cA_U \colon = g^*(\cA)$ is a maximal commutative {\'e}tale $\cO_U$-subalgebra of $\cA_U$. 

\begin{definition}
\label{definition:special2}
A right $\cA\boxtimes \cO_S$-module $\cK$ is called special of rank $r$ if the following hold:
\begin{enumerate}
\item[(i)] $\cK$ is coherent as an $\cO_{X\times S}$-module and the map $\Supp(\cK)\hra X\times S \to S$ is an isomorphism. Hence $\Supp(\cK)$ is the image of the graph of a morphism $N = N(\cK): S \to X$ and
$\cK$ is the direct image of a $N^*(\cA)$-module -- also denoted by $\cK$ -- by the graph $\Gamma_{N} = (N, \id_S): S \to X\times S$.
\item[(ii)] Consider $\cK$ as a sheaf on $S$ as in (i). For any {\'e}tale morphism $g:U\to X$ and maximal torus $\cT$ of $\cA_U$, $(g_S)^*(\cK)$ is a locally free $(N_U)^*(\cT)$-module of rank $r$. Here $g_S$ (resp.\ $N_U$) denote the base change of $g$ (resp.\ $N$) with respect to $N$ (resp.\ $g$).
\end{enumerate}
We denote by $\Coh^r_{\cA, \spe}$ the stack over $k$ such that for each $S\in \Sch/k$, $\Coh^r_{\cA, \spe}(S)$ is the groupoid of special $\cA\boxtimes \cO_S$-modules of rank $r$. The morphism $\cK\mapsto  N(\cK)$ will be denoted by $N: \Coh^r_{\cA, \spe}\to X$.
\end{definition}

\begin{remarks}
\label{remarks:special} 
\rm (a) By Lemma \ref{lemma:maxtoroidal2} it suffices to check condition (ii) for a fixed {\'e}tale covering $\{U_i \to U\}$ and maximal tori $\cT_i$ of $\cA_{U_i}$.

\noi (b) Let $\cK$ be as in \ref{definition:special2} satisfying (i) and assume that $N(\cK): S \to X$ factors through $X-|\Bad(\cA)|$. Then $\cK$ is special of rank $r$ if and only if $\cK$ is a locally free of rank $rd$ as an $\cO_S$-module.

\noi (c) Let $\cA'$ be another locally principal $\cO_X$-order of rank $d^2$ equivalent to $\cA$ and let $\cI$ be an invertible $\cA$-$\cA'$-bimodule. Tensoring with $\cI$ maps $\Coh^r_{\cA, \spe}$ isomorphically to $\Coh^r_{\cA', \spe}$. This follows easily from the fact that, locally on $X$, $\cA$ and $\cA'$ are isomorphic. More generally if $\cA$ are equivalent on some open subscheme $U\subseteq X$ and $\cI$ is a $\cA$-$\cA'$-bimodule which is invertible on $U$ then tensoring with $\cI$ yields an isomorphism $\cdot\otimes_{\cA} \cI: \Coh^r_{\cA, \spe}\times_X U\to \Coh^r_{\cA', \spe}\times_X U$.
\end{remarks}

Except in the appendix we need to consider only the case $r=1$. In the following we investigate the geometric properties of $\Coh_{\cA, \spe} \colon = \Coh^r_{\cA, \spe}$. Recall that a morphism $f:Y\to X$ is said to be semistable if its generic fiber is smooth and for any $y\in Y$ there exists an {\'e}tale neighbourhood $Y'$ of $y$, an open affine neighbourhood $\Spec R$ of $x=f(y)$ and a smooth $X$-morphism 
$Y' \stackrel{g}{\lra} \Spec R[T_1,\ldots, T_r]/(T_1\cdots T_r-\varpi)$ for some $r\ge 1$, where $\varpi$ is a local parameter at $x$. Equivalently, $Y$ is a smooth $k$-scheme, the generic fiber $Y_{\eta}$ is smooth over $F$ and the closed fiber $Y_x$ is a reduced divisor with normal crossings for all $x\in |X|$. Therefore if $f$ is semistable it is flat. 

We have the following simple Lemma whose proof will be left to the reader:

\begin{lemma}
\label{lemma:semistable}
Let $Y_1\stackrel{f}{\lra} Y_2\stackrel{g}{\lra}X$ be morphisms of schemes such that $f$ is smooth and surjective. Then $g$ is semistable if and only if $g\circ f$ is semistable.
\end{lemma}

Let $\cY$ be an algebraic stack over $k$. We will call a morphism $f:\cY\to X$ semistable if there exists a scheme $Y$ and a presentation $P:Y \to \cY$ (i.e.\ $P$ is smooth and surjective) such that $f\circ P : Y \to X$ is semistable. By \ref{lemma:semistable}, if this holds then any presentation $P': Y'\to \cY$ (with $Y'$ a scheme) has this property. In particular if $\cY$ is a scheme the two notions of semistability agree.

Our aim in this section is to prove the following result.

\begin{proposition}
\label{proposition:specialstack}
$\Coh_{\cA,\spe}$ is an algebraic stack over $\bFq$. The morphism $N:\Coh_{\cA,\spe}\to X$ is semistable of relative dimension $-1$. Its restriction to the open subset $X-\Bad(\cA)$ is smooth. Consequently, $\Coh_{\cA,\spe}$ is locally of finite type and smooth over $\bFq$.
\end{proposition}

{\em Proof.} The last assertion follows from (\cite{lau}, 3.2.1). Since the assertion is {\'e}tale local on $X$ we may assume that $X=\Spec R$ is affine with $R$ a principal ideal domain, $|\Bad(\cA)| = \{\fp\}$ and the generic fiber of $\cA$ is $\cong M_d(F)$. By \ref{corollary:gloprincipal}, we may also assume that $e_{\fp}(\cA) = d$. Let $\varpi$ be a generator of $\fp$. Then $\Gamma(\Spec R, \cA)$ is isomorphic to the $R$-subalgebra of $M_d(R)$ of matrices which are upper triangular modulo $\fp$. Hence $\Gamma(\Spec R, \cA)$ can be identified with to the $R$-algebra $R^d\{\Pi\}$ defined by the relations
\[
\Pi(x_1, \ldots, x_d) = (x_2, \ldots, x_d, x_1)\Pi, \qquad \Pi^d = \varpi.
\]
Let $\Coh^{\square}_{\cA,\spe}(S)$ denote the groupoid of pairs $(\cK, \alpha)$ where $\cK\in \Coh_{\cA,\spe}(S)$ and $\alpha: \cO_S^n\to N^*(\cK)$ is an isomorphism. The action of $\Pi$ on $\cK$ yields -- by transport of structure via $\alpha$ -- a map $\cO_S^n\to \cO_S^n$ of the form $(x_1, \ldots, x_d) \mapsto (x_2 a_1, \ldots, x_da_{d-1}, x_1 a_d)$ for some $(a_1, \ldots, a_d)\in \Gamma(S,\cO_S)$ such that $a_1\cdots a_d = N^*(\varpi)$. Thus the assignment $(\cK, \alpha) \mapsto (N, a_1, \ldots, a_d)$ defines an isomorphism
\[
\Coh^{\square}_{\cA,\spe}\cong \Spec R[T_1,\ldots, T_d]/(T_1\cdots T_d-\varpi)
\]
Finally the forgetful morphism $\Coh^{\square}_{\cA,\spe}\to \Coh_{\cA,\spe}$ is a presentation. In fact it induces an isomorphism $\bG^d_m\backslash\Coh^{\square}_{\cA,\spe}\cong\Coh_{\cA,\spe} $. Here the $\bG^d_m$ action on $\Coh^{\square}_{\cA,\spe}$ is defined by the natural $\bG^d_m(S)$-action on the set of isomorphisms $\alpha: \cO_S^n\to N^*(\cK)$.\enddemo

We finish this section with the following criterion for an $\cA\boxtimes \cO_S$-module $\cE$ to be a locally free. 

\begin{lemma} 
\label{lemma:stablyfree}
Let $U\subseteq X$ be a non-empty open subscheme such that $\cE|_{U\times S}$ is a locally free $\cA_U\boxtimes \cO_S$-module. The following conditions are equivalent.

\noi (i) $\cE$ is a locally free $\cA\boxtimes\cO_S$-module of rank $r$.

\noi (ii) For $x\in |X-U|$ and any pair of $k$-morphism $g: S'\to S $ and $x_{S'}:S' \to \Spec k(x) \to X$ the quotient $g^*(\cE)/g^*(\cE)(-\frac{1}{e_x(\cA)}x_{S'})$ is a special $\cA$-module of rank $r$.

\noi (iii) For $x\in |\Bad(\cA)|-U$ and any pair of $k$-morphism $g: S'\to S$ and $x_{S'}:S' \to \Spec k(x) \to X$ the quotient $g^*(\cE)/g^*(\cE)(-\frac{1}{e_x(\cA)}x_{S'})$ is a special $\cA$-module of rank $r$.
\end{lemma}

{\em Proof.} That (i) implies (ii) follows from Lemma \ref{lemma:stablyfree4} and the equivalence of (ii) and (iii) from \ref{remarks:special} (b) above. 

(ii) $\Rightarrow$ (i) We may assume that $S$ is affine, hence that $S$ and of finite type over $k$. For $y\in |X\times S|$ we have to show that $\cE\otimes \cO_{(X\times S),y}$ is a free $(\cA\otimes_{\cO_X} \cO_{(X\times S),y}$-module. It follows from (\cite{lafforgue}, I.2, lemme 4) that we may even replace $S$ by the image $s$ of $y\to X\times S \to S$. Thus it is enough to prove (i) if $S = \Spec \kbar$ is the algebraic closure of $k$. However in this case the assertion follows from \ref{corollary:stablyfree2} and \ref{lemma:stablyfree4}.\enddemo

We have the following generalization of (\cite{lau}, Lemma 1.2.6). 

\begin{lemma} 
\label{lemma:specialres}
Let $0\to \cE'\to \cE \to \cK \to 0$ be a short exact sequence of right $\cA\boxtimes \cO_S$-modules. We assume $\cK$ is coherent as an $\cO_{X\times S}$-module, the map $\Supp(\cK)\hra X\times S \to S$ is an isomorphism and $\cK$ is as an $\cO_S$-module locally free of rank $rd$. We also assume that $\cE$ is a locally free $\cA\boxtimes \cO_S$-module of rank $r$. Then the following conditions are equivalent. 

\noi (i) $\cE'$ is a locally free $\cA\boxtimes \cO_S$-module of rank $r$.

\noi (ii) $\cK$ is special of rank $r$.
\end{lemma}

{\em Proof.} Again by using Lafforgue's Lemma (\cite{lafforgue}, I.2.4) (applied to $\cA$ and maximal tori in $\cA$) it suffices to consider the case where $k$ is algebraically closed and $S = \Spec k$. The assertion follows then from Lemma \ref{lemma:stablyfree5}.\enddemo

\section{The moduli space of $\cA$-elliptic sheaves}
\label{section:aell}
\subsection{$\cA$-elliptic sheaves}
\label{subsection:aell}

In this chapter $X$ denotes a smooth projective geometrically connected curve over the finite field $\bFq$ of characteristic $p$, $F$ the function field of $X$. We also fix a closed point $\infty\in X$. Let $\cA$ be a locally principal $\cO_X$-order of rank $d^2$ and let $A$ be its generic fiber. We make the following
\begin{assumption}
\label{ass:egleichd}
$e_{\infty}(\cA)=d$. 
\end{assumption}

\begin{definition}
\label{definition:aell}
Let $S$ be an $\bFq$-scheme. An $\cA$-elliptic sheaf over $S$ with pole $\infty$ is a triple $E = (\cE,\ioinf, t)$, where $\cE$ is a locally free right $\cA\boxtimes\cO_S$-module of rank $1$, where $\ioinf: S\to X$ is an $\bFq$-morphism with $\ioinf(S) = \{\infty\}$ and where 
\[
t: \ta(\cE(-\frac{1}{d}\ioinf)) \lra \cE 
\]
is an injective $\cA\boxtimes \cO_S$-linear homomorphism such that the following condition holds:

\noi (*) The map $\,\Supp(\Coker(t))\hra X\times S \to S$ is an isomorphism. Considered as a sheaf on $S$, $\cK$ is a locally free $\cO_S$-module of rank $d$.

\noi Hence $\,\Supp(\Coker(t))$ is the image of the graph of a $\bFq$-morphism $\io_0: S \to X$ called the zero (or characteristic) of $E$.

\noi We denote by $\Ell^{\infty}_{\cA}$ the stack over $\bFq$ such that for each $S\in \Sch/k$, $\Ell^{\infty}_{\cA}(S)$ is the category whose objects are $\cA$-elliptic sheaves over $S$ and whose morphisms are isomorphisms between $\cA$-elliptic sheaves.
\end{definition}

For $n\in \frac{1}{d}\bZ$ we define  $\Ell^{\infty}_{\cA,n}$ to be the open and closed substack of $\cA$-elliptic sheaves $E = (\cE,\ioinf, t)$ with $\degcA(\cE) =n$. The functor which maps an $\cA$-elliptic sheaf $E = (\cE,\ioinf, t)$ over $S$ to its zero $\io_0: S\to X$ defines a morphism $\zero: \Ell^{\infty}_{\cA}\to X$ (called the {\it characteristic morphism}). Similarly $E = (\cE,\ioinf, t)\mapsto \ioinf$ defines a morphism $\pole: \Ell^{\infty}_{\cA}\to \Spec k(\infty)$. By Lemma \ref{lemma:specialres}, $\Coker(t)$ is a special $\cA$-module of rank 1. This fact allows us to compare the above condition (*) with the {\it condition sp{\'e}ciale} in (\cite{hausberger}, section 3) (see also \ref{remarks:degbeliebig} (b) below).
It follows that the characteristic morphism factors as 
\begin{equation}
\label{equation:zero}\zero: \Ell^{\infty}_{\cA}\lra \Coh_{\cA, sp} \stackrel{N}{\lra} X
\end{equation} 
We will see in the proof of Theorem \ref{theorem:gutesmodell} below that the first morphism is smooth.

\begin{remarks}
\label{remarks:aellstack}
\rm (a) The concept of an $\cA$-elliptic sheaf is due to Laumon, Rapoport and Stuhler (\cite{lrs}, section 2). The definition given above is different but, as will be explained in the appendix, equivalent to the one given in  (\cite{lrs}, section 2). In fact our Definition \ref{definition:aell} is slightly more general. Their notion corresponds to an $\cA$-elliptic sheaf where
(i) $A$ is a division algebra which is unramified at $\infty$, (ii) $\cA|_{X-\{\infty\}}$ is a maximal order in $A$ and (iii) the zero $\io_0$ is disjoint from $|\Bad(\cA)|$, i.e.\ $\io_0$ factors through $X-|\Bad(\cA)|\hookrightarrow X$ (the latter condition was weakened in \cite{stuhler} and \cite{hausberger} to require only that $\io_0$ factors through $(X-|\Bad(\cA)|)\cup\{\infty\} \cup \{x\in |X|\mid \, \inv_x(A) = \frac{1}{d}\}$). 

\noi (b) Let $\cA$ be the subsheaf of $M_d(\cO_X)$ of matrices which are upper triangular modulo $\infty$. In this case $\Ell^{\infty}_{\cA}$ is isomorphic to the stack $\Ell^{(d)}_X$ of elliptic sheaves of rank $d$ (hence above $X-\{\infty\}$ it is isomorphic to the stack of Drinfeld modules of rank $d$; compare (\cite{stuhler}, section 3)). In fact by Proposition \ref{proposition:comparison} of the appendix we have $\Ell^{\infty}_{\cA}\cong \PEll^{\infty}_{M_d(\cO_X)}$ and the latter is isomorphic to the stack $\Ell^{(d)}_X$ by Morita equivalence. 

\noi (c) If $A$ is a division algebra then $\cA$-elliptic sheaves are special cases of {\it right $\cA$-shtukas} of rank 1 (\cite{lafforgue}, 1.1). Recall that an $\cA$-shtuka of rank 1 is a diagram
\[
\xymatrix@-1.5pc{
\cE \ar[drr]^{j}\\
&&\cE'\\
\ta\cE \ar[urr]^t
}
\]
where $\cE, \cE'$ are locally free right $\cA\boxtimes\cO_S$-modules of rank $1$ and where $j$ and $t$ are injective $\cA\boxtimes \cO_S$-linear homomorphism such that the cokernels of $j$ and $t$ and of the dual morphisms $j^{\vee}$ and $t^{\vee}$ satisfy condition (*) above (actually, it follows from Lemma \ref{lemma:specialres} (compare also the proof of \ref{lemma:1/2gutesmodell} (b) below) that it is enough to require that the cokernels of $j$ and $t$ satisfies (*)). Hence we have $\Coker(j), \Coker(t) \in \Coh_{\cA,\spe}(S)$.
In fact if  $E = (\cE,\ioinf, t)\in \Ell^{\infty}_{\cA}(S)$ is an $\cA$-elliptic sheaf with zero $\io_0: S \to X$ then the diagram
\begin{equation}
\label{eqn:shtuka}
\xymatrix@-1.5pc{
\cE(-\frac{1}{d}\ioinf) \ar[drr]^{j}\\
&&\cE\\
 \ta(\cE(-\frac{1}{d}\ioinf))\ar[urr]^t
}
\end{equation}
is an $\cA$-shtuka with pole $\ioinf$ and zero $\io_0$. Therefore we have a $2$-cartesian square
\begin{equation}
\label{eqn:shtukaaell}
\begin{CD}
\Ell_{\cA}^{\infty} @>>> \Cht^1_{\cA}\\
@VV\pole V@VVV\\
\Spec k(\infty) @>>> \Coh_{\cA, sp}
\end{CD}
\end{equation}
Here the second vertical arrow is given by mapping an $\cA$-shtuka $(\cE, \cE', j, t)$ to $\Coker(j)$. The lower horizontal arrow is defined by $\cA/\cA(-\frac{1}{d}\infty_{k(\infty)})\in \Coh_{\cA, sp}(\Spec k(\infty))$. It is easy to see that it is representable and a closed immersion. Hence the morphism $\Ell_{\cA}^{\infty}\to \Cht^1_{\cA}$ given by (\ref{eqn:shtuka}) is a closed immersion.

\noi (d) One could consider $\cA$-elliptic sheaves more generally for a hereditary $\cO_X$-order $\cA$. However since any hereditary $\cO_X$-order is Morita equivalent to a locally principal $\cO_X$-order we do not obtain new moduli spaces in this way. 

\noi (e) If we consider $\Ell^{\infty}_{\cA}$ as a $k(\infty)$- rather than a $\bFq$-stack we can (and will) drop $\ioinf$ from the definition.
More precisely for $S\in \Sch/k(\infty)$ the objects of $\Ell^{\infty}_{\cA}(S)$ are just pairs $E = (\cE, t)$ such that $(\cE,\ioinf, t)$ is an $\cA$-elliptic sheaf as in \ref{definition:aell} where $\ioinf$ is the composite $S\to \Spec k(\infty)\hra X$.

\noi (f) Define an automorphism of stacks $\theta: \Ell^{\infty}_{\cA} \to \Ell^{\infty}_{\cA}$ by 
\begin{equation}
\label{eqn:shift}
\theta(\cE,\ioinf, t) = (\cE(\frac{1}{d}\ta\ioinf),\ta\ioinf, t(\frac{1}{d}\ta\ioinf))
\end{equation}
where $\ta\ioinf=\ioinf\circ \Frob_S$. We have $\theta(\Ell^{\infty}_{\cA,n}) = \Ell^{\infty}_{\cA,n+\frac{1}{d}}$ for all $n\in \frac{1}{d}\bZ$ and $\theta^{\deg(\infty)}(E) = E\otimes_{\cA} \cA(\frac{1}{d}\infty)$ for all $\cA$-elliptic sheaves $E$. 

\noi (g) Let $\cA'$ be a locally principal $\cO_X$-order which is Morita equivalent to $\cA$ and let $\cL$ be an invertible $\cA$-$\cA'$-bimodule. Then 
\[
E = (\cE,\ioinf,t) \mapsto E\otimes_{\cA} \cL\colon = (\cE\otimes_{\cA} \cL,\ioinf,t\otimes_{\cA}\id_{\cL})
\] 
defines an isomorphism between $\Ell^{\infty}_{\cA}$ and $\Ell^{\infty}_{\cA'}$. If $m = \degcA(\cL)$ then it maps the substack $\Ell^{\infty}_{\cA,n}$ isomorphically onto the substack $\Ell^{\infty}_{\cA',m +n}$. In particular $E\mapsto E\otimes_{\cA} \cL$ defines an action of the abelian group $\Pic(\cA)$ on $\Ell^{\infty}_{\cA}$.
\end{remarks}

We define $\Pic(\cA)[\theta]$ to be the group generated by its subgroup $\Pic(\cA)$ and the element $\theta$ which satisfies the relations
$\theta^{\deg(\infty)} = \cA(\frac{1}{d}\infty)$ and $\theta \cL = \cL \theta$ for all $\cL \in \Pic(\cA)$. Thus $\Pic(\cA)[\theta]$ acts on $\Ell^{\infty}_{\cA}$. The group $\Pic(\cA)[\theta]$ is an extension of $\bZ/\deg(\infty)\bZ\cong \Gal(k(\infty)/\bFq)$ by $\Pic(\cA)$. The map $\deg_{\cA}: \Pic(\cA)\to \frac{1}{d}\bZ$ extends to a homomorphism $\deg_{\cA}:\Pic(\cA)[\theta]\to \frac{1}{d}\bZ$ by defining $\deg_{\cA}(\theta) = \frac{1}{d}$.

\begin{definition}
\label{definition:atkinlehner}
The group of modular automorphisms $\cW(\cA,\infty)$ is defined as the kernel of $\deg_{\cA}:\Pic(\cA)[\theta]\to \frac{1}{d}\bZ$.
\end{definition}

$\cW(\cA,\infty)$ stabilizes the substack $\Ell^{\infty}_{\cA,n}$ for all $n\in \frac{1}{d}\bZ$. There exists a canonical homomorphism 
\begin{equation}
\label{eqn:al1}
\cW(\cA,\infty)\to \Gal(k(\infty)/\bFq)
\end{equation} 
so that $\pole: \Ell^{\infty}_{\cA}\to \Spec k(\infty)$ is $\cW(\cA,\infty)$-equivariant. The kernel of (\ref{eqn:al1}) is $\Pic_0(\cA)$ and the image is of order $\frac{\delta(\cA)\deg(\infty)}{d}$ (thus (\ref{eqn:al1}) is surjective if and only if $\degcA: \Pic(\cA) \to \frac{1}{d}\bZ$ is surjective).

\subsection{Level structure}
\label{subsection:levelaell}

We reformulate now the notion of a level structure on an $\cA$-elliptic sheaf given in (\cite{lrs}, 2.7 and 8.4) in our framework. Let $I= \sum_x \, n_x x$ be an effective divisor on $X$. We assume first that $\infty$ does not divide $I$, i.e.\ $n_{\infty} =0$.

\begin{definition}
\label{definition:aelllevelnoinf}
Suppose that $\infty\not\in |I|$. Let $E = (\cE,\ioinf, t)$ be an $\cA$-elliptic sheaf over an $\bFq$-scheme $S$ with zero $\io_0: S \to X$ disjoint from $I$, i.e.\ $\io_0(S) \cap I = \emptyset$. A level-$I$-structure on $E$ is an isomorphism of right $\cA_I\boxtimes \cO_S$-modules
\[
\al: \cA_I\boxtimes \cO_S \lra \cE|_{I\times S}
\]
compatible with $t$, i.e.\ the diagram
\[
\xymatrix@1{
\tcE|_{I\times S}\ar[rr]^{t|_{I\times S}} & &\cE|_{I\times S}\\
& \cA_I\boxtimes \cO_S\ar@/^/[ul]^{\ta\al}\ar@/_/[ur]_{\al}
}
\]
commutes. 
\end{definition}

We denote by $\Ell^{\infty}_{\cA, I}$ the stack of $\cA$-elliptic sheaves with level $I$-structure and for $n\in \frac{1}{d}\bZ$ by $\Ell^{\infty}_{\cA,I, n}$ the substack of $\cA$-elliptic sheaves of $\cA$-degree $n$ with level $I$-structure. Again we obtain morphisms $\zero: \Ell^{\infty}_{\cA, I}\to X- I$ and $\pole: \Ell^{\infty}_{\cA, I}\to \Spec k(\infty)$. The automorphism (\ref{eqn:shift}) of Remark \ref{remarks:aellstack} (f) extends canonically to an automorphism $\theta:  \Ell^{\infty}_{\cA, I} \to \Ell^{\infty}_{\cA, I}$. The right action of $\Pic(\cA)$ on $\Ell^{\infty}_{\cA}$ lifts to a right action of $\Pic_I(\cA)$ on $\Ell^{\infty}_{\cA, I}$ as follows. If $(\cL, \beta)$ is an invertible $\cA$-$\cA$-bimodule with level-$I$-structure and  $(E,\al)$ an $\cA$-elliptic sheaf with level-$I$-structure $(E,\al)$ over $S$ then we define $(E,\al)\otimes (\cL, \beta) \colon = (E\otimes\cL, \al\bullet\beta)$ with 
\begin{equation}
\label{equation:tensorlevel}
\al\bullet\beta: \cA_I \boxtimes \cO_S \stackrel{\beta\boxtimes \id}{\lra} \cL_I \boxtimes \cO_S = (\cA_I \boxtimes \cO_S)\otimes_{\cA} \cL \boxtimes \cO_S\stackrel{\al\otimes \id}{\lra}(\cE\otimes_{\cA} \cL)|_{I\times S}.
\end{equation}
As before we have $\theta^{\deg(\infty)}(E,\al) = (E,\al)\otimes (\cA(\frac{1}{d}\infty), \id)$.

Suppose now that $|I| = \{\infty\}$, i.e.\ $I = n \infty$ with $n>0$. Let $k(\infty)_d$ be a fixed extension of degree $d$ of $k(\infty)$. According to section \ref{subsection:fqalg} there exists a pair $(\cM_{\infty}, \phi_{\infty})$ consisting of a free right $\cA_{\infty}\otimes_{\bFq} k(\infty)_d$-module $\cM_{\infty}$ of rank $1$ and an isomorphism
\[
\phi_{\infty}:\sig(\cM_{\infty}\fP)\lra \cM_{\infty}
\]
where $\fP$ denotes the maximal invertible two-sided ideal of $\cA_{\infty}\otimes_{\bFq} k(\infty)_d$ corresponding to the inclusion $k(\infty)\hookrightarrow k(\infty)_d$. Let $\cM_I$ denote the sheaf of $\cA_I\otimes_{\bFq} k(\infty)_d$-modules associated to the $\cM_{\infty}/\cM_{\infty} \fp_{\infty}^n$. The map $\phi_{\infty}$ induces an isomorphism
\[
\phi_I:\ta(\cM_I(-\frac{1}{d}\infty_d)) \lra  \cM_I
\]
where $\infty_d$ denotes the morphism $\infty_{k(\infty)_d}:\Spec k(\infty)_d \to \Spec k(\infty)\hookrightarrow X$. 
\begin{definition}
\label{definition:aelllevel}
Let $E = (\cE,\ioinf, t)$ be an $\cA$-elliptic sheaf over an $\bFq$-scheme $S$ with zero $\io_0: S \to X$ disjoint from $I$.

\noi (a) Suppose that $I = n \infty$ with $n>0$. A level-$I$-structure on $E$ consist of a pair $(\la,\al)$ where $\la: S \to \Spec k(\infty)_d$ is an $\bFq$-morphism of schemes which lifts the pole $\ioinf$ and where $\al$ is an $\cA_{I}\boxtimes \cO_S$-linear isomorphism
\[
\al: (\id_I\times \la)^*(\cM_I) \lra \cE|_{I\times S}
\]
such that the diagram
\[
\xymatrix@+4pc{(\ta(\cE(-\frac{1}{d}\ioinf))|_{I\times S}\ar[r]^{t|_{I\times S}}& \cE|_{I\times S}\\
(\id_I\times \la)^*(\ta(\cM_I(-\frac{1}{d}\infty_d)))\ar[r]^{(\id_I\times \la)^*(\phi_I)}\ar[u]_{\ta(\alpha(-\frac{1}{d}\ioinf))} &(\id_I\times \la)^*(\cM_I)\ar[u]_{\alpha}
}
\]
commutes.

\noi (b) Suppose that $I$ is an arbitrary effective divisor on $X$ with $\infty\in |I|$ and write $I = n \infty + I^{\infty} = I_{\infty} + I^{\infty}$ with $n>0$ such that $\infty$ does not divide $I^{\infty}$. A level-$I$-structure on $E$ is a triple $(\al_f,\la,\al_{\infty})$ consisting of a level-$I^{\infty}$-structure $\al_f$ and a level-$I_{\infty}$-structure $(\la,\al_{\infty})$.
\end{definition}

Let $I$ be an effective divisor on $X$ with $\infty\in |I|$. Again we define $\Ell^{\infty}_{\cA, I}$ as the stack of $\cA$-elliptic sheaves with level-$I$-structure $(\cE, t,\al_f,\la,\al_{\infty})$ and denote for $n\in \frac{1}{d}\bZ$ by $\Ell^{\infty}_{\cA,I, n}$ the substack where $\degcA(\cE) =n$. There are canonical morphisms 
\[
\zero: \Ell^{\infty}_{\cA, I}\to X- I, \qquad \pole: \Ell^{\infty}_{\cA, I}\to \Spec k(\infty)_d
\]
(the latter is given by $(E, \al_f,\la,\al_{\infty})\mapsto \la$; it lifts the morphism $\pole: \Ell^{\infty}_{\cA}\to \Spec k(\infty)$). 

\paragraph{Modular automorphisms.} Next we are going to extend the definition of the automorphisms (\ref{eqn:shift}) and define a natural right action of a certain idele class group on $\Ell^{\infty}_{\cA, I}$ (thus lifting the action of $\Pic_I(\cA)$ when $\infty\not\in |I|$). Define $\theta: \Ell^{\infty}_{\cA, I}\to \Ell^{\infty}_{\cA, I}$ by
\begin{equation}
\label{eqn:shift2}
\theta(\cE, \ioinf, t, \al_f,\la,\al_{\infty}) = (\cE(\frac{1}{d}\ta\ioinf), \ta\ioinf, t(\frac{1}{d}\ta\ioinf), \al_f,\ta\la,\al_{\infty}^{\sharp} )
\end{equation}
where $\al_{\infty}^{\sharp}$ is the composite
\[
(\id_I\times \ta\la)^*(\cM_I)\stackrel{\phi_I(\frac{1}{d}\ta\infty_d )}{\lra} (\id_I\times \la)^*(\cM_I( \frac{1}{d}\ta\infty_d)) \stackrel{\al_{\infty}(\frac{1}{d}\ta\ioinf)}{\lra}\cE(\frac{1}{d}\ta\ioinf)|_{I\times S}.
\]
Write $I = n \infty + I^{\infty} = I_{\infty} + I^{\infty}$ with $n>0$ and $\infty\not\in |I^{\infty}|$. Let $\cD_{\infty}$ be a principal order in a central $F_{\infty}$-algebra $D_{\infty}$ of dimension $d^2$ such that 
\[
e(\cD_{\infty}) = e(\cA_{\infty}) \quad \mbox{and} \quad \inv(D_{\infty}) =\inv(A_{\infty}) + \frac{1}{d}. 
\]
We have seen in section \ref{subsection:fqalg} that
 \[
\cD_{\infty}\cong \End_{\cA_{\infty}\otimes_{\bFq} k(\infty)_d}(\cM_{\infty}, \phi_{\infty}). 
\]
We choose an isomorphism (thus making $(\cM_{\infty}, \phi_{\infty})$ into an invertible $\phi$-$\cD_{\infty}$-$\cA_{\infty}$-bimodule of slope $-\frac{1}{d}$). Let 
\[
U_I(\cA^{\infty}\times \cD_{\infty}) \colon = \Ker((\prod_{x\in |X|-\{\infty\}} \, \cA_x^*)\times  \cD_{\infty}^* \to \cA_{I^{\infty}}^*\times (\cD_\infty/\fp_\infty^n\cD_{\infty})^*)
\]
and define
\[
\cC_I(\cA^{\infty}\times \cD_{\infty})\colon = \mbox{$(\prod_{x\in |X|-\{\infty\}}' \, N(\cA_x))\times  N(\cD_{\infty})/U_I(\cA^{\infty}\times \cD_{\infty})F^*$}.
\]
For $a = (a_f, a_{\infty}) = (\{a_x\}_{x\ne \infty},a_{\infty})\in (\prod_{x\in |X|-\{\infty\}}' \, N(\cA_x))\times  N(\cD_{\infty})$ let $\dvsor(a) = (\sum_{x\in |X|-\{\infty\}} v_{\cA_x}(a_x) x) + v_{\cD_{\infty}}(a_\infty) \infty\in \DivcA$. Let $(E,\al_f,\la,\al_{\infty})\break \in \Ell^{\infty}_{\cA, I}(S)$. Left multiplication by $a_f$ on $\cA$ induces a level-$I^{\infty}$-structure $\al_f \cdot a_f$ on $\cE(\dvsor(a))$. More precisely $\al_f \cdot a_f$ is the composite
\[
\xymatrix@1{
\al_f \cdot a_f: \cA_{I^{\infty}}\boxtimes \cO_S \ar[r]^{a_f\cdot} & \cA(\dvsor(a))_{I^{\infty}}\boxtimes \cO_S \ar[rr]^{\al_f(\dvsor(a))} & & \cE(\dvsor(a))|_{I^{\infty}\times S}
}
\]
Similarly left multiplication by $a_\infty$ on $\cM_{\infty}$ yields a level-$I_{\infty}$-structure $\al_{\infty}\cdot a_\infty$ on $E\otimes \cA(\dvsor(a))$. One easily verifies that 
\[
(E,\al_f,\la,\al_{\infty})\cdot a \colon = (E \otimes \cA(\dvsor(a)),\al_f \cdot a_f, \la , \al_{\infty}\cdot a_\infty)
\] 
yields a right $(\prod_{x\in |X|-\{\infty\}}' \, N(\cA_x))\times  N(\cD_{\infty})$-action on $\Ell^{\infty}_{\cA, I}$ and that it factors through $\cC_I(\cA^{\infty}\times \cD_{\infty})$. 

The canonical projection $\cC_I(\cA^{\infty}\times \cD_{\infty})\to \cC_{I^{\infty}}(\cA)$ (given on the $\infty$-factor by $N(\cD_{\infty}) \stackrel{v_{\cD_{\infty}}}{\lra} \frac{1}{d}\bZ \cong N(\cA_{\infty})/\cA_{\infty}^*$) followed by the isomorphism $\cC_{I^{\infty}}(\cA)\to \Pic_{I^{\infty}}(\cA)$ from \ref{corollary:piclevel} yields also a $\cC_I(\cA^{\infty}\times \cD_{\infty})$-action on $\Ell^{\infty}_{\cA, I^{\infty}}$ and one checks that the forgetful morphism of stacks $\Ell^{\infty}_{\cA, I} \lra \Ell^{\infty}_{\cA, I^{\infty}}$ commutes with the $\cC_I(\cA^{\infty}\times \cD_{\infty})$-actions. 

By Remark \ref{remark:phief}, there exists a prime element $\varpi_{\infty}\in\cO_{\infty}$ such that the class $\xi = \xi_{\infty}\in\cC_I(\cA^{\infty}\times \cD_{\infty})$ of the idele $(\{1\}_{x\ne \infty},\varpi_{\infty})$ satisfies $\theta^{d\deg(\infty)}(E) = E \cdot \xi$ for all $E\in \Ell^{\infty}_{\cA, I}(S)$. 

If $\infty$ does not divide the level $I$ we define the group $\Pic_I(\cA)[\theta]$ similar to $\Pic(\cA)[\theta]$ in the last section. $\Pic_I(\cA)[\theta]$ contains $\Pic_I(\cA)$ as a subgroup and the element $\theta$ lies in the center and satisfies the relation $\theta^{\deg(\infty)} = (\cA(\frac{1}{d}\infty), \id) =\colon = \xi_{\infty}$. Let $\deg_{\cA}:\Pic_I(\cA)[\theta]\to \frac{1}{d}\bZ$ be given by $(\cL, \beta)\mapsto \degcA(\cL)$ on $\Pic_I(\cA)$ and $\deg_{\cA}(\theta) = \frac{1}{d}$.

Assume now that $\infty$ divides $I$ and write $I = n \infty + I^{\infty} = I_{\infty} + I^{\infty}$ as above. Define $\cC_I(\cA^{\infty}\times \cD_{\infty})[\theta]$ as the group generated by $\cC_I(\cA^{\infty}\times \cD_{\infty})$ and a central element $\theta$ satisfying the relation $\theta^{d\deg(\infty)}=\xi$. The homomorphism $\cC_I(\cA^{\infty}\times \cD_{\infty})\to \cC_{I^{\infty}}(\cA)\cong \Pic_{I^{\infty}}(\cA)\stackrel{\degcA}{\lra} \frac{1}{d}\bZ$ extends to a homomorphism $\degcA: \cC_I(\cA^{\infty}\times \cD_{\infty})[\theta]\to \frac{1}{d}\bZ$ by setting $\degcA(\theta) = \frac{1}{d}\bZ$. 

\begin{definition}
\label{definition:atkinlehnerlevel}
Let $I$ be an effective divisor on $X$. The group of modular automorphisms $\cW(\cA,I,\infty)$ of $\Ell^{\infty}_{\cA, I}$ is defined as follows:
\[
\cW(\cA,I,\infty)\colon = \left\{ \begin{array}{ll}
        \Ker(\deg_{\cA}:\Pic_I(\cA)[\theta]\to  \frac{1}{d}\bZ) & \mbox{if $\infty\not\in |I|$,}\\
        \Ker(\deg_{\cA}:\cC_I(\cA^{\infty}\times \cD_{\infty})[\theta]\to \frac{1}{d}\bZ) & \mbox{if $\infty\in |I|$.}
        \end{array} \right.
\]
\end{definition}

\begin{remarks}
\label{remarks:modaut}
\rm (a) If $\infty\not\in |I|$ (resp. $\infty\in |I|$) there exists a canonical homomorphism $\cW(\cA,I,\infty)\to \Gal(k(\infty)/\bFq)$ (resp.\ $\cW(\cA,I,\infty)\to \Gal(k(\infty)_d/\bFq)$) with kernel $\Pic_{I,0}(\cA)$ (resp.\ $\cC_I(\cA^{\infty}\times \cD_{\infty})_0$) such that $\Ell^{\infty}_{\cA,I}\to \Spec k(\infty)$ (resp.\ $\Ell^{\infty}_{\cA,I}\to \Spec k(\infty)_d$) is $\cW(\cA,I,\infty)$-equivariant.

\noi (b) For effective divisors $I,J$ on $X$ with $I\le J$ there exists a canonical projection $\cW(\cA,J,\infty)\to \cW(\cA,I,\infty)$ such that the forgetful morphism $\Ell^{\infty}_{\cA, J} \lra \Ell^{\infty}_{\cA, I}$ becomes $\cW(\cA,J,\infty)$-equivariant.

\noi (c) The map $x\mapsto x \theta^{-d\deg_{\cA}(x)}$ induces an isomorphism 
\[
\cC_I(\cA^{\infty}\times \cD_{\infty})/\xi_{\infty}^{\bZ} \quad\cong\quad\cW(\cA,I,\infty).
\]
This fact will be used in section \ref{subsection:uniformization}.
\end{remarks}

We have (compare (\cite{lrs}, 8.10) and (\cite{lafforgue}, I.3.5))

\begin{lemma}
\label{lemma:galoiscover}
Let $I<J$ be effective divisors on $X$. Over $X- J$ the forgetful morphism
\[
\Ell^{\infty}_{\cA, J} \lra \Ell^{\infty}_{\cA, I}
\]
is representable and is a finite, {\'e}tale Galois covering. Its Galois group is $\cong \Ker(\cW(\cA,J,\infty)\to \cW(\cA,I,\infty))$. If $\infty\not\in |J|-|I|$ it is $\cong\Ker(\cA^*_J \to \cA^*_I)$.
\end{lemma}

\begin{corollary}
\label{corollary:finitecover}
Let $\cA'$ be a locally principal $\cO_X$-suborder of $\cA$ with the same generic fiber $A$ and denote by $\iota:Y\hra X$ the reduced closed subscheme with $|Y| = \{x\in |X|\mid\, e_x(\cA') > e_x(\cA')\}$. Note that $\infty\not\in Y$. Let $I$ be an effective divisors disjoint from $Y$ and put $J\colon = I + Y$. Then over $X- J$ the forgetful morphism factors canonically as 
\begin{equation}
\label{equation:levely}
\Ell^{\infty}_{\cA, J} \lra \Ell^{\infty}_{\cA', I}\lra \Ell^{\infty}_{\cA, I}.
\end{equation}
Both maps are representable and finite and {\'e}tale. Moreover the first arrow is Galois.
\end{corollary}

{\em Proof.} Let $\cP \colon = \Image(\cA'|_Y \to\cA|_Y)$. Then the diagram 
\[
\begin{CD}
\cA'@>>> \iota_*(\cP)\\
@VVV@VVV\\
\cA@>>> \iota_*(\cA_Y)
\end{CD}
\]
is cartesian. Here we view $J$ as a closed subscheme of $X$ and denote by $\iota: J \to X$ the inclusion. For $E = (\cE, t, \ioinf, \alpha)$ in $\Ell^{\infty}_{\cA, J}(S)$ we decompose $\alpha$ into a level-$I$-structure $\alpha_I$ and a level-$Y$-structure $\alpha_Y$. Define $\cE'$ by the cartesian square
\[
\xymatrix@+0.5pc{
\cE'\ar[rr]\ar[d]& & \iota_*(\cP)\ar[d]\\
\cE\ar[r] & \iota_*(\cE|_{Y\times S})\ar[r]^{{\alpha_Y}^{-1}}& \iota_*(\cA_Y)
}
\]
Then $\cE'$ is a locally free $\cA'\boxtimes \cO_S$-module of rank 1. The first morphism in (\ref{equation:levely}) is induced by $\cE \mapsto \cE'$ whereas the second by $\cE'\mapsto \cE'\otimes_{\cA'} \cA$. The proof of the remaining assertions are left to the reader.\enddemo

\subsection{The coarse moduli scheme}
\label{subsection:coarsemoduli}

Our aim now is to prove the following theorem. 

\begin{theorem}
\label{theorem:gutesmodell}
(a) $\Ell^{\infty}_{\cA,I}$ is a Deligne-Mumford stack over $\bFq$. It is locally of finite type over $X$. The morphism $\zero: \Ell^{\infty}_{\cA,I}\to X$ is semistable of relative dimension $d-1$. 

\noi (b) The open and closed substack $\Ell^{\infty}_{\cA,I, 0}$ of $\Ell^{\infty}_{\cA,I}$ is of finite type over $X$. It admits a coarse moduli scheme which will be denote by $\El^{\infty}_{\cA,I}$. The structural map $\Ell^{\infty}_{\cA,I,0}\to \El^{\infty}_{\cA,I}$ is an isomorphism if $I\ne 0$. 

\noi (c) The morphism $\zero: \El^{\infty}_{\cA,I}\to X-I$ is quasiprojective and semistable of relative dimension $d-1$. It is smooth over the open subset $X-(\Bad(\cA)\cup I)$. In particular $\El^{\infty}_{\cA,I}$ is a smooth, quasiprojective $\bFq$-scheme.
\end{theorem}

This is known if $A$ is a division algebra or $A=M_d(F)$ and if we restrict $\Ell^{\infty}_{\cA,I}$ to the open subset $X-(\Bad(\cA)\cup I)$ (\cite{lrs}, Theorem 4.1 and \cite{drinfeld1}). The proof of \ref{theorem:gutesmodell} consists essentially of two parts. In the first part one shows that $\Ell^{\infty}_{\cA}\to X$ is a Deligne-Mumford stack and semistable. In the second part one proves that for $I\ne 0$, $\Ell^{\infty}_{\cA,I,n}$ is a quasiprojective scheme by showing that for large $m$ the map $\Ell^{\infty,\stab}_{\cA,mI,n}\to \Ell^{\infty}_{\cA,I,n}$ is surjective. Here $\Ell^{\infty,\stab}_{\cA,I}$ denotes the substack of $\cA$-elliptic sheaves whose underlying vector bundle is {\it $I$-stable}. It is a consequence of a theorem of Seshadri that $\Ell^{\infty,\stab}_{\cA,I}$ is a quasiprojective scheme. The key step in the proof of the surjectivity is to show that $\Ell^{\infty}_{\cA,I,n}$ is of finite type over $\bFq$ (see Lemma \ref{lemma:finitetype} below). Since the proof of both parts are mainly reproduction of the corresponding arguments in (\cite{lrs}, section 4 and 5; compare also (\cite{lafforgue}, Chapitre I), \cite{laumon} and (\cite{lau}, 1.3 and 1.4)) we will be rather brief and elaborate only on those steps were essential modification have to be made.

For the first part we follow (\cite{lau}, 1.2) and work with the factorization (\ref{equation:zero}), i.e.\ we consider $\Ell^{\infty}_{\cA,I}$ mostly over $\Coh_{\cA,\spe}$ rather than over $X$.  

Let $\Inje_{\cA,\spe}$ be the stack over $k(\infty)$ such that for each $S\in \Sch/k(\infty)$, $\Inje_{\cA,\spe}(S)$ is the groupoid of injective morphisms $j:\cE' \to\cE$ of locally free right $\cA\boxtimes\cO_S$-modules of rank $1$ with $\Coker(j)\in \Coh_{\cA,\spe}(S)$. 

\begin{lemma}
\label{lemma:injstack}
(a) The two morphism 
\[
\Inje_{\cA,\spe}\lra \Vect_{\cA,0}^1\times \Coh_{\cA,\spe}
\]
given by $(j:\cE' \to\cE)\mapsto (\cE, \Coker(j))$ and $(j:\cE' \to\cE)\mapsto (\cE', \Coker(j))$ are representable and quasiaffine of finite type and smooth of relative dimension $d$. Consequently $\Inje_{\cA,\spe}$ is algebraic, smooth and of finite type over $\bFq$.

\noi (b) The two morphism
 \[
\Inje_{\cA,\spe}\lra \Vect_X^{d^2}
\]
given by $\cE$ and $\cE'$ are representable and quasiprojective and in particular of finite type.
\end{lemma} 

The proof of (a) for the first morphism is literally the same as (\cite{lau}, 1.3.2). The statement for second morphism can be deduce from that for the first as in (\cite{lau}, 1.3.2). We need to remark only that for a short exact sequence $0 \to \cE' \to \cE \to \cK \to 0$ of right $\cA\otimes \cO_S$-modules with $\cE', \cE\in \Vect^1_{\cA}(S)$ and $\cK\in \Coh_{\cA,\spe}(S)$ the third term of the dual sequence of $\cA^{\opp}\otimes \cO_S$-modules $0 \to \cE^{\vee} \to {\cE'}^{\vee} \to \uExt^1_{\cA\otimes \cO_S}(\cK, \cA\boxtimes \cO_S) \to 0$ lies, by Lemma \ref{lemma:specialres}, in $\Coh_{\cA^{\opp},\spe}(S)$. 

(b) follows from (\cite{lafforgue}, I.2.2 and I.2.8).\enddemo

Consider the following obvious diagram of stacks
\begin{equation}
\label{equation:stackdiagram}
\begin{CD}
\Ell^{\infty}_{\cA} @>>> \Vect^1_{\cA}\otimes_{\bFq} k(\infty)\\
@VVV @VVV\\
\Inje_{\cA,\spe} @>>> (\Vect^1_{\cA}\times \Vect^1_{\cA})\otimes_{\bFq} k(\infty)\\
@VVV\\
\Coh_{\cA, \spe}
\end{CD}
\end{equation}
where the right vertical arrow in the ($2$-cartesian) square is the graph of the endomorphism $\Frob\circ \theta^{-1}: \Vect^1_{\cA}\otimes_{\bFq}k(\infty) \to \Vect^1_{\cA}\otimes_{\bFq}k(\infty)$ (for the definition of $\theta^{-1}$ compare \ref{remarks:aellstack} (f); if $\deg(\infty) = 1$ it is given by $\cE \mapsto \cE(-\frac{1}{d}\infty)$). By (\cite{lafforgue}, I.2.5) the stack $\Vect^1_{\cA}$ is algebraic, locally of finite type and smooth over $\bFq$. Together with Proposition \ref{proposition:specialstack}, Lemma \ref{lemma:galoiscover} and Lemma \ref{lemma:injstack} above the same argument as in (\cite{lrs}, section 4; see also (\cite{lafforgue}, I.2.5) and (\cite{lau}, I.3.5)) imply part (a) of

\begin{lemma}
\label{lemma:1/2gutesmodell}
(a) Let $I$ be an effective divisor on $X$. The morphism $\Ell^{\infty}_{\cA, I}\break \to \Coh_{\cA,\spe}$ is algebraic, locally of finite type and smooth of relative dimension $d$. The morphism $\zero: \Ell^{\infty}_{\cA,I}\to X$ is semistable of relative dimension $d-1$.

\noi (b) $\Ell^{\infty}_{\cA}$ is a Deligne-Mumford stack, locally of finite type and smooth over $\bFq$. Moreover if $I\ne 0$ then $\Ell^{\infty}_{\cA}$ is isomorphic to an algebraic space.
\end{lemma}

{\em Proof of (b).} Everything is clear if we replace ``Deligne-Mumford'' by ``algebraic''. To prove that $\Ell^{\infty}_{\cA}$ is indeed a Deligne-Mumford stack we use (\cite{lmb}, 8.1).  If we replace the lower vertical map in (\ref{equation:stackdiagram}) by $\Inje_{\cA,\spe}\to \Spec k(\infty)$ then, by Lemma \ref{lemma:injstack} and (\cite{laumon}, Lemma on p.\ 60), the diagonal morphism $\Ell^{\infty}_{\cA}\to \Ell^{\infty}_{\cA}\times_{\bFq} \Ell^{\infty}_{\cA}$ is unramified. Note that for $E\in \Ell^{\infty}_{\cA,I}(S)$ we have $\Aut(E) = \bFq^*$ if $I =0$ or $=1$ otherwise. Hence the last assertion follows from (\cite{lmb}, 8.1.1).\enddemo

\begin{lemma}
\label{lemma:finitetype}
For $n\in \frac{1}{d}\bZ$ and $I>0$ the stack $\Ell^{\infty}_{\cA,I, n}$ is isomorphic to a scheme of finite type over $\bFq$. 
\end{lemma}

{\em Proof.} In case where $A$ is a division algebra this follows from a result of Lafforgue (\cite{laumon}, 4.2) by applying Remark \ref{remarks:aellstack} (c) (compare also \cite{lrs}, 5.1). The general case will be reduced to this case by using property (ii) of the morphism (\ref{equation:tensorellsupersing}) of section \ref{subsection:main}. 

In the case where $A$ is an arbitrary central simple $F$-algebra we first remark that we can pass freely to some finite Galois covering. More precisely let $\cX\to \Ell^{\infty}_{\cA,I, n}$ be morphism which is representable and finite and an {\'e}tale Galois covering with Galois group $G$. If $\cX$ is a scheme of finite type then by \ref{lemma:1/2gutesmodell} (b), $\Ell^{\infty}_{\cA,I, n}$ is isomorphic to the quotient scheme $G\backslash\cX$. Thus we may extend the base field $\bFq$ and, by Lemma \ref{lemma:galoiscover}, also increase or decrease the level. Hence we may assume that $|I|\not\subset|\Bad(\cA)|$. Let $x\in |I|-|\Bad(\cA)|$ and let $\cA'$ be a locally principal $\cO_X$-suborder of $\cA$ with $\cA'|_{X-\{x\}} = \cA|_{X-\{x\}}$ and $e_x(\cA') =d$. Write $I = I' + mx$ with $x\not\in |I'|$. By \ref{corollary:finitecover} it is enough to show that $\Ell^{\infty}_{\cA',I', n}$ is a scheme of finite type. Thus we need to prove the assertion only under the following additional hypothesis:

\noi {\it There exists a point $\fp\in|X|-(\{ \infty\} \cup |I|)$ such that $e_{\fp}(\cA) = d$ and $\inv_{\fp}(A) = 0$.} 

Let $\cB$ be a locally principal $\cO_X$-order of rank $d^2$ with $\Bad(\cB) = \Bad(\cA)$ and such that for the generic fiber $B$ of $\cB$ we have $\inv_{\infty}(B)= \inv_{\infty}(A) + \frac{1}{d}$, $\inv_{\fp}(B)= - \frac{1}{d}$ and $\inv_x(B) = \inv_x(A)$ for all $x\in |X|-\{\fp, \infty\}$. In particular $B$ is a division algebra. The following fact will be proved in section \ref{subsection:main}: 

{\it There exists a finite extension $\bF/\bFq$ and embeddings $k(\infty)\hra \bF, k(\fp)\hra \bF$ such that
$\Ell^{\infty}_{\cA,I, n}\otimes_{k(\infty)} \bF \cong \Ell^{\fp}_{\cB,I, n}\otimes_{k(\fp)} \bF$.} 

\noi Since, by Lafforgue's Theorem (\cite{laumon}, 4.2), the stack $\Ell^{\fp}_{\cB,I, n}$ is a scheme of finite type the assertion follows.\enddemo

To finish the proof of \ref{theorem:gutesmodell} in the case $I\ne 0$ it remains to show that $\El^{\infty}_{\cA,I}$ is quasiprojective. The proof is the same as in (\cite{lrs}, section 5). However for the sake of completeness we reproduce it briefly. Firstly, by Lemma \ref{lemma:galoiscover} we are free to enlarge or shrink $I$. Thus we can assume that $\infty\not\in |I|$. By \ref{lemma:injstack} (b) the morphism $\El^{\infty}_{\cA, I} \to \Vect^{d^2}_{I,s}$ given by $E = (\cE, t, \ioinf, \alpha) \mapsto (\cE, \alpha)$ is representable and quasiprojective. Here $s \colon = \deg(\cA)$ and $\Vect^{d^2}_{I,s}$ denotes the stack of vector bundles of degree $s$ with level-$I$-structure. The open substack $\Vect^{d^2,\stab}_{I,s}$ of $I$-stable vector bundles -- hence also its pre-image $\El^{\infty,\stab}_{\cA, I,0}$ -- is a smooth quasiprojective scheme (\cite{lrs}, 4.3). For $n\in \bN$ let $G_n$ be the Galois group of $\El^{\infty}_{\cA, nI,0}\to \El^{\infty}_{\cA, I,0}$. Then $
 \El^{\infty, \stab}_{\cA, nI}/G_n$ is an increasing sequence of open quasiprojective subschemes of $\El^{\infty}_{\cA, I}$ which cover $\El^{\infty}_{\cA, I}$. By \ref{lemma:finitetype} the sequence becomes stationary, i.e.\ $\El^{\infty}_{\cA, I}=\El^{\infty, \stab}_{\cA, nI}/G_n$ for some $n$.

Finally let us consider the case $I=0$. Choose an auxiliary level $J$, $J>0$ with $\infty\not\in |J|$. By \ref{lemma:galoiscover} above $X-J$ we have 
\[
\Ell^{\infty}_{\cA,0}\quad \cong\quad (\cA_J)^*\backslash \El^{\infty}_{\cA, J}
\]
The subgroup $\bFq^*$ of $\cA_I^*$ acts trivially whereas the quotient $(\cA_I)^*/\bFq^*$ acts freely on $\El^{\infty}_{\cA, J}$. Hence above $X-J$, the quotient $((\cA_I)^*/\bFq^*)\backslash \Ell^{\infty}_{\cA, J}$ is the coarse moduli scheme of $\Ell^{\infty}_{\cA,0}$. For varying $J$ these coarse moduli schemes glue together and yield a coarse moduli scheme $\El^{\infty}_{\cA}$ of $\Ell^{\infty}_{\cA,0}$. Moreover $\Ell^{\infty}_{\cA,0}$ is isomorphic to the quotient stack $\bFq^*\backslash \El^{\infty}_{\cA}$. This completes the proof of \ref{theorem:gutesmodell}.

\begin{remarks}
\label{remarks:modulispace}
(a) Note that by \ref{remarks:aellstack} (f) we could have defined $\El^{\infty}_{\cA}$ also as the coarse moduli scheme of the quotient $\Ell^{\infty}_{\cA}/\theta^\bZ$ or of $\Ell^{\infty}_{\cA,n}$ for any $n\in \frac{1}{d}\bZ$.

\noi (b) Let $I\hookrightarrow X$ be a reduced closed subscheme with $\infty\not\in I$ and let $\cA$ be the subsheaf of $M_2(\cO_X)$ of matrices which are upper triangular modulo $I$. Then by using \ref{corollary:finitecover} and \ref{remarks:aellstack} (b) it is easy to see that the $\El^{\infty}_{\cA}$ is isomorphic to the (open) Drinfeld modular curve $Y_0(I)= Y^{\infty}_0(I)$.

\noi (c) If $A$ is a central division algebra which is unramified at $\infty$ and $\cA|_{X-\{\infty\}}$ is a maximal order in $A$ then $\zero: \El^{\infty}_{\cA}\to X$ is proper (see \cite{lrs}, Theorem 6.1 and \cite{hausberger}, 6.4). In the general case this is not true anymore even if $A$ is a division algebra. In fact if $d=2$ and $A$ is ramified only at $\infty$ and at $\fp\in |X|$ and if $\cA$ is a maximal $\cO_X$-order in $A$ then we will show in section \ref{subsection:main} that $\El^{\fp}_{\cA}$ is a twist of the affine curve $Y^{\infty}_0(\fp)\to X$.
\end{remarks}

\subsection{Invertible Frobenius bimodules}
\label{subsection:torsor}

We consider now two locally principal $\cO_X$-orders $\cA$ and $\cB$, both of rank $d^2$ with $\Bad(\cA) = \Bad(\cB)$ and assume that $e(\cA) =d = e(\cB)$. We denote by $A$ and $B$ the generic fibers of $\cA$ and $\cB$ respectively. Let $D = \sum_{x\in |X|}\, m_x x \in \Div(\cA)$ be a divisor such that $\sum_{x\in |X|}\, m_x = 0$. We consider the following moduli problem associated to $\cA, \cB$ and $D$.

\begin{definition}
\label{definition:frobbiglobal}
Let $S$ be an $\bFq$-scheme. An invertible Frobenius $\cA$-$\cB$-bimodule (or $F$-$\cA$-$\cB$-bimodule for short) of slope $D$ over $S$ is a tuple $L= (\cL, (x_S)_{x\in |D|},t)$ consisting of the following data:

\noi -- $\cL$ is an invertible $\cA\boxtimes \cO_S$-$\cB\boxtimes \cO_S$-bimodule which is locally free of rank 1 as a left $\cA\boxtimes \cO_S$- and right $\cB\boxtimes \cO_S$-module,

\noi -- For all $x\in |D|$, $x_S: S\to X$ is a morphism in $\Sch/\bFq$ which factors through $x\to X$ (the morphisms $x_S$ are called the poles of $L$),

\noi -- $t$ is an isomorphism of bimodule
\[
t :\ta(\cL(D_S))\lra \cL
\]
where $D_S\colon = \sum_{x\in |D|}\, m_x x_S$. 

We denote by $\SEl^D_{\cA,\cB}$ the stack over $\bFq$ of invertible $F$-$\cA$-$\cB$-bimodules of slope $D$.
\end{definition}

Note that $\deg_{\cA}(\cL) = \deg_{\cB}(\cL)$. Note also that $(\cA\boxtimes \cO_S)(D_S)\otimes_{\cA} \cL= \cL\otimes_{\cB}(\cB\boxtimes \cO_S)(D_S)$. Thus the notion $\cL(D_S)$ is unambiguous. 

There are canonical morphisms $\SEl^D_{\cA,\cB}\to \Spec k(x)$ for all $x\in |D|$. For $n\in \frac{1}{d}\bZ$ we let $\SEl^D_{\cA,\cB,n}$ be the substack of $\SEl^D_{\cA,\cB}$ of $F$-$\cA$-$\cB$-bimodules of fixed $\cA$-degree $n$. There is a canonical left $\Pic(\cA)$- and right $\Pic(\cB)$-action on $\SEl^D_{\cA,\cB}$ compatible with $\degcA$ (in fact the left and right action are the same if we identify the two groups under the canonical isomorphism $\Pic(\cA)\cong\DivcA/F^*\cong \Pic(\cB)$). 

For $x\in |D|$ we define an automorphism $\theta_x: \SEl^D_{\cA,\cB}\to \SEl^D_{\cA,\cB}$ by 
\begin{equation}
\label{eqn:shiftx1}
\theta_x(\cL, (x'_S)_{x'\in |D|},t) = (\cL(-m_x \tax_S), \tax_S, (x'_S)_{x'\in |D|, x'\ne x},t(-m_x \tax_S)).
\end{equation}
The automorphisms $\theta_x$ for different $x\in |D|$ commute with each other and with the $\Pic(\cA)$- and $\Pic(\cB)$-action. We have $\theta_x( \SEl^D_{\cA,\cB,n}) = \SEl^D_{\cA,\cB,n-m_x}$ for all $n\in \frac{1}{d}\bZ$ and $\theta_x^{\deg(x)}(L) = \cA(-m_x x)\otimes L  = L \otimes\cB(-m_x x)$ for all $L\in \SEl^D_{\cA,\cB}(S)$. Moreover if $|D| = \{x_1, \ldots, x_m\}$ and if we put $\Theta_D\colon = \theta_{x_1}\circ\ldots\circ \theta_{x_m}$ then $\Theta_D(L) = \Frob_S^*(L)$ for all $S\in \Sch/\bFq$ and $L\in \SEl^D_{\cA,\cB}(S)$. Hence $\Theta_D = \Frob_{\SEl^D_{\cA,\cB}}$.

\paragraph{Level structure.}

 Let $L= (\cL, (x_S)_{x\in |D|}, t)\in \SEl^D_{\cA,\cB}(S)$ and let $I$ be an effective divisor on $X$. To define a level-$I$-structure on $L$ we view $\cL$ as a right $\cB$-module and proceed as in section \ref{subsection:levelaell}. Assume first that $|I|\cap |D| = \emptyset$. Then a level-$I$-structure on $L$ is an isomorphism of right $\cB_I\boxtimes \cO_S$-modules $\beta: \cB_I\boxtimes \cO_S \lra \cL|_{I\times S}$ such that the diagram
\[
\xymatrix@1{
\ta\cL|_{I\times S}\ar[rr]^{t|_{I\times S}} & &\cL|_{I\times S}\\
& \cA_I\boxtimes \cO_S\ar@/^/[ul]^{\ta\beta}\ar@/_/[ur]_{\al}
}
\]
commutes. 

Next assume that $I = n x$ with $n>0$ for some $x\in |D|$. Put $e=e_x(\cA), m =m_x$ and let $k(x)_e$ be an extension of degree $e$ of $k(x)$. If $m> 0$ we denote by $M=(\cM_{x}, \phi_{x})$ a fixed invertible $\phi$-$\cB_x$-$\cA_x$-bimodule of slope $-m$ over $\cO_{x}\otimes_{\bFq} k(x)_e$. In case $m<0$, $M=(\cM_{x}, \phi_{x})$ denotes a $\phi$-$\cA_x$-$\cB_x$-bimodule of slope $m$ over $\cO_{x}\otimes_{\bFq} k(x)_e$. Thus if $m>0$ (resp.\ $m<0$) then $\phi_{x}$ is an isomorphism
\[
\phi_x:\sig(\cM_x\fP^{me})\lra \cM_x \quad (\mbox{resp.}\,\,\phi_x:\sig(\fP^{-me}\cM_x)\lra \cM_x)
\]
where $\fP$ denotes the maximal ideal of $\cA_x\otimes_{\bFq} k(x)_e$ corresponding to the inclusion $k(x)\hookrightarrow k(x)_e$. If $m>0$ (resp.\ $m<0$) then, as in \ref{subsection:levelaell}, the pair $(\cM_{x}, \phi_{x})$ induces a pair $(\cM_I, \phi_I)$ where $\cM_I$ is an $\cB_I$-$\cA_I$-bimodule (resp.\ $\cA_I$-$\cB_I$-bimodule) and $\phi_I$ is an isomorphism $\phi_I:\ta(\cM_I(-m x_e)) \lra  \cM_I$ (resp.\ $\phi_I:\ta(\cM_I(m x_e)) \lra  \cM_I$). Here $x_e$ denotes the map $\Spec k(x)_e \to x \hookrightarrow X$. A {\it level-$I$-structure on $L$} is a pair $(\mu,\beta)$ consisting of an $\bFq$-morphism $\mu: S \to \Spec k(x)_e$ which lifts $x_S$ and an isomorphism of right $\cB_{I}\boxtimes \cO_S$-modules $\beta$. If $m<0$, then $\beta$ is a map
\begin{equation}
\label{eqn:levelfrob1}
\beta: (\id_I\times \mu)^*(\cM_I) \lra \cL|_{I\times S}
\end{equation}
such that
\[
\xymatrix@1{\ta(\cL(D))|_{I\times S}\ar[rr]^{\phi_I}& &\cL|_{I\times S}\\
(\id_I\times \mu)^*(\ta(\cM_I(m x_e)) \ar[u]_{\beta}\ar[rr]^{\phi_I} & &(\id_I\times \mu)^*(\cM_I)\ar[u]_{\beta}
}
\]
commutes. If $m>0$ then
\begin{equation}
\label{eqn:levelfrob2}
\beta: \cB_I\boxtimes \cO_S\lra (\id_I\times \mu)^*(\cM_I) \otimes_{\cA_I\boxtimes \cO_S}\cL|_{I\times S}
\end{equation}
so that 
\[
\xymatrix@1{
\ta((\id_I\times \mu)^*(\cM_I) \otimes_{\cA}\cL|_{I\times S})\ar[rr]^{\phi_I\otimes t|_{I\times S}} & & (\id_I\times \mu)^*(\cM_I) \otimes_{\cA}\cL|_{I\times S}\\
& \cB_I\boxtimes \cO_S\ar@/^/[ul]^{\ta\beta}\ar@/_/[ur]_{\beta}
}
\]
commutes. 

For an arbitrary effective divisor $I$ on $X$ we write $I= I_0 + \sum_{x\in |I|\cap |D|}\, n_x x\break = I_0 + \sum_{x\in |I|\cap |D|}\, I_{x}$ with $|I_0| \cap |D|=\emptyset$ and $I_x = n_x x$ for $x\in |I|\cap |D|$. Then a level-$I$-structure on $L$ is a tuple $(\beta_0,(\mu_x,\beta_x)_{x\in |I|\cap |D|})$ consisting of a level-$I_0$-structure $\beta_0$ and level-$I_{x}$-structures $(\mu_{x},\beta_{x})$ for all $x\in |I|\cap |D|$. This yields stacks $\SEl^D_{\cA,\cB,I}, \SEl^D_{\cA,\cB,I,n}$ equipped with forgetful morphisms 
\[
\SEl^D_{\cA,\cB,I}\to \SEl^D_{\cA,\cB},\qquad  \SEl^D_{\cA,\cB,I}\to \Spec k(x)_{e_x(\cA)}\quad\mbox{for $x\in |I|\cap |D|$} 
\]
(the latter lifts the morphism $\SEl^D_{\cA,\cB}\to \Spec k(x)$). 

\paragraph{Modular automorphisms.}\label{paragraph:theta} Let $T\colon = \{x\in |D|\mid \, m_x>0\}$. If $|I|\cap T = \emptyset$ then there is a canonical left $\Pic_I(\cA)$-action on $\SEl^D_{\cA,\cB,I}$ lifting the $\Pic(\cA)$-action on $\SEl^D_{\cA,\cB}$. We want to extend this to a natural left action of an idele class group $\cC_I(\cA^T\times \cB_T)$ on $\SEl^D_{\cA,\cB,I}$ for arbitrary $I$ (similarly to the right action of $\cC_I(\cA^{\infty}\times \cD_{\infty})$-action on $\Ell^{\infty}_{\cA, I^{\infty}}$ defined in \ref{subsection:levelaell}). Write $I = I^T + I_T$ with $|I^T| \cap T = \emptyset$ and $|I_T|\subseteq T$. Put
\[
U_I(\cA^T\times \cB_T) \colon = \Ker(\prod_{x\in |X|-T} \, \cA_x^*\times  \prod_{x\in T}\cB_x^* \to \cA_{I^{T}}^*\times \cB_{I_T}^*)
\]
and define
\[
\cC_I(\cA^T\times \cB_T)\colon = \mbox{$\prod_{x\in |X|-T}'\, N(\cA_x)\times\prod_{x\in T}\, N(\cB_x)/U_I(\cA^T\times \cB_T)F^*$}.
\]
There is a canonical epimorphism $\cC_I(\cA^T\times \cB_T) \lra \Pic(\cA)$ given on the class $[g]$ represented by $g= (\{a_x\}_{x\in |X|-T},\{b_x\}_{x\in T})\in \prod_{x\in |X|-T}'\, N(\cA_x)\times\prod_{x\in T}\, N(\cB_x)$ by $\cA(\dvsor(g))$ where
\[
\dvsor(g) = \sum_{x\in |X|-T} v_{\cA_x}(a_x) x\,\, + \,\,\sum_{x\in T} v_{\cB_x}(b_x) x.
\]
The kernel of the composition $\degcA:\cC_I(\cA^T\times \cB_T) \to \Pic(\cA) \stackrel{\degcA}{\lra}\bQ$ will be denoted by $\cC_I(\cA^T\times \cB_T)_0$. 

Let $g=(a^T,b_T)= (\{a_x\}_{x\not\in T},\{b_x\}_{x\in T})\in \prod_{x\in |X|-T}'\, N(\cA_x)\times\prod_{x\in T}\, N(\cB_x)$ and $L = (L,\beta_0,(\mu_x,\beta_x)_{x\in |I|\cap |D|})\in \SEl^D_{\cA,\cB,I}(S)$. Left multiplication by $a^T$ on the target of $\beta_0$ and $\beta_x$ for $x\in |I^T|$ (respectively by $b_T$ on the source of $\beta_x$ for $x\in |I_T|$) yields a level-$I$-structure on $\cA(\dvsor(g)) \otimes L$. This defines a left action of $\prod_{x\in |X|-T}'\, N(\cA_x))\times\prod_{x\in T}\, N(\cB_x)$ on $\SEl^D_{\cA,\cB,I}$ which factors through $\cC_I(\cA^T\times \cB_T)$. 

\begin{abschnitt}
\label{abschnitt:theta} \rm Similar to (\ref{eqn:shift2}), for $x\in |D|$ there exists a canonical lift of (\ref{eqn:shiftx1}) to an automorphism $\theta_x: \SEl^D_{\cA,\cB,I} \to \SEl^D_{\cA,\cB,I}$ having the following properties:
\begin{itemize}
\label{itemize:theta}
\item[(i)] The following diagram commutes
\[
\begin{CD}
\SEl^D_{\cA,\cB,I} @> \theta_x >> \SEl^D_{\cA,\cB,I}\\
@VVV@VVV\\
\Spec k(x)_{*}@> \Frob_q >>\Spec k(x)_{*}
\end{CD}
\]
where $* = 1$ or $* =e_x(\cA)$ depending on whether $x\not\in |I|$ or $x\in |I|$.
\item[(ii)] For $n\in \frac{1}{d}\bZ$ we have $\theta_x( \SEl^D_{\cA,\cB,I ,n}) = \SEl^D_{\cA,\cB,I,n-m_x}$.
\item[(iii)] The automorphisms $\theta_x$ for different $x\in |D|$ commute with each other and with the $\cC_I(\cA^T\times \cB_T)$-action.
\item[(iv)] For $x\in |D|$ there exists $\xi_x\in \cC_I(\cA^T\times \cB_T)$ such that $\theta_x^{\deg(x)}(L) = \xi_x L$ (resp.\ $\theta_x^{e_x(\cA)\deg(x)}(L) = \xi_x L$) for all $L\in \SEl^D_{\cA,\cB,I}(S)$.
\item[(v)] If $|D| = \{x_1, \ldots, x_m\}$ put $\Theta_D\colon = \theta_{x_1}\circ\ldots\circ \theta_{x_m}$. Then $\Theta_D = \Frob_{\SEl^D_{\cA,\cB,I}}$.
\end{itemize}
\end{abschnitt}

Let $\cC_I(\cA^T\times \cB_T)[\,\theta_x, x\in |D|\,]$ denote the following group. It contains $\cC_I(\cA^T\times \cB_T)$ as a subgroup and it contains the set $\{\theta_x\mid \, x\in |D|\}$ and it is generated by the union of both sets. The elements $\theta_x$ for $x\in |D|$ are all central and satisfy the relations (iv) above. We extend the homomorphism $\deg_{\cA}: \cC_I(\cA^T\times \cB_T) \to \frac{1}{d}\bZ$ to $\cC_I(\cA^T\times \cB_T)[\,\theta_x, x\in |D|\,]$ by setting $\degcA(\theta_x) = -m_x$.

\begin{definition}
\label{definition:atkinlehnerfrob}
Suppose that $\SEl^D_{\cA,\cB,I ,0}\ne \emptyset$. The group of modular automorphisms of $\cW(\cA,\cB,I,D)$ of $\SEl^D_{\cA,\cB,I ,0}$ is defined as follows:
\[
\cW(\cA,\cB,I,D) = \Ker(\cC_I(\cA^T\times \cB_T)[\,\theta_x, x\in |D|\,] \lra  \frac{1}{d}\bZ).
\]
\end{definition}

\begin{remark}
\label{remark:atkinlehnerfrob}
\rm Assume that $\SEl^D_{\cA,\cB,I ,0}\ne \emptyset$. For all $x\in |D|$ there exists canonical homomorphisms $\cW(\cA,\cB,I,D)\to \Gal(k(x)_*/\bFq)$ where $* = \emptyset$ or $* =e_x(\cA)$ depending on whether $x\not\in |I|$ or $x\in |I|$. It is surjective since $\Theta_D$ is mapped to $\Frob_{k(x)_*}$ by property (v) above. The kernel of the homomorphism
\begin{equation}
\label{eqn:atkinlehnerfrob}
\cW(\cA,\cB,I,D) \lra \prod_{x\in |D|-|I|} \Gal(k(x)/\bFq) \times \prod_{x\in |D|\cap|I|} \Gal(k(x)_{e_x(\cA)}/\bFq)
\end{equation}
is $\cC_I(\cA^T\times \cB_T)_0$. 

It is easy to see that (\ref{eqn:atkinlehnerfrob}) is surjective provided that $\delta(\cA) = d$ (i.e.\  $\degcA: \Pic(\cA) \to \frac{1}{d}\bZ$ is surjective). 
\end{remark}

\paragraph{Tensor product and Inverse.}
 
There is also a notion of a tensor product of invertible Frobenius bimodules and of an inverse. These constructions are needed in the proof of Proposition \ref{proposition:torseur1} below. Let $\cC$ be a third locally principal $\cO_X$-order of rank $d^2$ with $\Bad(\cC) = \Bad(\cA)$. Let $D_1 = \sum_{x\in |X|}\, m^{(1)}_x x, D_2 = \sum_{x\in |X|}\, m^{(2)}_x x\in \Div(\cA)$ with $\sum_{x\in |X|}\, m^{(i)}_x = 0$ for $i=1,2$. Let $Y = \Spec \bFq$ if $|D_1|\cap |D_2|= \emptyset$ or $ Y = \Spec(\bigotimes_{x\in|D_1|\cap |D_2|} k(x))$ otherwise. We view $\SEl^{D_1}_{\cA,\cB}$ and $\SEl^{D_2}_{\cB,\cC}$ as stacks over $Y$. Let $S\in \Sch/Y$ and let $L= (\cL, (x_S)_{x\in |D|},t)\in \SEl^{D_1}_{\cA,\cB}(S)$, $M= (\cM, (x_S)_{x\in |D_2|},t)\in \SEl^{D_2}_{\cB,\cC}(S)$ (hence for $x\in |D_1|\cap |D_2|$, the morphisms $x_S$ for $L$ and $M$ agree and are equal to the canonical morphism $S \to \Spec k(x)\to X$). Define 
\[
L\otimes M = (\cL\otimes \cM, (x_S)_{x\in |D_1+ D_2|}, t \otimes_{\cB} t)\in \SEl^{D_1+ D_2}_{\cA,\cC}(S).
\]
Thus we get a morphism of stacks 
\begin{equation}
\label{eqn:tensfrobinv}
\otimes: \SEl^{D_1}_{\cA,\cB}\times_Y \SEl^{D_2}_{\cB,\cC} \lra \SEl^{D_1+ D_2}_{\cA,\cC}
\end{equation}
which is compatible with degrees.  

The inverse $L^{-1}$ of $L= (\cL, (x_S)_{x\in |D|},t)\in \SEl^{D}_{\cA,\cB}(S)$ is defined as\begin{equation}
\label{eqn:inversfrobinv}
L^{-1} = (\cL^{\vee}, (x_S)_{x\in |D|}, (t^{\vee})^{-1})\in \SEl^{-D}_{\cB,\cA}(S).
\end{equation}
We leave it to the reader to extend the Definition (\ref{eqn:tensfrobinv}) and (\ref{eqn:inversfrobinv}) to invertible Frobenius bimodules with level-$I$-structure (see also the next section where the tensor product of an $\cA$-elliptic sheaf with level-$I$-structure with a Frobenius bimodule with level-$I$-structure is defined).

\paragraph{Moduli spaces.} Let $D = \sum_{x\in |X|}\, m_x x \in \Div(\cA)$, $D\ne 0$ be such that $\sum_{x\in |X|}\, m_x = 0$ and let $I\in \Div(X)$ with $I\ge 0$. Our aim is to prove the following result.

\begin{proposition}
\label{proposition:torseur1}
(a) $\SEl^D_{\cA,\cB,I}\ne \emptyset$ if and only if 
\begin{equation}
\label{eqn:frobbiglobal1}
\sum_{x\in |X|}\, \inv_x(B) x = (\sum_{x\in |X|}\, \inv_x(A) x) + D \mod \Div(X).
\end{equation}

\noi (b) $\SEl^D_{\cA,\cB,I}$ is a Deligne-Mumford stack which is {\'e}tale over $\bFq$. The open and closed substack $\SEl^D_{\cA,\cB,I,n}$ is finite over $\bFq$ for all $n\in \frac{1}{d}\bZ$.

\noi (c) $\SEl^D_{\cA,\cB,I,0}$ admits a coarse moduli space $\SE^D_{\cA,\cB,I}$. The structural morphism  $\SEl^D_{\cA,\cB,I,0}\to \SE^D_{\cA,\cB,I}$ is an isomorphisms if $I\ne 0$.

\noi (d) Suppose that $\SEl^D_{\cA,\cB,I,0}\ne \emptyset$. Then $\SE^D_{\cA,\cB,I}\to \Spec \bFq$ is a $\cW(\cA, \cB, I, D)$-torsor. In particular  $\SE^D_{\cA,\cB,I}$ is a finite, {\'e}tale $\bFq$-scheme.
\end{proposition}

We begin with the proof of (a). Since $\SEl^D_{\cA,\cB,I}$ is locally of finite presentation it suffices to show that $\SEl^D_{\cA,\cB,I}(\Spec \Fqbar)\ne \emptyset$ if and only if (\ref{eqn:frobbiglobal2}) holds. We write $\barX$, $\barcA$ etc. for $X \otimes_{\bFq}\Fqbar$, $\cA\otimes \Fqbar$ etc. Let $\sigma \colon = \id_X\otimes \Frob_q: \barX\to \barX$ and let $\pi:\barX\to X$ be the projection. Define $\dvsor(\pi): \Div(\barX)\otimes \bQ \to \Div(X)\otimes\bQ$ by $\dvsor(\pi)(\sum_i\, n_i \bar{x}_i) =  \sum_i\, n_i \pi(\bar{x}_i)$ (Note that $\deg(\dvsor(\pi)(D)) \ne \deg(D)$ in general). Part (a) of Proposition \ref{proposition:torseur1} follows from the following slightly more general result.

\begin{lemma}
\label{lemma:frobbiglobal}
Let $\barD\in \Div_0(\barcA)$. The following conditions are equivalent:

\noi (i) There exists an invertible $\barcA$-$\barcB$-bimodule $\cL$ such that $\ta(\cL(\barD)) \cong \cL$.

\noi (ii) We have 
\begin{equation}
\label{eqn:frobbiglobal2}
\sum_{x\in |X|}\, \inv_x(B) x = (\sum_{x\in |X|}\, \inv_x(A) x) + \dvsor(\pi)(\barD)\mod \Div(X).
\end{equation}
\end{lemma}

{\em Proof.} That (i) implies (ii) can be easily deduced from the corresponding local result. To show the converse we consider first the special case $\dvsor(\pi)(\barD)\in \Div(X)$, i.e.\ $\cA\simeq \cB$. Then $\barD$ can be written as a sum of divisors of the form $\pi^*(D_1), D_1\in \Div_0(X)$ and of the form $\frac{1}{e_x(\cA)}(x- \sigma(x))$ for $x\in |X|$. Hence we can assume that either $\barD = \frac{1}{e_x(\cA)}(x- \sigma(x))$ or $\barD = pr^*(D_1)$. In the first case the assertion is obvious. In the second case it follows from the fact that the homomorphism of abelian varieties
\[
\id- \Frob: \Jac_X \to \Jac_X
\]
is an isogeny hence faithfully flat.

Returning to the general case note that by \ref{proposition:moritaglobal} at least $\barcA$ and $\barcB$ are Morita equivalent. Let $\cL$ be an arbitrary invertible $\barcA$-$\barcB$-bimodule. Then $\tcL$ is also invertible hence $\ta(\cL(\barD')) \cong \cL$ for some $\barD'\in \Div_0(\barcA)$. It follows that the congruence (\ref{eqn:frobbiglobal2}) holds with $\barD'$ instead of $\barD$ as well and therefore $\dvsor(\pi)(\barD-\barD')\in \Div(X)$. Hence by what we have shown above we may alter $\cL$ by some element of $\Pic(\barcA)$ so that $\ta(\cL(\barD)) \cong \cL$.\enddemo

To prove the other assertions of \ref{proposition:torseur1} we first note that for a connected $S\in \Sch/\bFq$ and $L\in \SEl^D_{\cA,\cB,I,n}(S)$ the group of automorphisms $\Aut(L)$ of $L$ is $=\bFq^*$ if $I =0$ or $=1$ otherwise. Hence for $I>0$ the presheaf $\SE^D_{\cA,\cB,I}$ defined by
\[
\SE^D_{\cA,\cB,I,n}(S)\colon = \mbox{isomorphism classes of objects of $\SEl^D_{\cA,\cB,I,0}(S)$}
\]
is a {\it fppf} sheaf and the canonical morphism $\SEl^D_{\cA,\cB,I,n} \lra \SE^D_{\cA,\cB,I,n}$ is an isomorphism. We put $\SE^D_{\cA,\cB,I} \colon = \SE^D_{\cA,\cB,I,0}$. For $I>0$, \ref{proposition:torseur1} (c), (d) follows from:

\begin{lemma}
\label{lemma:torsor}
Suppose that $I\ne 0$ and $\SEl^D_{\cA,\cB,I,0}\ne \emptyset$. Then $\SE^D_{\cA,\cB,I}$ is a $\cW(\cA, \cB, I, D)$-torsor.
\end{lemma}

{\em Proof.} Assume first that $D=0$,  $\cA = \cB$. It follows from (\cite{lafforgue}, I.3, Th{\'e}or{\`e}me 2) that the map 
\[
(f: S \to \Spec \bFq) \mapsto f^*: \Pic_{I}(\cA)\to \SEl^0_{\cA,\cA,I}(S)
\]
yields an isomorphism between $\SEl^0_{\cA,\cA,I}$ and the trivial $\Pic_{I}(\cA)$-torsor over $\bFq$. In particular $\SE^0_{\cA,\cA,I}$ is isomorphic to the trivial $\Pic_{I,0}(\cA)$-torsor.

Now let $D \ne 0$. To simplify the notation we assume $|D| \cap |I| = \emptyset$ so that $\cC_I(\cA^T\times \cB_T)_0\cong \Pic_{I,0}(\cA)$ (the proof in the general case is analogous). Let $S\in \Sch/\bFq$ be connected and let $L_1, L_2\in \SE^D_{\cA,\cB,I}(S)$. If $L_1$ and $L_2$ have the same poles then $\xi = L_2\otimes L_1^{-1}\in \SE^0_{\cA,\cA,I}(S)\cong \Pic_{I,0}(\cA)$ by the remark above, hence $\xi L_1 = L_2$. In general there exists suitable $r_x\in \bZ$ such that $L_2$ and $(\prod_{x\in |D|}\, \theta_x^{r_x})(L_1)$ have the same poles, hence $\xi(\prod_{x\in |D|}\, \theta_x^{r_x})(L_1) = L_2$ for some $\xi\in \Pic_{I}(\cA)$. Thus $w L_1 = L_2$ for $w = \xi(\prod_{x\in |D|}\, \theta_x^{r_x})\in \cG_0$. 

Let $w\in \cW(\cA, \cB, I, D)$, $L\in \SE^D_{\cA,\cB,I}(S)$ such that $w L = L$. Write $w = \xi \prod_{x\in |D|}\theta_x^{r_x}$ with $\xi \in \Pic_I(\cA)$ and $r_x\in \bZ$. By \ref{abschnitt:theta} (ii), for $x\in |D|$ and the pole $x_S$ of $w L=L$ we have $x_S\circ \Frob_S^{r_x} = x_S$, hence $\deg(x)\mid r_x$. By \ref{abschnitt:theta} (iv) it follows that $w \in \Pic_{I,0}(\cA)$. However $w L = L$ implies that $w$ corresponds to $(w L)\otimes L^{-1} = L\otimes L^{-1}\in \SE^0_{\cA,\cA,I,0}(S)$ under the canonical bijection $\Pic_{I}(\cA)\cong\SE^0_{\cA,\cA,I,0}(S)$, i.e.\ $w=1$. This proves that for a connected $S\in \Sch/\bFq$, $\SE^D_{\cA,\cB,I}(S)$ is either empty or $\cW(\cA, \cB, I, D)$ acts simply transitively on it.

To finish the proof we have to show that $\SEl^D_{\cA,\cB,I,0}\ne \emptyset$ implies that $\SE^D_{\cA,\cB,I}(\Spec \Fqbar)\ne \emptyset$. This is a consequence of the fact that $\SEl^D_{\cA,\cB,I,0}$ is locally of finite presentation.\enddemo

Similarly, one shows that $\SE^D_{\cA,\cB,I,n}$ is a $\cW(\cA, \cB, I, D)$-torsor for all $n\in \frac{1}{d}\bZ$ with $\SEl^D_{\cA,\cB,I,n}\ne \emptyset$. In particular each $\SE^D_{\cA,\cB,I,n }$ is a finite {\'e}tale $\bFq$-scheme. This proves (b) for $I\ne 0$.

It remains to consider the case $I=0$. Choose an auxiliary level $J\in \Div(X)$ with $J> 0$ and $|D| \cap |J| = \emptyset$. A similar argument as in \ref{subsection:coarsemoduli} shows that 
\[
\SEl^D_{\cA,\cB,n}\quad\cong \quad \cA_J^*\backslash \SE^D_{\cA,\cB,J,n}
\]
Hence $\SEl^D_{\cA,\cB}$ is a Deligne-Mumford stack. Moreover as in \ref{subsection:coarsemoduli} one shows that the quotient $\SE^D_{\cA,\cB}\colon = \cA_J^*/\bFq^*)\backslash \SE^D_{\cA,\cB,J}$ is a coarse moduli scheme of $\SEl^D_{\cA,\cB,0}$ and that $\SEl^D_{\cA,\cB,0} \cong \bFq^*\backslash \SE^D_{\cA,\cB}$. Finally, since $\cW(\cA, \cB, D) \cong \cW(\cA, \cB, J, D)/(\cA_J^*/\bFq^*)$ it follows from \ref{lemma:torsor} that $\SE^D_{\cA,\cB}$ is a $\cW(\cA,\cB,D)$-torsor over $\bFq$. This completes the proof of \ref{proposition:torseur1}.

\begin{remarks}
\label{remarks:torsor2}
\rm (a) Let $D = \sum_{x\in |X|}\, m_x x \in \Div(\cA)$, $D\ne 0$ be such that $\sum_{x\in |X|}\, m_x = 0$. Condition (\ref{eqn:frobbiglobal1}) is not sufficient for $\SEl^D_{\cA,\cB,I,0}\ne \emptyset$ (compare Remark \ref{remark:stronglyequivalent}). However if, additionally, we have $\sum_{x\in |X|}\, \bZ m_x = \frac{1}{d}\bZ$ then by taking suitable products of $\theta_x$'s we see that the center of $\cC_I(\cA^T\times \cB_T)[\,\theta_x, x\in |D|\,]$ contains elements of arbitrary degree $m\in \frac{1}{d}\bZ$. Therefore $\SEl^D_{\cA,\cB,I,m}\ne \emptyset$ implies $\SEl^D_{\cA,\cB,I,0}\ne \emptyset$.

\noi (b) Suppose that $\SEl^D_{\cA,\cB,I,0}\ne \emptyset$. One can describe the $\cW(\cA, \cB, I, D)$-torsor $\SE^D_{\cA,\cB,I}/\bFq$ explicitly as follows. For $L\in \SEl^D_{\cA,\cB,I,0}(\Fqbar)$ let 
\[
\psi_L:\cW(\cA, \cB, I, D)\times \Spec \Fqbar = \coprod_{w\in \cW(\cA, \cB, I, D)} \Spec \Fqbar\to \SE^D_{\cA,\cB,I}
\]
be given on the $w$-component by the morphism corresponding to $wL$. By \ref{abschnitt:theta} (v) the diagram
\[
\begin{CD}
\cW(\cA, \cB, I, D)\times \Spec \Fqbar @>\psi_L \times \id >> \SE^D_{\cA,\cB,I}\times \Spec \Fqbar\\
@VV \Theta_D^{-1}\times \Frob_q V@VV \id\times \Frob_q V\\
\cW(\cA, \cB, I, D)\times \Spec \Fqbar @>\psi_L \times \id >> \SE^D_{\cA,\cB,I}\times \Spec \Fqbar
\end{CD}
\]
commutes. Thus $\psi_L$ induces an isomorphism
\[
\SE^D_{\cA,\cB,I}\cong (\cW(\cA, \cB, I, D)\times \Spec \Fqbar)/<\Theta_D^{-1}\times \Frob_q>
\]
\end{remarks}

\subsection{Twists of moduli spaces of $\cA$-elliptic sheaves }
\label{subsection:main}

In this section $\cA$ denotes a locally principal $\cO_X$-order of rank $d^2$ with generic fiber $A$ such that $e_{\infty}(\cA) =d$. We also assume that there is a second closed point $\fp\ne \infty$ such that $e_{\fp}(\cA) =d$ and we put $D\colon =\frac{1}{d}\infty - \frac{1}{d}\fp$. Let $\cB$ be a locally principal $\cO_X$-order of rank $d^2$ with $\Bad(\cB) = \Bad(\cA)$ and such that for the generic fiber $B$ of $\cB$ we have
\[
\sum_{x\in |X|}\, \inv_x(B) x = (\sum_{x\in |X|}\, \inv_x(A) x) + D \mod \Div(X).
\]
In order to show that the moduli spaces $\El^{\infty}_{\cA, I}$ and $\El^{\fp}_{\cB, I}$ are twists of each other we are going to define a canonical tensor product $\Ell^{\infty}_{\cA, I}\times \SEl^D_{\cA, \cB, I}\to \Ell^{\fp}_{\cB, I}$.

We introduce more notation. Recall that the notion of level structure at $\infty$ for objects of $\Ell^{\infty}_{\cA, I}$ and $\SEl^D_{\cA, \cB, I}$ and at $\fp$ for objects of $\Ell^{\fp}_{\cB, I}$ and $\SEl^D_{\cA, \cB, I}$ depend on the choice of certain local Frobenius bimodules. In order to define the tensor product (\ref{equation:tensorellsupersing}) below these choices have to be compatibly matched. For $\infty$ let $M= (\cM_{\infty}, \phi_{\infty})$ be an invertible $\phi$-$\cB_\infty$-$\cA_\infty$-bimodule of slope $-\frac{1}{d}$. For $\fp$ we choose an invertible $\phi$-$\cA_\fp$-$\cB_\fp$-bimodule $N = (\cN_{\fp}, \phi_{\fp})$ also of slope $-\frac{1}{d}$. We use $M$ to define level structure at $\infty$ and $N$ to define level structure at $\fp$. By Remark \ref{remark:phief} there exists prime elements $\varpi_{\infty}\in \cO_{\infty}$ and $\varpi_{\fp}\in \cO_{\fp}$ such that
\begin{equation*}
\label{eqn:phief2}
\phi_{\infty}^{d\deg(\infty)} = \varpi_{\infty},\quad \phi_{\fp}^{d\deg(\fp)} = \varpi_{\fp}.
\end{equation*}
Now fix a level $I\in \Div(X), I\ge 0$. We put 
\[
k(\infty)_\star = \left\{ \begin{array}{ll}
        k(\infty) & \mbox{if $\infty\not\in |I|$,}\\
        k(\infty)_d & \mbox{if $\infty\in |I|$,}
        \end{array} \right.
\qquad
k(\fp)_\sharp = \left\{ \begin{array}{ll}
        k(\fp) & \mbox{if $\fp\not\in |I|$,}\\
        k(\fp)_d & \mbox{if $\fp\in |I|$.}
        \end{array} \right.
\]
The embedding $C_I(\cA^{\infty}\times \cB_{\infty})[\theta]\to C_I(\cA^{\infty}\times \cB_{\infty})[\theta_{\infty},\theta_{\fp}]$ given by $\theta\mapsto \theta_{\infty}^{-1}$ induces an embedding $\cW(\cA,I, \infty)\to\cW(\cA, \cB, I, D)$ (recall that $C_I(\cA^{\infty}\times \cB_{\infty})\cong \Pic_{I}(\cA)$ if $\infty$ does not divide $I$). We have a short exact sequence
\[
0\lra \cW(\cA,I, \infty)\lra\cW(\cA, \cB, I, D)\lra \Gal(k(\fp)_\sharp/\bFq) \lra 0
\]
(compare Remark \ref{remark:atkinlehnerfrob}). In the following we will consider $\cW(\cA,I, \infty)$ as a subgroup of $\cW(\cA, \cB, I, D)$. By Proposition \ref{proposition:torseur1}, $\SE^D_{\cA, \cB, I}$ is a $\cW(\cA,I, \infty)$-torsor over $\Spec k(\fp)_\sharp$. 

By \ref{abschnitt:theta} (iv) there exist $\xi_{\infty}, \xi_{\fp} \in C_I(\cA^{\infty}\times \cB_{\infty})$ such that 
\[
\theta_{\infty}^{[k(\infty)_\star:\bFq]} = \xi_{\infty}, \qquad\theta_{\fp}^{[k(\fp)_\sharp:\bFq]} = \xi_{\fp}.
\]
$\xi_{\infty}$ and $\xi_{\fp}$ are given as follows. Let $\Pi_{\infty}\in \cB_{\infty}$ (resp.\  $\Pi_{\fp}\in \cA_{\fp}$) be a generator of the radical of $\cB_{\infty}$ (resp.\ of $\cA_{\fp}$). If $\infty\not\in |I|$ (resp.\ $\infty\in |I|$) then $\xi_{\infty}$ denotes the class in $C_I(\cA^{\infty}\times \cB_{\infty})$ of the idele $(\{1\}_{x\ne \infty},\Pi_{\infty}^{-1})$ (resp.\ $(\{1\}_{x\ne \infty},\varpi_{\infty}^{-1})$) in $C_I(\cA^{\infty}\times \cB_{\infty})$. If $\fp\not\in |I|$ (resp.\ $\fp\in |I|$) then $\xi_{\fp}$ denotes the class in $C_I(\cA^{\infty}\times \cB_{\infty})$ of the idele $(\{1\}_{x\ne \fp},\Pi_{\fp})$ (resp.\ $(\{1\}_{x\ne \fp},\varpi_{\fp})$) in $C_I(\cA^{\infty}\times \cB_{\infty})$.

The tensor product
\begin{equation}
\label{equation:tensorellsupersing}
 \otimes:\Ell^{\infty}_{\cA, I}\otimes_{k(\infty)_{\star}} \SEl^D_{\cA,\cB,I} \lra \Ell^{\fp}_{\cB, I}, (E,L) \mapsto E\otimes L
\end{equation}
is a morphism of $\Spec k(\fp)_{\star}$-stacks having the following properties:
\begin{itemize} 
\item[(i)]
The morphism (\ref{equation:tensorellsupersing}) is compatible with $\degcA$ and $\deg_{\cB}$, i.e.\ for $m, n\in \frac{1}{d}\bZ$ it induces a morphism\[
\Ell^{\infty}_{\cA, I,m}\otimes_{k(\infty)_{\star}} \SEl^D_{\cA,\cB,I,n} \lra \Ell^{\fp}_{\cB, I, m+n}.
\]
\item[(ii)]
The morphism of stacks 
\[
\Ell^{\infty}_{\cA, I}\otimes_{k(\infty)_\star} \SEl^D_{\cA,\cB,I}\lra \Ell^{\fp}_{\cB, I}\otimes_{k(\fp)_{\sharp}}\SEl^D_{\cA,\cB,I}, (E,L) \mapsto (E\otimes L, L)
\]
is an isomorphism with quasi-inverse
\[
 \Ell^{\fp}_{\cB, I}\otimes_{k(\fp)_{\sharp}}\SEl^D_{\cA,\cB,I}\lra\Ell^{\infty}_{\cA, I}\otimes_{k(\infty)_\star} \SEl^D_{\cA,\cB,I}, (E,L)\mapsto (E\otimes L^{-1}, L).
\]
\item[(iii)] The following diagram commutes
\[
\xymatrix@-0.5pc{
\Ell^{\infty}_{\cA, I}\otimes_{k(\infty)_{\star}} \SEl^D_{\cA,\cB,I}\ar[dd]^{\theta\otimes \theta_{\infty}}\qquad\ar[drr]^{\otimes} \\
&& \Ell^{\infty}_{\cB, I}\\
\Ell^{\infty}_{\cA, I}\otimes_{k(\infty)_{\star}} \SEl^D_{\cA,\cB,I}\qquad\ar[urr]^{\otimes}
}
\]
\item[(iv)]
For $\xi\in \cC_I(\cA^{\infty}\times \cB_{\infty})$ the following diagram commutes
\[
\xymatrix@-0.5pc{
\Ell^{\infty}_{\cA, I}\otimes_{k(\infty)_{\star}} \SEl^D_{\cA,\cB,I}\ar[dd]^{\xi\otimes \xi^{-1}}\qquad\ar[drr]^{\otimes} \\
&& \Ell^{\infty}_{\cB, I}\\
\Ell^{\infty}_{\cA, I}\otimes_{k(\infty)_{\star}} \SEl^D_{\cA,\cB,I}\qquad\ar[urr]^{\otimes}
}
\]
.
\end{itemize}

To define (\ref{equation:tensorellsupersing}) let $S\in \Sch/\bFq$ and let $\infty_S: S\to X$, $\fp_S: S\to X$ be morphisms in $\Sch/\bFq$ which factor through $\infty\to X$ and $\fp \to X$ respectively. Let $E = (\cE, \infty_S,t)$ be an $\cA$-elliptic sheaf over $S$ with zero $z:S \to X$ and let $L = (\cL,\infty_S, \fp_S, t)$ be an invertible $F$-$\cA$-$\cB$-bimodule of slope $D$. Define 
\[
E\otimes L \colon = (\cE\otimes_{\cA} \cL, \fp_S, t\otimes_{\cA} t).
\]
Note that $\cE\otimes_{\cA} \cL(-\frac{1}{d}\fp_S) = \cE(-\frac{1}{d}\infty_S)\otimes_{\cA} \cL(\frac{1}{d}(\infty_S-\fp_S))$. One easily checks that $t\otimes_{\cA}t$ is an injective $\cB\boxtimes \cO_S$-linear homomorphism with $\Coker(t\otimes_{\cA}t) \cong \Coker(t) \otimes_{\cA} \cL$. It follows from \ref{remarks:special} (c) that $E\otimes L$ is a $\cB$-elliptic sheaf with pole $\fp$ and zero $z$. Thus we have defined (\ref{equation:tensorellsupersing}) if $I=0$. 

When considering additionally level-$I$-structure, it is enough to treat separately the three cases $\infty, \fp\not\in |I|$, $|I| = \{\infty\}$ and $|I| = \{\fp\}$. In the first case if $E$ carries a level-$I$-structure $\al$ and $L$ a level-$I$-structure $\beta$ then one defines a level-$I$-structure $\al\bullet\beta$ on $E\otimes L$ as in (\ref{equation:tensorlevel}). 

Next, assume $I = n\infty$, $n>0$ and let $E= (\cE, \infty_S,t)$, $L= (\cL, t, \infty_S, \fp_S)$ be as above. Let $(\al, \la)$, $(\mu,\beta)$ be level-$I$-structures on $E$ and $L$ respectively such that $\la= \mu: S \to \Spec k(\infty)_d$ lifts $\infty_S$. Thus 
\[
\al: (\id_I\times \la)^*(\cM_I) \stackrel{\cong}{\lra} \cE|_{I\times S},  \quad \beta: \cB_I\boxtimes \cO_S\stackrel{\cong}{\lra} (\id_I\times \mu)^*(\cM_I) \otimes_{\cA}\cL|_{I\times S}.
\]
Let $\al\bullet\beta$ be the composition
\[
\al\bullet\beta: \cB_I\boxtimes \cO_S\stackrel{\beta}{\lra} (\id_I\times \mu)^*(\cM_I)\otimes_{\cA}\cL|_{I\times S}
\stackrel{\al\otimes\id}{\lra} (\cE\otimes_{\cA}\cL)|_{I\times S}.
\]
Finally, let $I = n\fp$, $n>0$ and let $\al$ and $(\mu,\beta)$ be level-$I$-structures on $E$ and $L$, i.e.\
\[
\al: \cA_I\boxtimes \cO_S \lra \cE|_{I\times S}, \quad\beta: (\id_I\times \mu)^*(\cN_I) \lra \cL|_{I\times S}
\]
where $\mu: S \to \Spec k(\fp)_d$ is a lift of $\fp_S$. We set
\[
\al\bullet\beta: (\id_I\times \mu)^*(\cN_I)\stackrel{\beta}{\lra}\cL|_{I\times S} = \cA_I\otimes_{\cA}\cL|_{I\times S}\stackrel{\al\otimes\id}{\lra}(\cE\otimes_{\cA}\cL)|_{I\times S}.
\]
In both cases one easily checks that $\al\bullet\beta$ defines a level-$I$-structure on $E\otimes L$. Thus we have defined (\ref{equation:tensorellsupersing}). The straightforward but tedious verification of the properties (i)--(iv) will be left to the reader.

Recall that $\El^{\infty}_{\cA,I}$, $\El^{\fp}_{\cB, I}$ and $\SE^D_{\cA,\cB,I}$ denote the coarse moduli spaces of $\Ell^{\infty}_{\cA, I,0}$, $\Ell^{\fp}_{\cB, I,0}$ and $\SEl^D_{\cA,\cB,I,0}$ respectively (these are fine moduli spaces if $I\ne 0$). By (i)--(iv), (\ref{equation:tensorellsupersing}) induces an $\cW(\cA, I, \infty)$-equivariant isomorphism of $\bFq$-schemes
\begin{equation}
\label{eqn:mainiso}
\El^{\infty}_{\cA, I}\otimes_{k(\infty)_\star} \SE^D_{\cA,\cB,I}\lra \El^{\fp}_{\cB, I}\otimes_{k(\fp)_{\sharp}}\SE^D_{\cA,\cB,I}.
\end{equation}
Here the action of $\xi\in\cW(\cA, I, \infty)$ on the right is given by $\id\otimes \xi$ whereas on the left it is given by $\xi^{-1}\otimes \xi$. Consequently by passing to quotients under the action and using the fact that $\SE^D_{\cA, \cB, I}$ is a $\cW(\cA,I, \infty)$-torsor over $k(\fp)_\sharp$ we obtain:

\begin{theorem}
\label{theorem:main}
The isomorphism (\ref{eqn:mainiso}) induces an isomorphism of $k(\fp)_{\sharp}$-schemes
\[
(\El^{\infty}_{\cA, I}\otimes_{k(\infty)_\star} \SE^D_{\cA,\cB,I})/\cW(\cA, I, \infty) \quad \cong \quad \El^{\fp}_{\cB, I}
\]
\end{theorem}

We shall give now another formulation of this result. Note that 
\[
\Theta_D^{[k(\fp)_{\sharp}:\bFq]} = \theta_{\infty}^{[k(\fp)_{\sharp}:\bFq]} \xi_{\fp} = \Frob_{\SE^D_{\cA,\cB,I}/k(\fp)_{\sharp}}.
\]
In particular $\Theta_D^{-[k(\fp)_{\sharp}:\bFq]}$ lies in $\cW(\cA,I, \infty)$ and is equal to $\theta^{[k(\fp)_{\sharp}:\bFq]} \xi_{\fp}^{-1}$. Hence the following diagram commutes
\begin{equation}
\label{eqn:maindia}
\xymatrix@+0.5pc{
\El^{\infty}_{\cA, I}\otimes_{k(\infty)_\star} \SE^D_{\cA,\cB,I} \ar[rr]^{(\ref{eqn:mainiso})}\ar[d]^{\theta^{[k(\fp)_{\sharp}:\bFq]} \xi_{\fp}^{-1}\otimes \Frob_{\SE/k(\fp)_{\sharp}}} & & \El^{\fp}_{\cB, I}\otimes_{k(\fp)_{\sharp}}\SE^D_{\cA,\cB,I}\ar[d]^{\id \otimes \Frob_{\SE/k(\fp)_{\sharp}}}\\
\El^{\infty}_{\cA, I}\otimes_{k(\infty)_\star} \SE^D_{\cA,\cB,I}\ar[rr]^{(\ref{eqn:mainiso})} & & \El^{\fp}_{\cB, I}\otimes_{k(\fp)_{\sharp}}\SE^D_{\cA,\cB,I}}
\end{equation}
Fix $L\in \SE^D_{\cA,\cB,I}(\Fqbar)$. Its poles correspond to $\bFq$-embeddings $\lambda: k(\infty)_\star \to \Fqbar$, $\mu: k(\fp)_{\sharp} \to \Fqbar$. By taking base change of (\ref{eqn:maindia}) with respect to the morphism $\Spec \Fqbar \to \SE^D_{\cA,\cB,I}$ corresponding to $L$ we obtain:

\begin{theorem}
\label{theorem:main2}
Let $m = [k(\fp)_{\sharp}:\bFq]$. Thus $m= \deg(\fp)$ if $\fp\not\in |I|$ and $m= d\deg(\fp)$ otherwise.
The isomorphism $\cdot\otimes L: \El^{\infty}_{\cA, I}\otimes_{k(\infty)_\star, \la} \Fqbar \to \El^{\fp}_{\cB, I}\otimes_{k(\fp)_{\sharp},\mu }\Fqbar$ induces an isomorphism of $k(\fp)_{\sharp}$-schemes
\[
(\El^{\infty}_{\cA, I}\otimes_{k(\infty)_\star, \la} \Fqbar)/<\theta^{m} \xi_{\fp}^{-1}\otimes \Frob_q^{m}> \quad \cong \quad \El^{\fp}_{\cB, I}.
\]
\end{theorem}

\begin{remark}
\label{remark:admissiblepairs}
\rm A pair $(\lambda,\mu)\in \Hom_{\bFq}(k(\infty)_\star,\Fqbar)\times \Hom_{\bFq}(k(\fp)_{\sharp},\Fqbar)$ will be called {\it admissible for $(\cA, \cB,I)$} if there exists $L\in \SE^D_{\cA,\cB,I}(\Fqbar)$ with poles $\lambda$ and $\mu$. The surjectivity of the homomorphism $\cW(\cA, \cB, I, D) \to \Gal(k(\infty)_\star/\bFq)$ implies that for all $\lambda$ there exists a $\mu$ such that $(\lambda,\mu)$ is admissible. 
\end{remark}

\subsection{Application to uniformization}
\label{subsection:uniformization}

Let $\cA$ be locally principal $\cO_X$-order of rank $d^2$ with generic fiber $A$ such that $e_{\infty}(\cA) =d$. Let $I\in\Div(X)$ denote an effective divisor. For a closed point $x\in |X|-|I|$ we denote by $\widehat{\El}^{\infty}_{\cA,I}/\Spf(\cO_{x})$ the formal completion of $\El^{\infty}_{\cA,I}$ along the fiber at $x$ of the characteristic morphism $\El^{\infty}_{\cA,I}\to X- I$. Also for an arbitrary $x\in |X|$  we let $\El^{\infty,\an}_{\cA,I}/F_x$ denote the rigid analytic space associated to $\El^{\infty}_{\cA,I}\times_X \Spec F_x$. There exists two types of uniformization of $\El^{\infty}_{\cA,I}$, i.e.\ explicit descriptions of $\widehat{\El}^{\infty}_{\cA,I}/\Spf(\cO_{\infty})$ and $\El^{\infty,\an}_{\cA,I}/F_x$ as (finite unions of) certain quotients of Drinfeld's symmetric spaces and its coverings. These are called {\it uniformization at the pole} and {\it Cherednik-Drinfeld uniformization}. The first concerns the point $x = \infty$ (under the assumption $\inv
 _{\infty} A = 0$) whereas the second the points $\fp \in |X|-\{\infty\}$ with $\inv_\fp A = \frac{1}{d}$. By using Theorem \ref{theorem:main2} we show that the two types of uniformization are equivalent (see Proposition \ref{proposition:uniformization} below).

In order to introduce the quotients of symmetric spaces appearing in the uniformization results below we have to introduce more notation. Fix a closed point $x\in X$.  We denote by $\wcO_x$ the completion of the strict henselisation of $\cO_x$ and by $\wcF_x$ its quotient field. For each positive integer $m$ we denote by $F_{x,m}$ the unramified extension of degree $m$ of $F_x$ in $\wcF_x$ and let $\cO_{x,m}$ be its ring of integers. Note that the projection $\cO_{x,m} \to k(x)_m$ has a canonical section, i.e.\ $k(x)_m\subseteq F_{x,m}$. Similarly $\overline{k(x)} \subset \wcO_x$. Denote by $D_x$ the central division algebra over $F_x$ with invariant $\frac{1}{d}$ and let $\cD_x$ be the maximal order in $D_x$. We also fix a uniformizer $\varpi_x\in \cO_x$ and an element $\Pi_x\in \cD_x$ with $\Pi_x^d = \varpi_x$. Let $\sigma$ denote the automorphism on $\cO_{x,m}$ and $\wcO_x$ which induces the $\Frob_q$ on the residue fields.

Let $\Omega^d_x$ be Drinfeld's $d-1$-dimensional symmetric space over $F_x$ and $\widehat{\Omega}^{d}_x/ \Spf(\cO_x)$ its canonical formal model (see e.g.\ \cite{genestier}). The rigid analytic variety $\Omega^d_x$ parametrizes certain formal groups. The formal scheme $\widehat{\Omega}^{d}_x$ is equipped with a canonical $\Gl_d(F_x)$-action. 

We define an action of $\Gl_d(F_x)$ on $\widehat{\Omega}^{d}_x \widehat{\otimes}_{\cO_{x}} \cO_{x,m} = \widehat{\Omega}^{d}_x \otimes_{k(x)} k(x)_m$ and $\widehat{\Omega}^{d}_x \widehat{\otimes}_{\cO_{x}} \wcO_x = \widehat{\Omega}^{d}_x \otimes_{k(x)} \overline{k(x)}$ by letting $g\in \Gl_d(F_x)$ act canonically on $\widehat{\Omega}^{d}_x$ and by $\sigma^{-v_x(\det(g))}$ on $\cO_{x,m}$ and $\cO_x$ respectively. There exists a tower of finite {\'e}tale Galois coverings (\cite{genestier}, IV.1) $\ldots \Sigma^d_{n+1,x} \to \Sigma^d_{n,x}\to \ldots \to \Sigma^d_{1,x}\to \Sigma^d_{0,x} = \Omega^d_x\otimes_{F_x} F_{x,d}$. Each $\Sigma^d_{n,x}$ carries a $\Gl_d(F_x)/\varpi_x^{\bZ}$- and $D_x^*/\varpi_x^{\bZ}$-action and the covering maps $\Sigma^d_{n+1,x} \to \Sigma^d_{n,x}$ are equivariant. Finally, for $n\ge 0$ we equip $\Sigma^d_{n,x}\otimes_{F_{x,d}}\wcF_x = \Sigma^d_{n,x}\otimes_{k(x)_d}\overline{k(x)}$ with a $\Gl_d(F_x)$- and $\cD_x$-action by letting $g\in \Gl_d(F_x)$ (or $\in D_x$) act canonically on the first factor and by $\sigma^{-v_x(\Nrd(g))}$ on the second factor.

\paragraph{Rigid analytic Drinfeld-Stuhler varieties.}

Suppose that $\inv_{\infty} A = 0$ and fix an isomorphism $A_{\infty} \cong M_d(F_{\infty})$. We write $I = n\infty + I^{\infty}$ with $\infty\not\in |I^{\infty}|$. Assume first that $n=0$. Then we define the formal $\cO_{\infty}$-scheme $\Sh_{\cA, I}^{\infty}$ by
\[
\wSh_{\cA, I}^{\infty}\colon = A^*\backslash \left(A^*(\bA^{\infty})/U_I(\cA^{\infty})\times \widehat{\Omega}^{d}_{\infty}\right).
\]
Next assume $\infty\in |I|$ and write $I = n\infty + I^{\infty}$ with $\infty\not\in |I^{\infty}|$. Then put
\[
\Sh_{\cA, I}^{\infty} \colon =A^*\backslash \left(A^*(\bA^{\infty})/U_{I^{\infty}}(\cA^{\infty})\times \Sigma^d_{n,\infty}\right).
\]
This is rigid analytic space over $F_{\infty}$.

There exists a canonical right action of the group $\cC_I(\cA^{\infty}\times \cD_{\infty})$ on $\wSh_{\cA, I}^{\infty}$ and $\Sh_{\cA, I}^{\infty}$ which is defined as follows. Let $a = (\{a_x\}_{x\ne \infty},d_{\infty})\in (\prod_{x\in |X|-\{\infty\}}' \, N(\cA_x))\times  N(\cD_{\infty})$ and assume first $n=0$. Then the right action of the class $[a]\in \cC_I(\cA^{\infty}\times \cD_{\infty})$ of $a$ on $\wSh_{\cA, I}^{\infty}$ is given by right multiplication by $\{a_x\}_{x\ne \infty}$ on $A^*(\bA^{\infty})/U_I(\cA^{\infty})$. Now assume that $n>0$. Then $[a]$ acts on $\Sh_{\cA, I}^{\infty}$ by right multiplication of $\{a_x\}_{x\ne \infty}$ on $A^*(\bA^{\infty})/U_{I^{\infty}}(\cA^{\infty})$ and letting $d_{\infty}^{-1}$ act on $\Sigma^d_{n,\infty}$. 

There are canonical morphism 
\begin{align}
\label{eqn:poleanalytic1}
& \pole:\wSh_{\cA, I}^{\infty} \to \Spec k(\infty)& \mbox{if $n=0$,}\\
\label{eqn:poleanalytic2} & \pole:\Sh_{\cA, I}^{\infty} \to \Spec k(\infty)_d & \mbox{if $n>0$}.
\end{align}
For that let $l=l_{\infty}: A^*(\bA^{\infty}) \to \bZ$ be the composite
\[
l_{\infty}: A^*(\bA^{\infty})\stackrel{\Nrd}{\lra} F^*(\bA^{\infty})\stackrel{\dvsor}{\lra} \bigoplus_{x\ne \infty} \bZ x \stackrel{\deg}{\lra}\bZ.
\]
Note that for $a\in A^*\subset A^*(\bA^{\infty})$ we have $l_{\infty}(a) = -\deg(\infty) v_{\infty}(\Nrd(a))$. Now assume $n=0$ and let
\begin{equation}
\label{eqn:poleanalytic3}
A^*(\bA^{\infty})/U_I(\cA^{\infty})\times \widehat{\Omega}^{d}_{\infty} \lra \Spec k(\infty)
\end{equation}
be given on the component $\eta U_I(\cA^{\infty})\times \widehat{\Omega}^{d}_{\infty}$ by 
\begin{equation}
\label{eqn:poleanalytic4}
\begin{CD}
\widehat{\Omega}^{d}_{\infty}@>>> \Spec k(\infty) @>\Frob_q^{-l(\eta)}>> \Spec k(\infty). 
\end{CD}
\end{equation}
Clearly, (\ref{eqn:poleanalytic4}) factors through $\wSh_{\cA, I}^{\infty}$, hence it induces (\ref{eqn:poleanalytic1}).

Now suppose $n>0$. Since $k(\infty)_d\subseteq F_{\infty,d}$, we get a map $\Sigma^d_{n,\infty}\to \Spec F_{\infty,d} \to \Spec k(\infty)_d$. Note that for $g\in \Gl_d(F_{\infty})$ the diagram
\begin{equation}
\label{eqn:poleanalytic5}
\begin{CD}
\Sigma^d_{n,\infty}@> g\cdot >>\Sigma^d_{n,\infty}\\
@VVV@VVV\\
\Spec k(\infty)_d@> \Frob_q^{-v_{\infty}(\det(g))}>>\Spec k(\infty)_d
\end{CD}
\end{equation}
commutes. We define
\begin{equation}
\label{eqn:poleanalytic6}
A^*(\bA^{\infty})/U_I(\cA^{\infty})\times \Sigma^d_{n,\infty}\to \Spec k(\infty)_d
\end{equation}
on the component corresponding to $\eta U_I(\cA^{\infty})\in A^*(\bA^{\infty})/U_I(\cA^{\infty})$ by
\[
\begin{CD}
\Sigma^d_{n,\infty} @>>> \Spec k(\infty)_d @>\Frob_q^{-l(\eta)}>> \Spec k(\infty)_d.
\end{CD}
\] 
The commutativity of (\ref{eqn:poleanalytic5}) implies that (\ref{eqn:poleanalytic6}) factors through $\Sh_{\cA, I}^{\infty}$, i.e.\ it yields the  map (\ref{eqn:poleanalytic2}).

\paragraph{Cherednik-Drinfeld varieties.} Let $\wxi  = \{\wxi_x\}_{x\ne \infty}\in (\prod_{x\in |X|-\{\infty\}}' \, N(\cA_x))$ and let $\xi\in \cC_I(\cA^{\infty}\times \cD_{\infty})$ be the idele class represented by $(\{\wxi_x\}_{x\ne \infty}, 1)\in (\prod_{x\in |X|-\{\infty\}}' \, N(\cA_x))\times  N(\cD_{\infty})$. We assume that $\xi$ is a central element in $\cC_I(\cA^{\infty}\times \cD_{\infty})$ and that $m = -d\deg_{\cA}(\xi) = -l_{\infty}(\wxi)\ne 0$. 
We define 
\begin{align*}
& \wCD_{\cA, I,\infty}^{\,\xi}\colon = A^*\backslash \left(A^*(\bA^{\infty})/U_I(\cA^{\infty})\wxi^{\bZ}\times \widehat{\Omega}^{d}_{\infty}\otimes_{k(\infty)} \overline{k(\infty)}\right) & \mbox{if $n=0$,}\\
& \ChDr_{\cA, I, \infty}^{\xi} \colon = A^*\backslash \left(A^*(\bA^{\infty})/U_{I^{\infty}}(\cA^{\infty})\wxi^{\bZ}\times \Sigma^d_{n,\infty}\otimes_{k(\infty)_d} \overline{k(\infty)}\right) & \mbox{if $n>0$}.
\end{align*}
As above one defines a right action of $\cC_I(\cA^{\infty}\times \cD_{\infty})$ on $\wCD_{\cA, I,\infty}^{\,\xi}$ and $\ChDr_{\cA, I,\infty}^{\,\xi}$ by letting $a = (\{a_x\}_{x\ne \infty},d_{\infty})\in (\prod_{x\in |X|-\{\infty\}}' \, N(\cA_x))\times  N(\cD_{\infty})$ act by right \break multiplication by $\{a_x\}_{x\ne \infty}$ on $A^*(\bA^{\infty})/U_I(\cA^{\infty})\wxi^{\bZ}$ and letting $d_{\infty}^{-1}$ act on \break $\widehat{\Omega}^{d}_{\infty}\widehat{\otimes}_{\cO_{\infty}}\wcO_{\infty}$ (if $n=0$) and $\Sigma^d_{n,\infty}\otimes_{F_{\infty,d}}\wcF_{\infty}$ (if $n>0$). Note that $\xi$ acts trivially.

Let $k(\xi)$ denote the fixed field of $\Frob_q^m$ in $\overline{k(\infty)}$. There are canonical morphisms 
\begin{align}
\label{eqn:polecd1}
& \wCD_{\cA, I,\infty}^{\,\xi} \to \Spec k(\xi)& \mbox{if $n=0$,}\\
\label{eqn:polecd2} & \ChDr_{\cA, I,\infty}^{\xi} \to \Spec k(\xi) & \mbox{if $n>0$}.
\end{align}
Their definition is similar to the definition of (\ref{eqn:poleanalytic1}) and (\ref{eqn:poleanalytic2}). For example (\ref{eqn:polecd1}) is induced by the maps
\[
\begin{CD}
\eta (U_I(\cA^{\infty})\wxi^{\bZ})\times \widehat{\Omega}^{d}_{\infty}\otimes_{k(\infty)} \overline{k(\infty)} @>>> \Spec k(\xi) @>\Frob_q^{-l(\eta)}>> \Spec k(\xi). 
\end{CD}
\]
The rigid analytic varieties $\Sh_{\cA, I}^{\infty}$ and $\ChDr_{\cA, I,\infty}^{\xi}$ are twists of each other. More precisely we have the following result.

\begin{lemma}
\label{lemma:twisting}
(a) There exists a canonical isomorphism of formal schemes over $\Spf(\cO_{\infty})$ (resp.\ rigid analytic varieties over $F_{\infty}$)
\begin{align}
\label{eqn:twisting1}
& \wSh_{\cA, I}^{\infty}\otimes_{k(\infty)} \overline{k(\infty)}/< \xi\otimes \Frob_q^{m}> \quad \cong \quad \wCD_{\cA, I,\infty}^{\,\xi} & \mbox{for $n=0$, resp.}\\
 \label{eqn:twisting2}& \Sh_{\cA, I}^{\infty}\otimes_{k(\infty)_d} \overline{k(\infty)}/< \xi\otimes \Frob_q^{m}> \quad \cong \quad \ChDr_{\cA, I,\infty}^{\xi} & \mbox{for $n>0$}.
\end{align}
Here $\wSh_{\cA, I}^{\infty}\otimes_{k(\infty)} \overline{k(\infty)}$ (resp.\ $\Sh_{\cA, I}^{\infty}\otimes_{k(\infty)_d} \overline{k(\infty)_d}$) denotes the base change to $\overline{k(\infty)}$ of the morphism (\ref{eqn:poleanalytic1}) (resp.\ (\ref{eqn:poleanalytic2})).

\noi (b) Let $\xi_{\infty}\in \cC_I(\cA^{\infty}\times \cD_{\infty})$ be the class of the idele $(\{1\}_{x\ne \infty},\Pi_{\infty}^{-1})$ (if $n=0$) resp.\ of $(\{1\}_{x\ne \infty},\varpi_{\infty}^{-1})$ (if $n>0$). Then we have
\begin{align}
\label{eqn:twisting3}
& \wCD_{\cA, I,\infty}^{\,\xi}\otimes_{k(\xi)} \overline{k(\xi)}/< \xi_{\infty}\otimes \Frob_q^{\deg(\infty)}> \,\cong \, \wSh_{\cA, I}^{\infty} & \mbox{for $n=0$, resp.}\\
 \label{eqn:twisting4} & \ChDr_{\cA, I,\infty}^{\xi}\otimes_{k(\xi)} \overline{k(\xi)}/< \xi_{\infty}\otimes \Frob_q^{d\deg(\infty)}> \, \cong \, \Sh_{\cA, I}^{\infty}
& \mbox{for $n>0$}.
\end{align}
\end{lemma}

{\em Proof.} We prove only the existence of (\ref{eqn:twisting1}). The other cases are similar and will be left to the reader. For $\eta U_I(\cA^{\infty})\in A^*(\bA^{\infty})/U_I(\cA^{\infty})$ we denote the base change of the map (\ref{eqn:poleanalytic4}) to $\overline{k(\infty)}$  by $(\eta U_I(\cA^{\infty})\times \widehat{\Omega}^{d}_{\infty})\otimes_{k(\infty)} \overline{k(\infty)}$. Let
\begin{equation}
\label{eqn:twisting5}
\begin{CD}
(\eta U_I(\cA^{\infty})\times \widehat{\Omega}^{d}_{\infty})\otimes_{k(\infty)} \overline{k(\infty)}@> \id \otimes \Frob_q^{l(\eta)} >>  \widehat{\Omega}^{d}_{\infty}\otimes_{k(\infty)} \overline{k(\infty)} 
\end{CD}
\end{equation}
and let
\begin{equation}
\label{eqn:twisting6}
\left(A^*(\bA^{\infty})/U_I(\cA^{\infty})\times \widehat{\Omega}^{d}_{\infty}\right)\otimes_{k(\infty)}\overline{k(\infty)}\to A^*(\bA^{\infty})/U_I(\cA^{\infty})\times (\widehat{\Omega}^{d}_{\infty}\otimes_{k(\infty)} \overline{k(\infty)})
\end{equation}
be made up of all the morphisms (\ref{eqn:twisting3}). One easily checks that it is $A^*$-equivariant and that the following diagram commutes:
\begin{equation*}
\begin{CD}
(A^*(\bA^{\infty})/U_I(\cA^{\infty})\times \widehat{\Omega}^{d}_{\infty})\otimes \overline{k(\infty)}@> (\ref{eqn:twisting6}) >> A^*(\bA^{\infty})/U_I(\cA^{\infty})\times (\widehat{\Omega}^{d}_{\infty}\otimes \overline{k(\infty)})\\
@VV \cdot \xi\otimes \Frob_q^{m} V @VV \cdot \wxi\otimes \id V\\
(A^*(\bA^{\infty})/U_I(\cA^{\infty})\times \widehat{\Omega}^{d}_{\infty})\otimes \overline{k(\infty)}@> (\ref{eqn:twisting6}) >> A^*(\bA^{\infty})/U_I(\cA^{\infty})\times (\widehat{\Omega}^{d}_{\infty}\otimes \overline{k(\infty)})
\end{CD}
\end{equation*}
Hence (\ref{eqn:twisting6}) induces the isomorphism (\ref{eqn:twisting1}).\enddemo

Note that, since $\xi_{\infty}$ acts trivially on $\wSh_{\cA, I}^{\infty}$ (resp.\ $\Sh_{\cA, I}^{\infty}$), the $\cC_I(\cA^{\infty}\times \cD_{\infty})$-action on $\wSh_{\cA, I}^{\infty}$ (resp.\ $\Sh_{\cA, I}^{\infty}$) induces a right $\cW(\cA,I,\infty)$-action (by Remark \ref{remarks:modaut} (c)). In terms of the latter, Lemma \ref{lemma:twisting} (a) can be reformulated as follows:
\begin{align}
\label{eqn:twisting7}
& \wSh_{\cA, I}^{\infty}\otimes_{k(\infty)} \overline{k(\infty)}/< \xi\theta^m\otimes \Frob_q^{m}> \quad \cong \quad \wSh_{\cA, I,\infty}^{\,\xi}& \mbox{for $n=0$, resp.}\\
 \label{eqn:twisting8}& \Sh_{\cA, I}^{\infty}\otimes_{k(\infty)_d} \overline{k(\infty)}/< \xi\theta^m\otimes \Frob_q^{m}> \quad \cong \quad \Sh_{\cA, I,\infty}^{\xi}& \mbox{for $n>0$}.
\end{align}

\paragraph{Uniformization at the pole.} Suppose that $\inv_{\infty} A = 0$ and assume first that $\infty$ does not divide the level $I$.  Then there exists an isomorphism of formal schemes over $\Spf(\cO_{\infty})$ 
\begin{equation}
\label{eqn:uniinf1}
\widehat{\mbox{Ell}}_{\cA,I}^{\infty}/\Spf(\cO_{\infty}) \quad\cong \quad\wSh_{\cA, I}^{\infty}
\end{equation}
which is compatible with the $\cW(\cA,I,\infty)$-action and the morphisms $\pole$. 
 
Now assume $\infty\in |I|$. Then we expect
\begin{equation}
\label{eqn:uniinf2}
\El^{\infty,\an}_{\cA,I}/F_{\infty}\quad\cong \quad\Sh_{\cA, I}^{\infty}.
\end{equation}
Again, (\ref{eqn:uniinf2}) should be compatible with the $\cW(\cA,I,\infty)$-action and the morphisms $\pole$. 

We say that $\El^{\infty}_{\cA,I}$ admits {\it uniformization at the pole} if (\ref{eqn:uniinf1}) (resp.\ (\ref{eqn:uniinf2})) holds. Suppose that $\infty\not\in |I|$. (\ref{eqn:uniinf1}) has been proved in (\cite{stuhler}, 4.4) if $A$ is a division algebra or $A= M_d(F)$. As in loc.\ cit.\ the general case can be easily deduced from (\cite{stuhler2}, Corollary, p.\ 531 and Theorem 1, p.\ 538) or (\cite{genestier}, III.3.1.1). If $\infty\in |I|$ then the uniformization (\ref{eqn:uniinf2}) is proved in \cite{drinfeld3} in the case of Drinfeld modular varieties (i.e.\ $A = M_d(F)$).

\paragraph{Cherednik-Drinfeld uniformization.}
\label{paragraph:cherednikdrinfeld}
Let $\fp\in |X|-\{\infty\}$ and assume that $\inv_{\fp} A = \frac{1}{d}$. Let $\cB$ be a locally principal $\cO_X$-order of rank $d^2$ with $\Bad(\cB) = \Bad(\cA)$ and such that the local invariants of the generic fiber $B$ of $\cB$ are given by $\inv_{\infty}(B) = \inv_{\infty}(A) + \frac{1}{d}$, $\inv_{\fp}(B) = 0$ and $\inv_x(B) = \inv_x(A)$ for all $x\in |X|-\{\infty, \fp\}$. We fix an isomorphism $B_{\fp} \cong M_d(F_{\fp})$ and isomorphisms $\cB_x \cong \cA_x$ for all $x\in |X|-\{\infty, \fp\}$. Using the latter we can identify the groups $\cC_I(\cA^{\infty}\times \cB_{\infty})$ and $C_I(\cB^{\fp}\times \cA_{\fp})$. Since $\xi_{\infty}$ acts trivially on $\wCD_{\cB, I,\fp}^{\,\xi_{\infty}}$ (resp.\ $\ChDr_{\cB, I,\fp}^{\xi_{\infty}}$) we obtain a right $\cW(\cA,I,\infty)\cong C_I(\cB^{\fp}\times \cA_{\fp})/\xi_{\infty}^{\bZ}$-action on $\wCD_{\cB, I,\fp}^{\,\xi_{\infty}}$ (resp.\ $\ChDr_{\cB, I,\fp}^{\xi_{\infty}}$). We also fix an 
isomorphism $k(\infty)_{\star} \cong k(\xi_{\infty})\subset \overline{k(\fp)}$ such that the pair $(k(\infty)_{\star}\cong k(\xi_{\infty})\hookrightarrow \overline{k(\fp)}, k(\fp)_{\sharp} \hookrightarrow\overline{k(\fp)})$ is admissible in the sense of Remark \ref{remark:admissiblepairs} ($k(\infty)_{\star}$ and $k(\fp)_{\sharp}$ are defined as in the last section) and define 
\begin{align*}
& \pole: \wCD_{\cB, I,\fp}^{\,\xi_{\infty}} \stackrel{(\ref{eqn:polecd1})}{\lra}\Spec k(\xi_{\infty}) \cong \Spec k(\infty)_{\star} & \mbox{if $n=0$,}\\
 & \pole:\ChDr_{\cB, I,\fp}^{\xi_{\infty}} \stackrel{(\ref{eqn:polecd2})}{\lra} \Spec k(\xi_{\infty}) \cong \Spec k(\infty)_{\star} & \mbox{if $n>0$}.
\end{align*}

Assume that $\fp\not\in |I|$ (resp.\ $\fp\in |I|$). Then we expect that there is a canonical isomorphism of formal $\cO_{\fp}$-schemes (resp.\ of rigid analytic spaces over $F_{\fp}$)
\begin{align}
\label{eqn:uniinffp1} & \widehat{\mbox{Ell}}_{\cA,I}^{\infty}/\Spf(\cO_{\fp})\quad\cong \quad\wCD_{\cB, I,\fp}^{\,\xi_{\infty}} & \mbox{if $\fp\not\in |I|$,}\\
\label{eqn:uniinffp2} & \El^{\infty,\an}_{\cA,I}/F_{\fp} \quad\cong \quad\ChDr_{\cB, I,\fp}^{\xi_{\infty}} & \mbox{if $\fp\in |I|$}\end{align}
compatible with $\cW(\cA,I,\infty)$-actions and the morphisms $\pole$.

We say that $\El^{\infty}_{\cA,I}$ admits {\it Cherednik-Drinfeld uniformization at $\fp$} if (\ref{eqn:uniinffp1}) (resp.\ (\ref{eqn:uniinffp2})) holds. Both (\ref{eqn:uniinffp1}) and (\ref{eqn:uniinffp2}) are proved in (\cite{hausberger}, 8.1 and 8.3) in the case $\deg(\infty) =1$, $\inv_{\infty} A = 0$ and $\infty\not\in |I|$. Under these assumptions $\wCD_{\cB, I,\fp}^{\,\xi_{\infty}}$ and $\ChDr_{\cB, I,\fp}^{\xi_{\infty}}$ have the following simpler description
\begin{align*}
& \wCD_{\cB, I,\fp}^{\,\xi_{\infty}}= B^*\backslash \left(B^*(\bA^{\fp,\infty})/U_I(\cB^{\infty,\fp})\times \widehat{\Omega}^{d}_{\fp}\otimes_{k(\fp)} \overline{k(\fp)}\right) & \mbox{if $n=0$,}\\
& \ChDr_{\cB, I,\fp}^{\xi_{\infty}}= B^*\backslash \left(B^*(\bA^{\fp,\infty})/U_I(\cB^{\infty,\fp})\times \Sigma^d_{n,\fp}\otimes_{k(\fp)_d} \overline{k(\fp)}\right) & \mbox{if $n>0$.}
\end{align*}
Here $n$ denotes now the exact multiple of $\fp$ occurring in $I$. 

By combining Theorem \ref{theorem:main2}, (\ref{eqn:twisting7}), (\ref{eqn:twisting8}) and Lemma \ref{lemma:twisting} (b) we obtain:

\begin{proposition}
\label{proposition:uniformization}
Let $\fp \in |X| -\{\infty\}$ and let $\cA$ and $\cB$ be locally principal $\cO_X$-orders of rank $d^2$ with $\Bad(\cB) = \Bad(\cA)$ such that the local invariants of the generic fibers $A$ and $B$ are given by $\inv_{\infty}(A) =0, \inv_{\infty}(B) = \frac{1}{d}$, $\inv_{\fp}(A)= \frac{1}{d}, \inv_{\fp}(B) = 0$ and $\inv_x(B) = \inv_x(A)$ for all $x\in |X|-\{\infty, \fp\}$. The following conditions are equivalent:

\noi (i) $\El^{\infty}_{\cA,I}$ admits uniformization at the pole.

\noi (ii) $\El^{\fp}_{\cB,I}$ admits Cherednik-Drinfeld uniformization at $\infty$.
\end{proposition}

By applying \ref{proposition:uniformization} to the results of \cite{stuhler} and \cite{hausberger} we obtain further cases 
where $\El^{\infty}_{\cA,I}$ admits uniformization at the pole or Cherednik-Drinfeld uniformization. For example if $\inv_{\infty}(A) =0$, $\infty\in |I|$ and if there exists a point $\fp \in |X| -\{\infty\}$ such that $\inv_{\fp}(A)= \frac{1}{d}$ and $\deg(\fp)=\frac{1}{d}$ then $\El^{\infty}_{\cA,I}$ admits uniformization at the pole. Conversely Cherednik-Drinfeld uniformization for $\El^{\infty}_{\cA,I}$ holds whenever if $\fp$ does not divide the level. 

\section{Appendix}
\label{section:appendix}

\subsection{Commutative subalgebras in semisimple algebras}

Let $k$ be a perfect field and $A$ a finite-dimensional semisimple $k$-algebra. We collect a few facts about maximal separable and commutative subalgebras of $A$ for which we could not find any references.

Let $Z$ denote the center of $A$. By Wedderburn's Theorem we have $Z \cong k_1 \times \ldots \times k_r$ for some finite separable extensions $k_i/k$. For a finite $Z$-module $M$, $\rank_{Z} M$ denotes the (not necessarily constant) rank of the corresponding locally free $\cO_{\Spec Z}$-module.

\begin{lemma} 
\label{lemma:maxcommsepsubalgebra1}
Let $T$ be a commutative separable $k$-subalgebra of $A$. The following conditions are equivalent.

\noi (i) $T = Z_A(T) = \{x\in A\mid \, t x = x t \quad\forall t\in T\}$.

\noi (ii)  $T$ is a maximal commutative separable $k$-subalgebra of $A$.

\noi (iii) $T\supseteq Z$ and $(\rank_{Z} T)^2 = \rank_{Z} A$.

Moreover if $A =  \End_{k_1}(V_1)\times \ldots \times \End_{k_r}(V_r)$ where $V_i$ a finite-dimensional $k_i$-vector space for $i=1, \ldots, r$, then (i) -- (iii) are equivalent to 

\noi (iv) $V_1\oplus\ldots\oplus V_r$ is a free $T$-module of rank $1$.
\end{lemma}

A commutative separable $k$-subalgebra $T$ of $A$ satisfying the equivalent conditions (i) -- (iii) above will be called a {\it maximal torus} of $A$. 

\begin{lemma} 
\label{lemma:maxcommsepsubalgebra2} 
Let $T_1, T_2$ be two maximal tori of $A$. Then there exists a finite extension $k'/k$ such that $T_1\otimes_k k'$ and $T_2\otimes_k k'$ are conjugated in $A\otimes_k k'$.
\end{lemma}

A finite $A$-module $M$ is called a generator of $\Mod_A$ if the functor 
\[
\Hom_A(M, \cdot): \Mod_A\lra \Mod_k
\]
is faithful. $M$ is called a minimal generator if $\dim_k(M)$ is minimal. Assume now that $A$ is split, i.e.\ $A = \End_{k_1}(V_1)\times \ldots \times \End_{k_r}(V_r)$ as in condition (iv) of Lemma \ref{lemma:maxcommsepsubalgebra1} and let $T$ be a maximal torus in $A$. We have

\begin{lemma} 
\label{lemma:maxcommsepsubalgebra3} 
Let $M$ be a finite $A$-module. The following conditions are equivalent.

\noi (i) $M$ is a minimal generator.

\noi (ii) $M\cong V_1\oplus \ldots\oplus V_r$

\noi (iii) $M$ is a free $T$-module of rank 1.
\end{lemma}

\subsection{$\cA$-elliptic sheaves according to Laumon-Rapoport-Stuhler}
\label{subsection:aelllrs}

The aim of this section is to show that under suitable assumptions on $\cA$ the moduli stack $\Ell^{\infty}_{\cA}$ defined in section \ref{subsection:aell} is isomorphic to the stack defined in (\cite{lrs}, 2.4). 

Firstly, we establish an equivalence between certain parabolic vector bundles and locally free modules of a hereditary algebra. We use the following notations and assumptions. Let $k$ be a perfect field of cohomological dimension $\le 1$ and let $X$ be a smooth connected curve over $k$ and $F$ is the function field of $X$. We also fix a closed point $\infty\in X$. To simplify the notation we assume that $\deg(\infty) = 1$ (see Remark \ref{remarks:degbeliebig} below for the case $\deg(\infty) >1$).

Let $\cA'$ be a locally principal $\cO_X$-order of rank $d^2$ with generic fiber $A'$. We assume that $e_{\infty}(\cA') =1$, i.e.\ $\cA'_{\infty} \cong M_d(\cO_{\infty})$. To begin with we introduce the notion of a parabolic $\cA'$-modules and parabolic vector bundles with $\cA'$-action (compare \cite{yokogawa}). A {\it filtered object} in a category $\cC$ is a functor $\Cd: \bZ \to \cC$. Morphisms of filtered objects are natural transformations. Here we regard the ordered set $\bZ$ as a category in the usual way. The set of objects is $\bZ$ and for $i,j\in \bZ$ we have 
\[
\sharp(\Mor(i,j)) = \left\{ \begin{array}{ll} 1 & \mbox{if $i\le j$}\\
                  0 & \mbox{otherwise.}
                \end{array}\right.
\]
For $i\in\bZ$ the morphism $C_i\to C_{i+1}$ will be denoted by $j_i = j^{C}_i$. For a filtered object $\Cd$ in $\cC$ and $n\in \bZ$ the shifted filtered object $C[n]_{\star}$ is defined as the composite $\bZ\stackrel{+n}{\lra}\bZ\to\cC$. A morphism $\phi: \Cd\to \Dd$ of filtered objects induces a morphism $\phi[n]:C[n]_{\star}\to D[n]_{\star}$. 

Recall that for $S\in \Sch/k$ we have set $_{\cA'}\Mod(S) \colon = {_{\cA'\boxtimes \cO_S}\Mod}$ (resp.\ $\Mod_{\cA'}(S)\colon = \Mod_{\cA'\boxtimes \cO_S}$). 

\begin{definition}
\label{definition:parvecact}
Let $S$ be a $k$-scheme.

\noi (a) For $e\in \bZ$ with $e\ge 1$ let $\PMod_{\cA',e}(S)$ denote the category of pairs $(\cFd, \psi_{\star})$ consisting of a filtered $\Mod_{\cA'}(S)$-object $\cFd$ and an isomorphism $\psi_{\star}: \cF[e]_{\star}\to \cFd(\infty)\colon = \cFd\otimes_{\cO_{X\times S}} (\cO_X(\infty)\boxtimes \cO_S)$ such that the restriction of $j_i:  \cF\to \cF_{i+1}$ to $X-\{\infty\}\times S$ is an isomorphism and such that the following diagram commutes
\begin{equation}
\label{eqn:par}
\xymatrix@-0.5pc{
&& \cF_{i+e}\ar[dd]_{\psi_i}\\
\cF_i\ar[urr]^{j_{i+e-1}\circ\ldots \circ j_i} \ar[drr]_{\id \otimes \io} \\
&& \cF_i(\infty)
}
\end{equation}
where $\io: \cO_{X\times S} \hookrightarrow \cO_X(\infty)\boxtimes \cO_S$ is the inclusion. Morphisms in $\PMod_{\cA',e}(S)$ are morphisms of filtered objects compatible with the isomorphisms $\psi$.

\noi (b) Let $\PCoh^r_{\cA',sp,e}(S)$ denote the groupoid of $(\cKd, \psi_{\star})$ in $\PMod_{\cA',e}(S)$ such that $\cK_i\in \Coh^r_{\cA',sp}(S)$ and $N(\cKd) \colon =N(\cK_i) = N(\cK_{i+1})$ for all $i\in \bZ$.

\noi (c) For $e,r\in \bZ$ with $e,r\ge 1$ and $e\mid rd$. We denote by $\PVect^r_{\cA',e}(S)$ the full subcategory of $(\cFd, \psi_{\cFd})$ in $\PMod_{\cA',e}(S)$ such that $\cF_i\in \Vect^r_{\cA'}(S)$ for all $i\in \bZ$ and such that $\Coker(j_{\star}: \cFd\to \cF[1]_{\star})\in \PCoh^s_{\cA',sp,e}(S)$ with $s= \frac{rd}{e}$.

Similarly one defines $_{\cA'}\PMod_{e}(S)$ and $_{\cA'}\PVect^r_e(S)$ using left $\cA'\boxtimes \cO_S$-modules.
\end{definition}

Note that for $(\cFd, \psi_{\star})$ in $\PMod_{\cA',e}(S)$ with $\cF_i\in \Vect^r_{\cA'}(S)$ for all $i\in \bZ$, the commutativity of diagram (\ref{eqn:par}) implies that $j_i: \cF_i\to \cF_{i+1}$ is injective and $\Coker(j_i)$ is a sheaf on $\infty \times S$. For $\cA'= \cO_X$ we write $\Mod_X$, $\PMod_{X,e}$, $\Vect_X$ etc. for $\Mod_{\cO_X}$, $\PMod_{\cO_X,e}(S)$, $\Vect_{\cO_X}$ etc. 

Let $\cEd\in \PMod_{\cA',e}(S)$ and $\cFd\in {_{\cA'}\PMod_{e}}(S)$. We are going to define now a tensor product $(\cEd\otimes_{\cA'}\cFd)_{\star}$. For $i\in \bZ$ we set
\[
T_i(\cEd, \cFd) \colon = \bigoplus_{\la+ \mu =i, \la,\mu\in\bZ} \, \cE_{\la}\otimes_{\cA'} \cF_{\mu}.
\]
For $i\in \bZ$ with we define homomorphisms 
\begin{eqnarray*}
\alpha_i:T_i(\cEd, \cFd)\lra T_{i+1}(\cEd, \cFd),\quad\beta_{i}:T_{i}(\cEd, \cFd)\lra T_{i+1}(\cEd, \cFd)
\end{eqnarray*} 
as the direct sums of the inclusions $j_{\la}\otimes \id: \cE_{\la}\otimes_{\cA'} \cF_{\mu}\to\cE_{\la+1}\otimes_{\cA'} \cF_{\mu}$ (resp.\ $\id\otimes j_{\mu}:\cE_{\la}\otimes_{\cA'} \cF_{\mu}\to\cE_{\la}\otimes_{\cA'} \cF_{\mu+1}$). Also let 
\[
\gamma_i: T_i(\cEd, \cFd)\lra T_i(\cEd, \cFd)
\]
be the isomorphism given on the summand $\cE_{\la}\otimes_{\cA'} \cF_{\mu}$ by 
\[
\cE_{\la}\otimes_{\cA'} \cF_{\mu} \cong \cE_{\la}(\infty)\otimes_{\cA'} \cF_{\mu}({-\infty}) \stackrel{\psi^{-1}\otimes \psi^{-1}}{\lra} \cE_{\la+e}\otimes_{\cA'} \cF_{\mu-e}.
\]
Finally let 
\[
\delta_{i}: T_{i-1}(\cEd, \cFd)\oplus T_{i}(\cEd, \cFd)\lra T_{i}(\cEd, \cFd)
\]
be given on the summand $T_{i-1}(\cEd, \cFd)$ by $\alpha_{i-1}-\beta_{i-1}$ and by $id-\gamma_{i}$ on $T_{i}(\cEd, \cFd)$. We define
\begin{equation}
\label{eqn:tenpara}
\sum_{\la+ \mu =i, \la,\mu\in\bZ} \, \cE_{\la}\otimes_{\cA'} \cF_{\mu} = \Coker(\delta_i)
\end{equation}
There are canonical morphisms
\begin{equation}
\label{eqn:tenparaincl}
\sum_{\la+ \mu =i-1, \la,\mu\in\bZ} \, \cE_{\la}\otimes_{\cA'} \cF_{\mu} \lra \sum_{\la+ \mu =i, \la,\mu\in\bZ} \, \cE_{\la}\otimes_{\cA'} \cF_{\mu}.
\end{equation}
The isomorphisms 
\[
\cE_{\la + d}\otimes_{\cA'} \cF_{\mu}\stackrel{\psi\otimes \id}{\lra}(\cE_{\la}\otimes_{\cA'} \cF_{\mu})(\infty),\qquad \cE_{\la }\otimes_{\cA'} \cF_{\mu+d}\stackrel{\id\otimes \psi}{\lra}(\cE_{\la}\otimes_{\cA'} \cF_{\mu})(\infty)
\]
induces an isomorphism
\begin{equation}
\label{eqn:tenparaphi}
\sum_{\la+ \mu =i +d, \la,\mu\in\bZ} \, \cE_{\la}\otimes_{\cA'} \cF_{\mu}\lra (\sum_{\la+ \mu =i, \la,\mu\in\bZ} \, \cE_{\la}\otimes_{\cA'} \cF_{\mu})(\infty).
\end{equation}

\begin{definition}
\label{definition:tenpara}
The tensor product $(\cEd\otimes_{\cA'}\cFd)_{\star}\in\PMod_{X,e}(S)$ is defined as the collection of $\cO_{X\times S}$-modules 
\[
(\cEd\otimes_{\cA'}\cFd)_{i} = \sum_{\la+ \mu =i, \la,\mu\in\bZ} \, \cE_{\la}\otimes_{\cA'} \cF_{\mu}
\]
(for $i\in\bZ$) together with the maps (\ref{eqn:tenparaincl}) and (\ref{eqn:tenparaphi}).
\end{definition}

\begin{lemma} 
\label{lemma:tenpara}
Let $\cEd\in \PVect^r_{\cA',e}(S)$ and $\cFd\in {_{\cA'}\PVect^r_e(S)}$. Then $(\cEd\otimes_{\cA'}\cFd)_{\star}$ lies in $\PVect^{rd^2}_{X,e}(S)$. In particular $\sum_{\la+ \mu =i, \la,\mu\in\bZ} \, \cE_{\la}\otimes_{\cA'} \cF_{\mu}$ is a locally free $\cO_{X\times S}$-module of rank $rd^2$ for all $i\in \bZ$.
\end{lemma}

{\em Proof.} For $\cA' = \cO_X$ this follows immediately from the fact that a parabolic $\cO_{X\times S}$-module is a parabolic vector bundle if and only if it is parabolically flat (\cite{yokogawa}, Proposition 3.1). The general case can be deduce from this special case by Morita equivalence. More precisely since the question is local we can replace $X$ by an {\'e}tale neighbourhood of $\infty$ and therefore can assume that $\cA' \cong M_d(\cO_X)$. Let $\cI$ be an invertible $\cA'$-$\cO_X$-bimodule and $\cJ$ its inverse. Since $(\cEd\otimes_{\cA'}\cFd)_{\star} = ((\cEd\otimes_{\cA'}\cI)\otimes_{\cO_X} (\cJ\otimes_{\cA'}\cFd))_{\star}$ the assertion follows from the case $\cA' = \cO_X$.
\enddemo

Let $\cA$ be another locally principal $\cO_X$-order of rank $d^2$ and assume that $e_{\infty}(\cA) =d$ and $\cA|_U$ and $\cA'|_U$ are Morita equivalent where $U = X-\{\infty\}$. There exists an increasing families of $\cA$-$\cA'$-bimodules $\{\cI_i\mid i\in \bZ\}$ and of $\cA'$-$\cA$-bimodules $\{\cJ_i\mid i\in \bZ\}$ such that $(\cI_i)|_U = (\cI_{i+1})_U=\colon\cI_U$ and $(\cJ_i)|_U = (\cJ_{i+1})_U=\colon\cJ_U$ for all $i\in \bZ$, $\cI_U$ is an invertible $\cA_U$-$\cA'_U$-bimodule with inverse $\cJ_U$ and such that $\{(\cI_i)_{\infty}\mid i\in \bZ\}$ and $\{(\cJ_i)_{\infty}\mid i\in \bZ\}$ are as in \ref{itemize:morher}. It follows from Corollary \ref{corollary:invlocfree} that $\cI_i$ and $\cJ_j$ are locally free $\cA'$-modules of rank 1. Also we have
\[
 \cA(\frac{1}{d}\infty)\otimes_{\cA}\cI_{i} = \cI_{i+1}, \quad \cJ_{i}\otimes_{\cA}\cA(\frac{1}{d}\infty) = \cJ_{i +1}
\]
for all $i\in \bZ$.

\begin{proposition} 
\label{proposition:hermaxglob}
Put $\Vect_{\cA} =  \Vect^1_{\cA}$ and $\PVect_{\cA'} = \PVect^1_{\cA',d}$. The morphisms
\begin{equation}
\label{eqn:herparab1}
\cdot\otimes_{\cA}\cId: \Vect_{\cA} \lra \PVect_{\cA'}, \quad (\cdot\otimes_{\cA'}\cJd)_{d-1}: \PVect_{\cA'}\lra \Vect_{\cA}
\end{equation}
given by $\cE \mapsto \cE\otimes_{\cA}\cId \colon = \{\cF\otimes_{\cA}\cI_i\mid i \in\bZ\}$ and $\cEd \mapsto\sum_{\la+\mu =d-1} \, \cE_{\la}\otimes_{\cA'} \cJ_{\mu}$ are mutually inverse isomorphisms of stacks. Define $\theta: \Vect_{\cA}\to\Vect_{\cA}$ and $\theta': \PVect_{\cA'} \to \PVect_{\cA'}$ by $\theta(\cE) = \cE(\frac{1}{d}\infty)$ and $\theta'(\cEd, \psi_{\star}) = (\cE[1]_{\star}, \psi[1]_{\star})$. Then the diagrams
\begin{equation}
\label{eqn:herparab2}
\begin{CD}
\Vect_{\cA} @>\otimes_{\cA}\cId >> \PVect_{\cA'}@. \hspace{2cm} @. \PVect_{\cA'}@>(\cdot\otimes_{\cA'}\cJd)_{d-1} >>\Vect_{\cA}\\
@VV \theta V @VV \theta' V @. @VV \theta' V@VV \theta V\\
\Vect_{\cA} @>\otimes_{\cA}\cId >> \PVect_{\cA'}@.\hspace{2cm} @. \PVect_{\cA'}@>(\cdot\otimes_{\cA'}\cJd)_{d-1} >>\Vect_{\cA}\end{CD}
\end{equation}
are 2-commutative.
\end{proposition}

{\em Proof.} In view of \ref{lemma:tenpara} we only have to show that the second morphism is well-defined. By \ref{lemma:tenpara} and \ref{lemma:stablyfree} we have to prove that for $\cEd\in \PVect_{\cA'}(S)$ the quotient 
\[
(\sum_{\la+\mu =0, \la,\mu\in \bZ} \, \cE_{\la}\otimes_{\cA'} \cJ_{\mu})/(\sum_{\la+\mu =-1, \la,\mu\in \bZ} \, \cE_{\la}\otimes_{\cA'} \cJ_{\mu})\cong \sum_{\la+\mu =0, \la,\mu\in \bZ} \, \barcE_{\la}\otimes_{\cA'} \cJ_{\mu}
\]
is a special $\cA$-module on $S \cong \infty \times S$ where $\barcE_{\star} \colon = \Coker(\cEd[-1] \hookrightarrow \cEd)\in _{\cA'}\PMod_{d}(S)$. However this follows from:

\begin{lemma}
\label{lemma:tensorspecial} 
The assignment $\cKd \mapsto\sum_{\la+\mu =d-1} \, \cK_{\la}\otimes_{\cA'} \cJ_{\mu}$ defines a morphism
\[
(\cdot\otimes_{\cA'}\cJd)_{d-1}: \PCoh^1_{\cA',sp,d}\lra \Coh^1_{\cA,\spe}.
\]
\end{lemma}

{\em Proof.} By Lafforgue's Lemma (\cite{lafforgue}, I.2.4) (applied to a maximal tori in $\cA$) it suffices to consider the case where $S = \Spec k$ and $k$ is an algebraically closed field. If $N(\cKd)\ne \infty$ then the assertion follows from Remarks \ref{remarks:special} (b), (c). Now assume $N(\cKd) = \infty$. Let $\cO = \cO_{\infty}$, $K= F_{\infty}$. Since the question is local with respect to the {\'e}tale topology we can replace $X$ by $\Spec \cO$ where $\cO = \cO_{\infty}$. Then $\cA'\cong M_{d}(\cO)$ and $\cA \cong \End(\cLd)$ for a lattice chain $\cLd$ of period $d$ in $K^d$. Morita equivalence allows us to replace $\cA'$ by $\cO$, i.e.\ we can assume that $\cA'= \cO$. Then $\cM_i \colon = \Gamma(\Spec \cO, \cK_i)$ is a onedimensional $k$-vector space for all $i\in \bZ$ and we have to show that 
\begin{equation}
\label{eqn:compareaell}
 \sum_{i+j = 0} \, \cM_i\otimes_{\cO} \cJ_j
\end{equation}
is a free $\barcT = \cT\otimes k$-module of rank 1 where $\cT\cong \cO^d$ is any maximal torus in $\cA$. If $1 = e_1+\ldots+ e_d$ is a decomposition of $1\in\cT$ into primitive idempotents we obtain a corresponding decomposition of (\ref{eqn:compareaell}) into 
\[
\sum_{i+j = 0} \, \cM_i\otimes_{\cO} \cJ^{(\nu)}_j, \qquad \nu=1, \ldots, d
\]
where $\cJ^{(\nu)}_j \colon = \cJ_j e_{\nu}$. Since $\cJ_j$  is free of rank 1 as a $\cT$-module, $\cJd^{(\nu)}$ is a {\it shifted parabolic line bundle} (compare \cite{yokogawa}) for each $\nu\in \{1, \ldots, d\}$. Therefore $(\cMd\otimes_{\cO}\cJd^{(\nu)})_{\star}\cong \cMd[m]$ for some $m\in \bZ$. Consequently
\[
\sum_{i+j = 0} \, \cM_i\otimes_{\cO} \cJ^{(\nu)}_j\cong \cM_m
\]
is a onedimensional $k$-vector space. It follows that (\ref{eqn:compareaell}) is a free $\cT\otimes_{\cO} R$-module of rank 1.\enddemo

Now assume that $k=\bFq$ and that $\cA'$ is a maximal $\cO_X$-order in a central division algebra $A'$ of dimension $d^2$ with $A'_{\infty}\cong M_d(F_{\infty})$. Let us recall the definition of an $\cA'$-elliptic sheaf given in (\cite{lrs}, 2.2) and (\cite{stuhler}, 4.4.1) (here we do not require $\deg(\infty)=1$). 

\begin{definition}
\label{definition:aelllrs} Let $S\in \Sch/\bFq$. An $\cA'$-elliptic sheaf $E' = (\cE_i, j_i, t_i)_{i\in\bZ}$ with pole $\infty$ in the sense of \cite{lrs} consists of a commutative diagram
\[
\begin{CD}
\ldots @>>> \cE_{i-1} @> j_{i-1} >> \cE_{i} @> j_i >> \cE_{i+1} @>>> \ldots\\
@. @AA t_{i-2} A @AA t_{i-1} A @AA t_i A \\
\ldots @>>> \tcE_{i-2} @> j_{i-1} >> \tcE_{i-1} @> j_{i-1} >> \tcE_{i} @>>> \ldots
\end{CD}
\]
where $\cE_i$ are locally free $\cO_{X\times S}$-modules of rank $d^2$ additionally equipped with a right action of $\cA'$ compatible with the $\cO_X$-action. The maps are injective $\cA'\boxtimes \cO_S$-linear homomorphisms. 

Furthermore the following conditions should hold:

\noi (i) (Periodicity) $\cE_{i + e\deg(\infty)} = \cE_i(\infty) \colon = \cE_i\otimes_{\cO_{X\times S}} (\cO(\infty)\boxtimes \cO_S)$ where the canonical embedding of $\cE_i$ on the right side corresponds on the left to the composition $\cE_i \stackrel{j}{\hra} \ldots \stackrel{j}{\hra} \cE_{i + d' deg(\infty)}$.

\noi (ii) The quotient sheaf $\cE_{i}/j_{i-1}(\cE_{i-1})$ is a locally free sheaf of rank $d$ on the graph of a morphism $\io_{\infty, i}: S \to X$.

\noi (iii) There exists a morphism $z: S \to X- |\Bad(\cA')|$ -- called the zero or characteristic of $E'$ -- such that for all $i\in \bZ$, $\Coker(t_i)$ is supported on the graph of a morphism $z$ and is a direct image of a locally free $\cO_S$-module of rank $d$ by $\Gamma_{z}= (z,\id_S): S \lra X \times S$.
\end{definition}

We first remark that condition (iii) implies that $\cE_{i}$ is actually a locally free $\cA'\boxtimes \cO_S$-module. This follows from (\cite{lafforgue}, I.4, proposition 7) or can be deduced from Lemma \ref{lemma:stablyfree} together with (\cite{lrs}, 2.6). Secondly condition (i) implies that $\io_{\infty, i}(S) = \{\infty\}$ and we have 
\[
\io_{\infty, i}\circ \Frob_S = \io_{\infty, i+1}
\]
for all $i \in \bZ$. For that consider the two filtrations of $\cE_{i+1}/t_{i-1}(\tcE_{i-1})$
\[
0 \seq \cE_{i}/t_{i-1}(\tcE_{i-1})\seq \cE_{i+1}/t_{i-1}(\tcE_{i-1}),\,\, 0 \seq t_i(\tcE_i)/t_{i-1}(\tcE_{i-1})\seq \cE_{i+1}/t_{i-1}(\tcE_{i-1}).
\]
The first shows that the support of $\cE_{i+1}/t_{i-1}(\tcE_{i-1})$ is $\Gamma_{z}+\Gamma_{\io_{\infty,i+1}}$ and the second that it is $\Gamma_{z}+\Gamma_{\io_{\infty,i}\circ \Frob_S}$.

Suppose again that $\deg(\infty) = 1$. Hence the stack $\PEll^{\infty}_{\cA'}(S)$ of $\cA'$-elliptic sheaves as defined in \ref{definition:aelllrs} is isomorphic to the stack of triples $E'=(\cEd, t_{\star})$ where $\cEd= (\cEd, \psi_{\star})\in \PVect_{\cA'}(S)$ and $t_{\star}: \tcE[-1]_{\star} \to \cEd$ is a morphism in $\PVect_{\cA'}(S)$ such that (iii) above holds for $\Coker(t_{\star})$. 

We show that the isomorphisms (\ref{eqn:herparab1}) yield isomorphisms between $\PEll_{\cA'}^{\infty}$ and $\Ell_{\cA}^{\infty}|_{X-|\Bad(\cA')|} = \Ell_{\cA}^{\infty}\times_X (X-|\Bad(\cA')|)$. Define
\[
\cdot\otimes_{\cA}\cId: \Ell_{\cA}^{\infty}|_{X-|\Bad(\cA')|}(S) \lra \PEll_{\cA'}^{\infty}(S)
\]
by $(\cE, t)\mapsto (\cE\otimes_{\cA}\cId, t\otimes_{\cA}\cId)$. The commutativity of the first diagram (\ref{eqn:herparab2}) shows that $t\otimes_{\cA}\cId$ is a map $\tcE\otimes_{\cA}\cId[-1]\to \cE\otimes_{\cA}\cId$. That $E'$ has property (iii) above follows from Remark \ref{remarks:special} (c). Conversely, we define
\[
(\cdot\,\otimes_{\cA'}\cJd)_{d-1}:\PEll_{\cA'}^{\infty}(S)\lra \Ell_{\cA}^{\infty}|_{X-|\Bad(\cA')|}(S)
\]
by $(\cEd,t_{\star}) \mapsto ((\cEd\otimes_{\cA'}\cJd)_{d-1}, (t_{\star}\otimes_{\cA'}\cJd)_{d-1})$. Again the commutativity of the second diagram of (\ref{eqn:herparab2}) implies that $(t_{\star}\otimes_{\cA'}\cJd)_{d-1}$ is a morphism $\ta(\cE(-\frac{1}{d}\ioinf)) \to \cE$ where $\cE = (\cEd\otimes_{\cA'}\cJd)_{d-1}$. Finally condition (*) of Definition \ref{definition:aell} follows Lemma \ref{lemma:tensorspecial}. We deduce from \ref{proposition:hermaxglob}:

\begin{proposition} 
\label{proposition:comparison}
Let $S$ be a $k$-scheme. The morphisms
\begin{align}
\label{eqn:comparison1}
& \cdot\otimes_{\cA}\cId: \Ell_{\cA}^{\infty}|_{X-|\Bad(\cA')|}  \lra  \PEll_{\cA'}^{\infty},\\ \label{eqn:comparison2}
& (\cdot\otimes_{\cA'}\cJd)_{d-1}: \PEll_{\cA'}^{\infty}\lra  \Ell_{\cA}^{\infty}|_{X-|\Bad(\cA')|}
\end{align}
are mutually inverse isomorphisms.
\end{proposition}

\begin{remarks}
\label{remarks:degbeliebig} 
\rm (a) In order to extend \ref{proposition:comparison} to the case $\deg(\infty) >1$ we have to modify Definition \ref{definition:parvecact} (b) as follows. For $S\in \Sch/k$ let $\PVect_{\cA'}(S)$ denote the category of triples $(\cFd, \psi_{\cFd},\ioinf)$ where $\ioinf: S \to X$ is a $k$-morphism which factors through $\infty\to X$ and $(\cFd, \psi_{\cFd})$ is an element of $\PMod_{\cA',d\deg(\infty)}(S)$ such that $\cF_i\in \Vect^1_{\cA'}(S)$ for all $i\in \bZ$ and such that the sheaf $\Coker(j_i)$ is a locally free sheaf of rank $d$ on the graph of $\ioinf\circ \Frob_S^i: S \to X$. To define isomorphisms similar to (\ref{eqn:herparab1}) we consider increasing families of $\cA\boxtimes k(\infty)$-$\cA'\boxtimes k(\infty)$-bimodules $\{\cI_i\mid i\in \bZ\}$ and  $\cA'\boxtimes k(\infty)$-$\cA\boxtimes k(\infty)$-bimodules $\{\cJ_i\mid i\in \bZ\}$ with the following properties: 
\begin{enumerate}
\item[(i)]
 $\cA(\frac{1}{d}\infty_i)\otimes_{\cA}\cI_{i} = \cI_{i+1}, \quad \cJ_{i}\otimes_{\cA}\cA(\frac{1}{d}\infty_i) = \cJ_{i +1}$ for all $i\in \bZ$. Here $\infty_0$ denotes the canonical morphism $\Spec k(\infty)\to X$ and $\infty_i\colon =\infty_0\circ \Frob^i:\Spec k(\infty) \to X$.
\item[(i)] $\cI_U= (\cI_i)|_{U\times_{\bFq} k(\infty)}$ is an invertible $\cA_U\boxtimes k(\infty)$-$\cA'_U\boxtimes k(\infty)$-bimodule with inverse $\cJ_U= (\cJ_i)|_{U\times_{\bFq} k(\infty)}$.
\item[(iii)] For all $i\in \bZ$, $\cI_i$ and $\cJ_j$ are locally free $\cA'\boxtimes k(\infty)$-modules of rank 1. 
\end{enumerate}
As in \ref{proposition:hermaxglob} one defines isomorphisms
\begin{align*}
& \cdot\,\otimes_{\cA}\cId:\Vect_{\cA}\times_{\bFq} k(\infty) \lra \PVect_{\cA'},\\
& (\cdot\,\otimes_{\cA'}\cJd)_{d-1}: \PVect_{\cA'}\lra \Vect_{\cA}\times_{\bFq} k(\infty)
\end{align*}
which then yield the isomorphisms (\ref{eqn:comparison1}), (\ref{eqn:comparison2}) above.

\noi (b) Let $\fp$ be a closed point of $X$ such that $\inv_{\fp}(A') = \frac{1}{d}$. In \cite{hausberger}, Hausberger constructed a flat proper model of $\Ell^{\infty}_{\cA'}$ over $(X-|\Bad(\cA')|)\cup \{\fp\}$ by extending the definition of the moduli problem \ref{definition:aelllrs} of Laumon-Rapoport-Stuhler to characteristic $\fp$. By using (\cite{hausberger}, 2.16) it is easy to see using that his {\it condition sp{\'e}ciale} (\cite{hausberger}, section 3) corresponds to our condition (*).  
\end{remarks}

\bigskip
\begin{tabbing} 
 \hspace{1cm}\= Michael Spie{\ss}\\
 \>Fakult{\"a}t f{\"u}r Mathematik\\
 \>Universit{\"a}t Bielefeld\\
\>Postfach 100131\\
\>D-33501 Bielefeld\\
\>Germany\\
\end{tabbing}

\noindent \hspace{.9cm} e-mail: mspiess@mathematik.uni-bielefeld.de
\end{document}